\documentclass{article}
\usepackage{amsmath,amsthm,amsfonts,amssymb,amscd,latexsym, fullpage,textcomp,amscd}
\date{ }
\begin{document}

\theoremstyle{plain}
\newtheorem{thm}{\bf Theorem}[section]
\newtheorem{lem}[thm]{\bf Lemma}
\newtheorem{cor}[thm]{\bf Corollary}
\newtheorem{prop}[thm]{\bf Proposition}
\theoremstyle{remark}
\newtheorem{Case}{\bf Case}
\newtheorem{rem}[thm]{\bf Remark}
\newtheorem{claim}[thm]{\bf Claim}

\theoremstyle{definition}
\newtheorem{Def}{\bf Definition}
\newtheorem*{pf}{\bf Proof}
\newtheorem{Conj}{\bf Conjecture}

\newcommand{\nc}{\newcommand}

\newcommand{\BB}{{\mathcal B}} 
\newcommand{\CC}{{\mathcal C}}
\newcommand{\DD}{{\mathcal D}}
\newcommand{\EE}{{\mathcal E}}
\newcommand{\FF}{{\mathcal F}}
\newcommand{\GG}{{\mathcal G}}
\newcommand{\HH}{{\mathcal H}}
\newcommand{\II}{{\mathcal I}}
\newcommand{\JJ}{{\mathcal J}}
\newcommand{\KK}{{\mathcal K}}
\newcommand{\LL}{{\mathcal L}}
\newcommand{\MM}{{\mathcal M}}
\newcommand{\NN}{{\mathcal N}}
\newcommand{\OO}{{\mathcal O}}
\newcommand{\PP}{{\mathcal P}}
\newcommand{\QQ}{{\mathcal Q}}
\newcommand{\RR}{{\mathcal R}}
\newcommand{\TT}{{\mathcal T}}
\newcommand{\UU}{{\mathcal U}}
\newcommand{\VV}{{\mathcal V}}
\newcommand{\WW}{{\mathcal W}}
\newcommand{\ZZ}{{\mathcal Z}}
\newcommand{\XX}{{\mathcal X}}
\newcommand{\YY}{{\mathcal Y}}
\nc{\bba}{{\mathbb A}}
\nc{\bbb}{{\mathbb B}}
\nc{\bbc}{{\mathbb C}}
\nc{\bbd}{{\mathbb D}}
\nc{\bbe}{{\mathbb E}}
\nc{\bbf}{{\mathbb F}}
\nc{\bbg}{{\mathbb G}}
\nc{\bbh}{{\mathbb H}}
\nc{\bbi}{{\mathbb I}}
\nc{\bbj}{{\mathbb J}}
\nc{\bbk}{{\mathbb K}}
\nc{\bbl}{{\mathbb L}}
\nc{\bbm}{{\mathbb M}}
\nc{\bbo}{{\mathbb O}}
\nc{\bbp}{{\mathbb P}}
\nc{\bbq}{{\mathbb Q}}
\nc{\bbr}{{\mathbb R}}
\nc{\bbs}{{\mathbb S}}
\nc{\bb}{{\mathbb T}}
\nc{\bbu}{{\mathbb U}}
\nc{\bbv}{{\mathbb V}}
\nc{\bbw}{{\mathbb W}}
\nc{\bbx}{{\mathbb X}}
\nc{\bby}{{\mathbb Y}}
\nc{\bbz}{{\mathbb Z}}
\nc{\fA}{{\mathfrak A}}
\nc{\fB}{{\mathfrak B}}
\nc{\fC}{{\mathfrak C}}
\nc{\fD}{{\mathfrak D}}
\nc{\fE}{{\mathfrak E}}
\nc{\fF}{{\mathfrak F}}
\nc{\fG}{{\mathfrak G}}
\nc{\fH}{{\mathfrak H}}
\nc{\fI}{{\mathfrak I}}
\nc{\fJ}{{\mathfrak J}}
\nc{\fK}{{\mathfrak K}}
\nc{\fL}{{\mathfrak L}}
\nc{\fM}{{\mathfrak M}}
\nc{\fN}{{\mathfrak N}}
\nc{\fO}{{\mathfrak O}}
\nc{\fP}{{\mathfrak P}}
\nc{\fQ}{{\mathfrak Q}}
\nc{\fR}{{\mathfrak R}}
\nc{\fS}{{\mathfrak S}}
\nc{\fT}{{\mathfrak T}}
\nc{\fU}{{\mathfrak U}}
\nc{\fV}{{\mathfrak V}}
\nc{\fW}{{\mathfrak W}}
\nc{\fZ}{{\mathfrak Z}}
\nc{\fX}{{\mathfrak X}}
\nc{\fY}{{\mathfrak Y}}
\nc{\fa}{{\mathfrak a}}
\nc{\fb}{{\mathfrak b}}
\nc{\fc}{{\mathfrak c}}
\nc{\fd}{{\mathfrak d}}
\nc{\fe}{{\mathfrak e}}
\nc{\ff}{{\mathfrak f}}
\nc{\fh}{{\mathfrak h}}
\nc{\fj}{{\mathfrak j}}
\nc{\fk}{{\mathfrak k}}
\nc{\fl}{{\mathfrak{l}}}
\nc{\fm}{{\mathfrak m}}
\nc{\fn}{{\mathfrak n}}
\nc{\fo}{{\mathfrak o}}
\nc{\fp}{{\mathfrak p}}
\nc{\fq}{{\mathfrak q}}
\nc{\fr}{{\mathfrak r}}
\nc{\fs}{{\mathfrak s}}
\nc{\ft}{{\mathfrak t}}
\nc{\fu}{{\mathfrak u}}
\nc{\fv}{{\mathfrak v}}
\nc{\fw}{{\mathfrak w}}
\nc{\fz}{{\mathfrak z}}
\nc{\fx}{{\mathfrak x}}
\nc{\fy}{{\mathfrak y}}

\nc{\al}{{\alpha }}
\nc{\be}{{\beta }}
\nc{\ga}{{\gamma }}
\nc{\de}{{\delta }}
\nc{\vap}{{\tepsilon }}

\nc{\ze}{{\zeta }}
\nc{\et}{{\eta }}
\nc{\vth}{{\vartheta }}

\nc{\io}{{\iota }}
\nc{\ka}{{\kappa }}
\nc{\la}{{\lambda }}
\nc{\vpi}{{     \varpi          }}
\nc{\vrho}{{    \varrho         }}
\nc{\si}{{      \sigma          }}
\nc{\ups}{{     \upsilon        }}
\nc{\vphi}{{    \varphi         }}
\nc{\om}{{      \omega          }}

\nc{\Ga}{{\Gamma }}
\nc{\De}{{\Delta }}
\nc{\nab}{{\nabla}}
\nc{\Th}{{\Theta }}
\nc{\La}{{\Lambda }}
\nc{\Si}{{\Sigma }}
\nc{\Ups}{{\Upsilon }}
\nc{\Om}{{\Omega }}

\nc{\zz}{{\mathbb Z}}
\newcommand{\N}{{\mathbb N}}
\newcommand{\etat}{\tilde{\eta}}
\newcommand{\sln}{\mathfrak{sl} _n }
\newcommand{\slnr}{{\mathfrak {sl} _n (\mathbb R)}}
\newcommand{\sun}{\mathfrak{su} _n}
\newcommand{\cc}{{\mathbb C}}
\newcommand{\rr}{{\mathbb R}}
\newcommand{\ac}{{\check{\alpha}}}
\newcommand{\orb}{{\mathcal O}}
\newcommand{\gcs}{{\mathcal{GC} _{G/K}}}
\newcommand{\cds}{{\mathbb{C} \mathcal{D} _{G/K}}}
\newcommand{\ds}{{\mathcal{D} _{G/K}}}
\newcommand{\ggcs}{{\mathcal{GC} _{G/K} ^G}}
\newcommand{\gcds}{{\mathbb{C} \mathcal{D} _{G/K} ^G} }
\newcommand{\gds}{{\mathcal{D} _{G/K} ^G}}
\newcommand{\bc}{{\check{\beta}}}
\newcommand{\tep}{{\varepsilon_\sharp}}
\newcommand{\ep}{{\varepsilon_\sharp}}
\newcommand{\vep}{\varepsilon}
\newcommand{\epp}{T_{\overline{\varepsilon}}}
\newcommand{\ctsd}{\mathfrak{g} _\mathbb{C} ^*}
\newcommand{\gts}{$\mathfrak{g}\oplus\mathfrak{g}^{\ast}$}
\newcommand{\inv}{^{-1}}
\newcommand{\fg}{{\mathfrak g}}
\newcommand{\noi}{\noindent}
\newcommand{\lra}{{\longrightarrow}}
\newcommand{\pbd}{f^\star \mathcal{D}}
\newcommand{\cep}{{\overline{\varepsilon}}}
\newcommand{\cdt}{\tilde{\mathcal{D}}}
\newcommand{\cts}{{\mathfrak{g} _{\mathbb{C}}}}
\newcommand{\csts}{\mathfrak{k} _{\mathbb{C}}}
\newcommand{\cgts}{{\mathfrak{g} _{\mathbb C} \oplus \mathfrak{g} _\mathbb{C} ^* } }

\newcommand{\pionep}{{ (\pi_ 1)_\sharp}}
\newcommand{\pitwop}{{(\pi _2 )_\sharp}}
\newcommand{\omonep}{{(\omega _1)_\sharp}}
\newcommand{\omtwop}{{(\omega_2)_\sharp}}
\nc{\st}{{\; | \; }}

\nc{\trm}{\textreferencemark}
\nc{\tih}{{\tilde \HH}}
\nc{\vb}{\text{vector bundle}}
\nc{\vbm}{{\VV ect \BB und _M}}
\nc{\vs}{{Vect_{fd}}}
\nc{\cxn}{\text{connection}}
\nc{\bla}{{\mathfrak g ^* \rtimes \mathfrak g}}
\nc{\cbla}{{\mathfrak g _\cc ^* \rtimes \mathfrak g _\cc}}
\nc{\gc}{\text{generalized complex }}
\nc{\gcstr}{\text{generalized complex structure }}
\nc{\gcstrs}{\text{generalized complex structures }}

\title{Generalized Complex and Dirac Structures on Homogeneous Spaces}
\author{Brett Milburn}
\maketitle

\begin{abstract}
\noi Generalized complex geometry \cite{hit} has been a subject of recent interest in mathematics and physics and is a general setting for differential geometry.  The aim of this paper is to study generalized complex geometry and Dirac geometry~\cite{cou}, \cite{cow} on homogeneous spaces.  We offer a characterization of equivariant Dirac structures on homogeneous spaces, which is then used to construct new examples of generalized complex structures.  We consider Riemannian symmetric spaces, quotients of compact groups by closed connected subgroups of maximal rank, and nilpotent orbits in $\fs \fl _n (\mathbb R)$.  For each of these cases, we completely classify equivariant Dirac structures.  Additionally, we consider equivariant Dirac structures on semisimple orbits in a semisimple Lie algebra.  Here equivariant Dirac structures can be described in terms of root systems or by certain data involving parabolic subagebras.  \\
\end{abstract}

\tableofcontents

\section{Introduction} 

\noi Hitchin's notion of ``generalized complex geometry'' \cite{hit} is a common generalization of complex and symplectic structures which has found several applications in physics and mathematics~\cite{kap},~\cite{hit},~\cite{gua}.  For instance, the generalized K\"ahler structure turns out to be precisely the setting for $N=2$ supersymmetric sigma models~\cite{kap}.  Generalized complex structures are a special case of complex Dirac structures, a concept defined by Courant and Weinstein~\cite{cou},~\cite{cow} that also includes Poisson structures, integrable distributions, and presymplectic structures. \\


\noi After a brief exposition on the essential ideas in generalized complex geometry in \S \ref{gcgeometry}, we partially classify equivariant generalized complex structures on 
homogeneous spaces and, more generally, equivariant (complex) Dirac structures on homogeneous spaces in  \S \ref{GChomo} and \S \ref{GCclassification}.  This gives a description of such Dirac structures in terms of linear algebra data.  The $G$-equivariant generalized complex structures $L$ on $G/K$ are in bijection with pairs $(E,\vep)$ of a subalgebra 
$E\subset \fg _\cc$ and $\vep \in \wedge ^2 E^*$ satisfying certain conditions. 
 This bijection allows us to provide some new examples of generalized complex structures.  \\

\noi Here is a list of results for particular classes of homogeneous spaces.  \\

\noi \S \ref{26nov9} \emph{$G$ compact and $K$ connected of maximal rank} (i.e. $K$ contains a Cartan subgroup).  We completely classify equivariant generalized complex structures.  In this setting, there are examples of generalized complex structures which are neither symplectic nor complex. \\

\noi  \S \ref{sectionorbits}\emph{Semisimple coadjoint orbits in semisimple real Lie algebras}.  We describe generalized complex structures in terms of simpler ``combinatorial'' data that involves only the root system.\\

\noi \S \ref{sectionnilorbits} \emph{Real nilpotent ortbits}.  Here we restrict ourselves to split semisimple Lie algebras.  For the Lie algebra $\mathfrak {sl} _n(\rr)$, the only equivariant generalized complex structures are B-transforms of symplectic structures.  
We hypothesize that this is true for any split semisimple Lie algebra and show that the claim reduces
to distinguished orbits in simple Lie algebras.  \\  

\noi  \S \ref{sectionriemann}\emph{Riemannian Symmetric Spaces}.  Again we completely classify equivariant generalized complex structures but  
these turn out to yield little that is new.  Every generalized complex structure on a Riemannian symmetric space is 
essentially a product  of complex and (B-transforms of) symplectic structures. \\


\section{Dirac and Generalized Complex Geometry}\label{gcgeometry} 

\noi We introduce the basic definitions and notational conventions used in
this paper.  For a systematic development of generalized complex
structures as well as some of their applications, we refer the reader to ~\cite{gua}.  The
notation of ~\cite{gua} will most strongly be followed.\\

\noi For a manifold $M$, generalized geometry is concerned with the bundle
$\mathcal V _M := TM \oplus
T^*M$.  There is a natural bilinear form on $\mathcal V _M$, given by
the obvious pairing $\langle X + \xi , Y + \eta \rangle = X(\eta) + Y (\xi)$ for sections $X,Y$ of $TM$ and $\xi $, $\eta$ of $T^*M$.  Furthermore,
$\mathcal V _M$ is equipped with the \emph{Courant bracket} defined by 
\[ [X + \xi , Y + \eta ] = [X,Y] + \io _X d \eta +\frac{1}{2} d( \io _X \eta) - \io _Y d \xi - \frac{1}{2} d(\io _Y \xi),\]
where $\io$ denotes contraction in the first variable ($ \io _x \phi = \phi (x,-,...)$).  
The Courant bracket $[ \; , \; ]$ and the bilinear form $\langle \; , \; \rangle$
extend $\cc$-bilinearly to $(\mathcal V _M )_\cc = \mathcal V _M
\otimes \cc$.  


\begin{Def} A \emph{generalized almost complex structure} on $M$ is a map $\JJ : \VV _M \lra \VV _M$ such that $\JJ$ is orthogonal with respect to the inner product $\langle \; , \; \rangle$ and $\JJ ^2 = -1$.  Just as with complex structures, one may consider the $i$-eigenbundle, $D$, of $\JJ $ in $(\VV _M)_\cc$.  The Courant bracket defines an integrability condition ($[D , D] \subset D$) for $\JJ$ to be called a \emph{ generalized complex structure}.  This follows the analogy with almost complex stuctures; an almost complex structure is a complex structure precisely when its $i$-eigenbundle is integrable with respect to the Lie bracket. 
\end{Def} 

\noi The two canonical examples of generalized complex structures come from complex and symplectic structures.  Since $\VV _M = TM \oplus T^*M$, we can express any map $\VV _M \lra \VV _M$ as a block matrix in terms of this decomposition, and we will follow this convention throughout the text.  If $J$ is a complex structure, 

\[
\begin{bmatrix}
J & 0 \\
0 & -J^*
\end{bmatrix} 
\]

\noi is a generalized complex structure.  The $i $-eigenbundle of $\JJ$ is $E \oplus Ann(E)$, where $E$ is the $i$-eigenbundle of $J$, and $Ann(E)$ is the annihilator of $E$ in $T^*M$.  \\

\noi For a symplectic structure $\om$, we get a generalized complex structure 
\[
\begin{bmatrix}
0 & -\om _\sharp \inv \\
\om _\sharp & 0
\end{bmatrix} 
\]
\noi where $\om _\sharp (x) := \om (x , -)$.  The $i$-eigenbundle is the graph of $i\om _\sharp $ in $(\VV _M)_\cc$.  The fact that the symplectic form $\om$ is closed implies that this generalized almost complex structure is integrable, hence a generalized complex structure.   \\

\noi The $i$-eigenbundle $D$ of a generalized complex structure $\JJ$ turns out to be an integrable maximal isotropic subbundle of $(\VV _M)_\cc$, also known as a \emph{complex Dirac structure}.  
Thus, the study of generalized geometry now lies in the framework of Dirac structures.   With this in mind, we recall the following working definitions for our paper.  

\begin{Def} For any manifold, $M$, 
\begin{enumerate}
\item  A \emph{real almost Dirac structure} on $M$ is a maximal isotropic
  subbundle $D$ of  $\mathcal V _M$. A real almost Dirac structure is called a \emph{real Dirac structure} if it is
  integrable with respect to the Courant bracket.  Similarly, a \emph{complex almost Dirac structure} is a maximal isotropic
  subbundle $D \subset (\mathcal V _M )_\cc $, and a \emph{complex Dirac
  structure} is an integrable complex almost Dirac structure.   
\item A complex Dirac structure $D$ is said to be of \emph{constant rank} if the projection map $pr : D \lra TM$ is of constant rank.  
\end{enumerate}
\end{Def}

\noi A generalized (almost) complex structure $\JJ$ is equivalent to is a complex (almost) Dirac structure $D$
  such that $D \cap \overline D = 0$.  Note that the integrability is a closed condition and that being generalized complex is an open condition.  Henceforth we will think of generalized complex structures as complex Dirac structures.  \\

\noi Since the complexification of any Dirac structure is a complex Dirac
structure, both generalized complex structures and Dirac structures
are complex Dirac structures.  Thus, the set of complex Dirac structures contains real Dirac structures and generalized complex structures.  Henceforth (almost) \emph{Dirac structure} will always mean complex (almost) Dirac structure, and we will specify whether it is also real Dirac (i.e. if $\overline  {\mathcal D } = \mathcal D$) if there is any ambiguity. \\

\noi Most of the complex Dirac structures
considered in this paper will be of constant rank. 
It is checked in ~\cite{gua} that any complex Dirac structure of constant rank is of the form 
\[ L(E , \vep ) := \{ X + \xi \in (\mathcal V _M)_\cc \; | \; X \in E \; and \;  \iota _X \vep = \xi _{|E} \}, 
\]
\noi where $E$ is a subbundle of $TM$ and $\vep \in \Ga (M , \wedge ^2 E^*)$.  Complex and symplectic structures, for example, are of this form.  For a subbundle $E$ of $TM$, we define the differential \\
$d_E : \Ga (M, \wedge ^2 E^*) \lra \Ga (M, \wedge ^3 E^*)$ by the following forumula.  For sections $X,Y,Z$ of $E$ and  $\vep \in \Ga (M , \wedge ^2 E^*)$,  
\[d_E \vep (X,Y,Z) = \vep (X,[Y,Z]) + \vep (Y,[Z,X]) + \vep (Z,[X,Y]) + X\vep (Y,Z) - Y\vep (X,Z) + Z\vep (X,Y). \]
\noi In other words, $d_E \vep $ is the restriction to $\wedge ^3 E $ of the ordinary De Rham differential of any extension $\tilde \vep \in \wedge ^2 T^*M $ of $\vep$.  Gualtieri \cite{gua} proves the following useful lemma.

\begin{lem}\label{21sept1} A complex almost  Dirac structure of constant rank $L(E ,
  \vep )$ is integrable if and only if $E$ is integrable and $d_E \vep
  =0$.  
\end{lem}  


\begin{rem} There are several formulations of the integrability condition for an almost Dirac structure $L$ on $M$.  For example, the Lie algebroid derivative $d$ defined in~\cite{gua} extends to $d : \wedge ^k \VV _M ^*\lra \wedge ^{k+1} \VV _M ^*$.  Then it is easily verified that since $L = Ann(L)$ and $\VV _M $ is naturally isomorphic to $\VV _M ^*$ via the inner product, $L$ is integrable if and only if $dL \subset L \wedge \VV _M$.  Gualtieri~\cite{gua} shows that there is a similar condition when one describes $L$ in terms of pure spinors.  Since equivariant Dirac structures on a homogeneous space are always of constant rank, lemma~\ref{21sept1} will the the most useful criterion for determining integrability.  The following lemma also allows one to determine integrability in special cases by more familiar criterion--such as having a Poisson bivector.  
\end{rem}

\begin{lem}We recall from~\cite{gua} that the following procedures create Dirac structures 
from geometric structures on $M$.  The first three are real, and (4) and (5) are special cases of generalized complex structures (as we have seen).
\begin{enumerate}
\item 
To an integrable distribution $\mathcal{D}\subset TM$, assign
$[\mathcal{D}\oplus Ann(\mathcal{D})]_\cc$.
\item 
To a Poisson structure $\pi \in \Ga (M, \wedge ^2 TM)$, assign $L(\pi , T^*M)$.
\item 
To a presymplectic structure $\om \in \Om ^2 (M)$, assign $L(TM , \om)$.

\item 
To a complex structure
$J$, assign $T^{(1,0)}M \oplus T^{* ,(0,1)}M$, where $T^{(1,0)} M$ and $T^{(0,1)}M$ denote the holomorphic and antiholomorphic tangent bundles respectively with respect to $J$. 
\item To a symplectic structure $\om \in \Om ^2 (M)$, assign $L(TM_\cc , i \om)$.  
\end{enumerate}
\end{lem}

\begin{rem} A symplectic form $\omega$ on $M$  determines a complex Dirac structure
in one of two ways: $L(TM, \omega) $ and $L(T_\cc M, -i\omega)$.  The former
is a real Dirac structure, and the latter is a generalized complex
structure. 
\end{rem}

\noi This way of representing a Dirac structure $L$ as $L(E,\vep)$ turns out to be extremely useful for our purposes.  
For any vector bundle $E$ and $\vep \in \bigwedge ^2 E^*$, the convention used in this paper is for $\varepsilon _\sharp$ to denote the map $E \lra E^*$ determined by $\vep$.  That is, for $X,Y \in E$, $(\vep _\sharp X)(Y) = \vep (X,Y)$.  \\

\noi For a Dirac structure $D$, if the projection $pr : D \lra T^*M$ has constant rank, then there is some subbundle $U \subset T^*M$ and some $\pi \in \Ga (M, \bigwedge ^2 U^*)$ such that $D$ is of the form 
\[ L(\pi , U) := \{ X + \xi \st X_{|U} = \io _\xi \pi \}.  \] 
If $U = T^*M$, then for $\pi \in \Ga (M , \wedge ^2 TM)$, $L(\pi , T^*M)$ is a Dirac structure if and only if $\pi$ is a Poisson bi-vector ~\cite{gua}, ~\cite{vai}.  Now presymplectic structures, complex structures, and Poisson structures can all be considered Dirac structures.  \\

\noi We recall the notions of pullback and pushforward of linear Dirac structures ~\cite{gua}. For a map $F: V \lra W$ of vector spaces and a subspace $D \subset V \oplus V^*$, define
\[F_\star D = \{ FX + \xi \in W \oplus W^* \st X + F^* \xi \in D \} \]
and for a subspace $D \subset W \oplus W^*$, define 
\[ F^\star D = \{ X + F^* \xi \in V \oplus V \st FX + \xi \in D \} = (F^*)_\star D. \]

\noi Now let $f : M \lra N$ be any map of
manifolds.  For a Dirac structure, $D$, on $N$, the pullback $f^\star D$  is
defined pointwise by $(f^\star D)_p = (df_p)^\star D_{f(p)} $. It is
not necessarily itself a Dirac structure.  \\

\subsection{Twisted Courant Bracket and Automorphisms}

\noi In addition to the standard Courant bracket on $\VV _M$, \v Severa and Weinstein noticed a twisted Courant bracket $[ \; , \; ] _H$ for each closed 3-form $H$ on $M$, defined as
\[ [X + \xi , Y + \eta  ] _H = [X + \xi , Y + \eta  ] + H(X,Y,-) 
\]
and $\VV _{M,H} := (TM \oplus T^*M , [ \; , \; ] _H , \langle \; , \; \rangle )$ so that $\VV _M$ is just $\VV _{M,0}$.  For any 2-form $B$ on $M$, there is an automorphism of the vector bundle $TM \oplus T^*M$, 

\[
\begin{bmatrix}
1 & 0 \\
B_\sharp & 1
\end{bmatrix}, 
\]
denoted $e^B$.  
Indeed, $e^B$ is an isomorphism $\VV _{M,H+dB} \lra \VV _{M, H}$.  In other words, $e^B$ is orthogonal with respect to $ \langle \; , \; \rangle $, and $[e^B u , e^B v]_H = e^B [u,v]_{H +dB}$ for all sections $u,v$ of $TM \oplus T^*M$.  When $B$ is closed, then $e^B$ is an automorphism of $\VV _{M,H}$.  This is what we call a \emph{B-transform}, i.e. an automorphism of $\VV _{M,H}$ of the form $e^B$ for a closed 2-form $B$.   In fact, the automorphism group of $\VV _{M,H}$ is the semidirect product of the group of diffeomorphisms $M \lra M$ and closed 2-forms $Z^2(M)$ \cite{gua}.  B-transforms are thought of as the symmetries of the Courant bracket.  An \emph{H-twisted Dirac structure} $D \subset \VV _{M,H}$ is simply a maximal isotropic subundle which is integrable with respect to the $H$-twisted Courant bracket.

\section{GENERALIZED COMPLEX AND DIRAC STRUCTURES ON HOMOGENEOUS SPACES}\label{GChomo}

\noi We partially describe equivariant Dirac and generalized complex structures on a homogeneous space $G/K$ by giving equivalent data involving only the Lie algebra.  The main results are in Theorem~\ref{dhsthm} and Proposition~\ref{b}, where we parameterize equivariant Dirac structures and generalized complex structures by pairs $(E,\vep)$ of a Lie subalgebra $E$ of $\cts$ and $\vep \in \wedge ^2 E^*$ satisfying some conditions.



\subsection{Distributions on Homogeneous Spaces}\label{dists}
\noi In this subsection, we recall some well-known background facts for operations on Dirac structures. We would like to classify G-invariant distributions on a homogeneous space G/K.  By pulling back a distribution on G/K to a distribution on G and then considering the subspace of $\mathfrak{g}$ determined by this distribution on $G$, this will give for each distribution on $G/K$ a subspace of $\mathfrak{g}$ which uniquely determines the distribution on $G/K$.  In particular, this can be used for complex distributions given by complex structures on G/K.\\

\noi Let $f : M \longrightarrow N$ be a submersion and $D \subset TN$ a distribution on N.  Define $f^{-1} D$ by $(f^{-1} D) _p = df_p ^{-1}( D_{f(p)})$.

\begin{prop}\label{distsprop1} If D is a distribution on N and $f: M
  \longrightarrow N $ is a submersion, then 
\begin{enumerate}
\item $f\inv D$ is a distribution on $M$, and it is integrable if and only if $D$ is integrable. 
\item If $M$ and $N$ are G-spaces and $f$ is equivariant, then $D$ is G-invariant if and only if $f \inv D$ is G-invariant.
\end{enumerate}  
\end{prop}

\begin{pf} 
\begin{enumerate}
\item Both of these statements are local.  Since $f$ is a submersion, locally $f$ looks like $p = pr_1 : U \times V \longrightarrow U$.  Now $p^{-1} D = D \oplus TV$, i.e. $p ^{-1} D = pr_1^* D \oplus pr_2 ^* TV$, and both claims are clear.
\item Let $L_g ^M$ denote left multiplication by $g \in G$ on $M$, and let $L_g ^N$ denote left multiplication by g on $N$. Let $p \in M$.  Since $f \circ L_g ^M = L_g ^N \circ f$, $df_{g.p} \, (dL ^M _g)_p \, (f\inv D)_p  = (dL^N _g)_{f(p)} \,  df_p (f\inv D _p) = (dL^N _g)_{f(p)} \, D_{f(p)}$.  If $f\inv D$ is left-invariant, then $ D_{g.f(p)} = df_{g.p} (f\inv D)_{g.p}= df_{g.p} (dL^M _g)_p (f\inv D)_p  =  (dL^N _g)_{f(p)} D_{f(p)}$, and $D$ is equivariant.  Now suppose that $D$ is equivariant.  Then $df_{g.p} (dL^M _g)_p (f\inv D)_p  = (dL^N_g)_{f(p)} D_{f(p)} = D_{g.f(p)}$, which shows that $(dL_g)_p (f\inv D)_p  \subset f \inv D _{g.p}$, so $(dL_g)_p (f\inv D)_p  = f \inv D _{g.p}$, and $f\inv D$ is left-invariant.  $\square$
\end{enumerate}
\end{pf}

\begin{cor}\label{d} Let $G$ be any Lie group and $K$ any closed, connected subgroup.  There is a bijection between left-invariant integrable distributions on $G$ and Lie subalgebras of $\fg$, and there is a bijection between G-invariant integrable distributions on G/K and Lie subalgebras of $\mathfrak{g}$ containing $\mathfrak{k}$.
\end{cor}

\begin{pf}   
There is a bijection between subspaces $\fu$ of $\fg$ and left-invariant distributions $\DD$ on $G$ given by $\fu = \DD _e$.  The Lie algebra structure of $\fg$ is defined by the Lie bracket on left-invariant sections of $TG$.  Therefore, the integrable left-invariant distributions correspond to subalgebras.  Now by Proposition~\ref{distsprop1} and the first part of this proof, sending a distribution $D$ on $G/K$ to $(\pi \inv D) _e$ gives an injective map from G-invariant integrable distributions on $G/K$ into subalgebras of $\fg$ containing $\fk$.  It remains to show that any subalgebra $\fu$ of $\fg$ containing $\fk$ is of the form $(\pi \inv D )_e$ for some integrable $G$-invariant distribution $D$ on $G/K$.  Given a subalgebra $\fu \supset \fk$, define $D_{\pi (g)} = (dl_g)_{\pi (e)} \fu /\fk$.  Since $\fu / \fk$ is a subrepresentation of $K$ in $\fg / \fk$, $D$ is a well-defined G-invariant distribution on $G/K$.  By Proposition~\ref{distsprop1}, $\pi \inv D$ is a left-invariant distribution on $G$.  Since $(\pi \inv D)_e = \fu$ is a subalgebra, $\pi \inv D$ is integrable and therefore so is $D$.  $\square$
\end{pf}


\begin{rem} If $K$ were disconnected, Corollary \ref{d} would say that there is a bijection between $G$-invariant integrable distributions on $G/K$ and $K$-invariant subalgebras of $\fg$ containing $\fk$.  
\end{rem}

\subsection{Homogeneous Dirac Structures}


\begin{lem}\label{pdsprop2} Let $f: M\longrightarrow N$ be a submersion and $\mathcal
  D $ be an almost Dirac structure on N.
\begin{enumerate}
\item If $\DD = L(E,\vep)$ is of constant rank, then $\pbd = L(f\inv E, f^* \varepsilon) $, so $\pbd$ is an
  almost Dirac structure.
\item $\pbd$ is integrable if and only if $\mathcal{D}$ is integrable.
\end{enumerate}
\end{lem}

\begin{pf}(1) We need to check that for all $ p \in M$,
$ L(f\inv E, f^\ast \varepsilon)_p  = (d f_p )^\star L(E, \varepsilon)_{f(p)}$. 
 Hence we must show that for a linear map of vector spaces $\vphi : V \lra W$ and $\DD = L(E,\vep) \in W\oplus W^*$, $\vphi^\star L(E,\vep) = L(\vphi \inv E, \vphi^*\vep)$.
  First we show that $\vphi^\star \DD $ is maximal isotropic, as is pointed out in~\cite{gua}.  The fact that $\langle X + \vphi^* \xi , Y + \vphi^* \eta\rangle = \langle \vphi X + \xi , \vphi Y + \eta \rangle$ implies that $\vphi^\star \DD \subset (\vphi^\star \DD)^\perp$.
  It is obvious that $\vphi ^*Ann(E) \subset Ann(\vphi \inv E)$.
  Since there is an injection $\vphi : V/\vphi \inv E \hookrightarrow W/E$,
 there is a surjection $\vphi ^* : Ann(E) \twoheadrightarrow Ann(\vphi \inv E)$
, so $Ann(\vphi \inv E) = \vphi ^* Ann(E)$.  Observe that $pr_V \vphi ^\star \DD = \vphi \inv E$
.  This combined with the fact that $Ann(\vphi \inv E) = \vphi ^* Ann(E) \subset \vphi ^\star \DD$
 implies that $dim(\vphi ^\star \DD) \geq dim(V)$.  But because
 $\vphi ^\star \DD \subset (\vphi ^\star \DD)^\perp$, one has $\vphi ^\star \DD = (\vphi ^\star \DD)^\perp
$, and this is a subspace of dimension $dim(V)$.  Therefore, $\vphi ^\star \DD$ is maximal isotropic.  \\

\noi Since $f^\star \DD$ is maximal isotropic, to show that $f^\star \DD = L(f\inv E , f^* \vep) $, it suffices that $f^\star \DD \subset L(f\inv E , \vep) $.  Starting with any $X + f^* \xi \in f ^\star \DD$, one has $fX + \xi \in \DD= L(E , \vep)$, so $X \in f\inv E$, and $\xi _{|E} = \io _{fX} \vep$.  This implies that $(f^* \xi )_{|f\inv E} = f^*(\io _{fX} \vep) = \io _X   (f^* \vep)$, which in turn shows that $X + f^* \xi \in L(f\inv E , \vep)$.  Therefore, $L(f\inv E, f^* \varepsilon) = \pbd$.  Since $f\inv E$ is a distribution on $M$, $L(f\inv E, \vep)$ is an almost Dirac stucture.  \\

\noindent (2) Note that if $\pbd$ is integrable, then being a submersion, locally $f: M \longrightarrow N$ looks like \\$ pr _1 : U\times V \longrightarrow U$ and $\pbd _{(p,q)} = \mathcal{D} _p \oplus T_qV$. So locally $f^\star \DD$ is a product of Dirac structures.  It is a general fact~\cite{gua} that a product $\DD _1 \times \DD _2$ of Dirac structures on a product manifold $U\times V$ is integrable.  This implies that $\DD _1 \times \DD _2$ is integrable if and only if $\DD_1$ and $\DD_2 $ are integrable. In this case, $\DD_2 = TV$.  Therefore, $f^\star D$ is integrable if and only if $\DD$ is integrable.  $\square$


    

\end{pf}

\begin{rem} In general, $f^\star \DD$ need not even be a smooth subbundle of $\mathcal{V}
  _M$, and the requirement that $f$ be a submersion is sufficient but not necessary to ensure
  that $\pbd$ is a Dirac structure.\\

\noindent  By way of counterexample, let $f: \mathbb{R} \longrightarrow \mathbb{R} ^2$ be $ x \mapsto (x,0)$, and let $\mathcal{D} = L(E,0)= E \oplus AnnE$, where $E = \mathbb{R}\cdot (\frac{\partial}{\partial x} + x\frac{\partial}{\partial y} )$. 
Obviously $\mathcal{D}$ is  maximal isotropic.  Since $E$ is a 1-dimensional distribution, it is integrable.  Also $d_E \vep = 0$.  It follows from Lemma~\ref{21sept1} that $\DD$ is integrable.  Therefore, $\DD$ is a Dirac structure.  However, the pullback at $0$ is $(\pbd )_0 = df_0 ^\star (\mathcal{D} _{(0,0)} )= \mathbb{R} \cdot
  \frac{d}{dx} _{|x=0}$, whereas if $x\neq 0$, $(\pbd )_x = df^\star
  (\mathcal{D} _{(x,0)}) = \mathbb{R} \cdot dx_{|x}$.  Therefore $\pbd$ is not a smooth subbundle of $\mathcal{V} _M$.

\end{rem}


\begin{prop}The Courant bracket on $(\mathcal V _G)_\cc$ gives $\cgts$ the
  structure of a Lie algebra.  This is the semidirect product
  $\cbla$, where Lie algebra 
  $\mathfrak g _\cc$ acts on $\mathfrak g^* _\cc$ by the coadjoint representation.  
\end{prop}
\begin{pf}
It must be shown that the Courant bracket is closed on left-invariant sections of $(\mathcal{V}_G )_\cc$.  The Courant bracket is invariant under diffeomorphisms $G \tilde \longrightarrow G$.  Left-invariant sections of $(\VV _G )_\cc$ are invariant under the diffeomorphisms $L_g$, $g\in G$.  It follows that the Courant bracket preserves left-invariance. For left-invariant sections, the formula for the Courant bracket is easily seen to be the same as $\cbla$. \\

  

\noi Since $[ \; , \; ]$ is clearly anti-symmetric, the only thing to be checked is that the restriction of $[ \; , \; ] $ to left-invariant sections of $(\VV _G )_\cc$ satisfies the Jacobi identity.  It is known~\cite{gua} that for sections A,B,C of $(\mathcal{V} _G)_\cc$, $Jac(A,B,C) = d(Nij(A,B,C))$, where $Nij(A,B,C) = 1/3(\langle [A,B], C \rangle + \langle [B,C], A \rangle + \langle [C,A], B \rangle )$. Then if A, B, C are constant sections (meaning left-invariant), Nij(A,B,C) is constant and therefore Jac(A,B,C)= 0. $\square$


\end{pf}

\begin{Def}We say that a maximal isotropic subalgebra of $\cbla$ is a \emph{Dirac Lie subalgebra} of $\cbla$ or a \emph{linear Dirac structure} on $\cgts $, and a maximal isotropic subspace of $\fg _\cc \oplus \fg ^* _\cc$ is called a \emph{linear almost Dirac structure}
\end{Def}

\begin{prop}\label{dhsprop3} 
 \begin{enumerate}
\item There is a bijection between linear (almost) Dirac
  structures on $\cbla$ and left-invariant (almost) Dirac structres on
  G.  
\item For a subspace $E \subset \fg _\cc$ and $\vep \in \wedge ^2 E^*$, let $\tilde E$ be the left-invariant distribution on $G$ determined by $E$, and let $\tilde \vep \in \wedge ^2 \tilde E ^*$ be the left-invariant 2-form determined by $\vep$.  The linear almost Dirac structure $L(E,\vep)$ corresponds to the almost Dirac structure $L(\tilde E, \tilde \vep)$ on $G$.  
\end{enumerate}
\end{prop}

\begin{pf}
We will show that if $L$ is a linear almost Dirac subspace of
  $\cgts $, then it determines a unique left-invariant almost Dirac
  structure $\LL$ on G such that $\LL _e = L$ and such that $\LL$ is integrable if and only if $L$ is a Lie
  subalgebra of $\cbla$.  We also show that, on the other hand, any left-invariant Dirac
  structure $\LL$ induces a linear Dirac structure $\LL _e$ of $\cbla$. \\

\noi That $L$ determines a G-invariant almost Dirac structure $\LL $ is obvious.  If $\LL $ is integrable, it is apparent from the definition of the Courant bracket on $\cbla$ that $L$ is a subalgebra, for if $X,Y \in L$, then $[X,Y]= [\tilde{X} , \tilde{Y}]_e \in L$ because $\LL $ is integrable.  \\

\noi Now suppose that $L$ is a subalgebra.  For sections $ X + \xi , Y + \eta$ of $(\VV _G)_\cc$ and $f \in \CC ^\infty (G) _\cc $,  it is a general fact of Courant brackets \cite{gua} that $[X + \xi , f(Y + \eta)] = f[A,B] + X(f)(Y+\eta) - \langle  X + \xi , Y + \eta \rangle df$.  Therefore, if $\LL$ is the left-invariant almost Dirac structure such that $\LL _e =L$, to check integrability of $\LL$, it suffices to check it on a frame of $\LL$.  With this in mind, choose any basis $X_1, ...X_n$ of $L$, which provides a left-invariant frame $\tilde X_1 , ..., \tilde X_n $ of $\LL$.  By definition of the Lie algebra structure on $\cbla$, $ [\tilde X_i , \tilde X_j] = \widetilde {[X_i , X_j]}$.  But since $L$ is a subalgebra, $[X_i , X_j] \in L$ and so $ [\tilde X_i , \tilde X_j] $ is a section of $\LL$.  Therefore, $\LL$ is integrable. \\

\noi The final claim is transparent because $L(\tilde E , \tilde \vep )_e = L(E,\vep)$. $\square$ 

\begin{Def} A subalgebra $E \subset \fg _\cc$ and $\vep \in \wedge ^2 E^*$ determine left-invariant $\tilde E$ and $\tilde \vep$ as in Proposition~\ref{dhsprop3}.  For $X,Y,Z \in E$, let $\tilde X , \tilde Y , \tilde Z$ be the left-invariant vector fields on $G$ with respective values $X,Y,Z$ at $e \in G$.  Define the differential $d_E : \wedge ^2 E^* \lra \wedge ^3 E^*$ by $d_E \vep (X,Y,Z) := d_{\tilde E} \tilde \vep (\tilde X , \tilde Y, \tilde Z)_e$.  
\end{Def}

\begin{rem} The differential $d_E$ is given by the following formula:  
\[ d_E \vep (X,Y,Z) = \varepsilon (X, [Y,Z]) + \varepsilon (Y, [Z,X]) + \varepsilon (Z,[X,Y]). \]
\noi With this formula, $d_E$ is the Lie algebra differential for $E$ \cite{wei}.
\end{rem} 

\begin{prop}\label{dhsprop2}  $L(E,\varepsilon)$ is a Dirac subalgebra if and only if E is a subaglebra of $\mathfrak{g} _\cc$ and $d_E \vep =0$.
\end{prop}

\begin{pf} This follows from Propositions~\ref{21sept1} and~\ref{dhsprop3}.

\end{pf}

\begin{rem}\label{h} If the Lie algebra cohomology of $E$ in degree $2$ vanishes (i.e. $H^2(E, \mathbb C) = 0$ ) , then for $\varepsilon \in \wedge ^2
  E^* $ we have
$d_E\vep =0$ if and only if 
  $\varepsilon = \phi \circ [\; , \; ]$ for some $\phi \in E^* $.  If E is a semisimple Lie algebra, then $H^2 (E
  ,\mathbb C) = 0 $ ~\cite{wei}. 

\end{rem}

\end{pf}

\subsubsection{Classification of Homogeneous Dirac Structures}\label{13may1}
Throughout the remainder of this section, assume that $K$ is connected.  
We establish a bijection between the G-invariant Dirac structures on
  $G/K$ and the set of Dirac subalgebras $L$ of
  $\cbla$ containing $\mathfrak k  _\cc$. This correspondence sends a Dirac
  structure $\mathcal D$ on G/K to $(\pi ^\star \mathcal D) _e$, and its inverse sends a subalgebra $L \subset \cbla$ to the G-invariant Dirac structure determined by $(d\pi _e)_\star L$.

\begin{thm}\label{dhsthm}
 Let $G$ be a Lie group and $K$ be a closed, connected subgroup.  There is a bijection between the G-invariant (almost) Dirac structures on
  $G/K$ and the set of (almost) Dirac subalgebras $L$ of
  $\cbla$ containing $\mathfrak k _\cc$. 
Any Dirac subalgebra of $\cbla$ is of the form
  $L(E,\varepsilon)$. The G-invariant Dirac stuctures on G/K are thus parameterized by
  pairs $(E, \varepsilon) $, where E is a Lie subalgebra of $\mathfrak
  g _\cc$ containing $\mathfrak k _\cc$, $d_E \varepsilon =0$, and $\vep _\sharp$ vanishes on $\mathfrak k$.

\end{thm}

\noi Observe that the integrability of an almost Dirac structure corresponds to whether the linear almost Dirac structure (i.e. maximally isotropic subspace) $L \subset \cbla$ is a subalgebra or simply a subspace.  Equivalently, integrability is the same as asking that $E \subset \fg _\cc$ is a subalgebra, rather than simply a subspace, and that $d_E \vep =0$.  \\

\begin{Def} Let $G$ be a Lie group and $K$ a closed subgroup.  A \emph{Dirac pair} is a pair $(E,\vep)$ satisfying the conditions of of Theorem \ref{dhsthm}.
\end{Def}
\noi Theorem~\ref{dhsthm} will be proven over the course of several lemmas.

\begin{lem}\label{kinv} $\,$ 
\begin{enumerate}
\item If $L$ is a Dirac subalgebra of $\cbla$, then $\mathfrak k _\cc \subset L
  $ if and only if $pr_{\mathfrak g ^*} L \subset Ann(\fk _\cc)$.
\item If $L= L(E,\vep)$ is a Dirac Lie subalgebra of $\cbla$ containing $\mathfrak k _\cc$, then
  $\varepsilon$ is Ad(K)-invariant. In other words, $\varepsilon([k,X],Y)+\varepsilon(X,[k,Y]) = 0$ for all $k\in \mathfrak{k} _\cc$ and $ X,Y \in \mathfrak{g}$.  Also $\vep _\sharp$ vanishes on $\fk _\cc$.
\end{enumerate}
\end{lem}

\begin{pf} $\,$
\begin{enumerate}
\item Since L is maximal isotropic, $L = L^\perp$, so $L \subset
  \mathfrak g _\cc \oplus Ann(\mathfrak k _\cc) $ if and only if $L = L ^\perp \supset
  (\mathfrak g _\cc \oplus Ann(\mathfrak k _\cc) )^\perp = \mathfrak k _\cc$.

\item If $L = L(E,\varepsilon ) $ is a Dirac subalgebra of $\cbla$
  containing $\mathfrak k _\cc$, then by part 1 of this lemma,
  $pr_{\mathfrak g ^* _\cc} L \subset Ann(\fk _\cc)$. For any $X \in E$, there
  exists  $X + \eta \in L$. Because L is isotropic, for any  $k
  \in \mathfrak k _\cc \subset L$, $0 =\langle k , X + \eta \rangle = \eta
  (k) = \varepsilon (X,k)$. Therefore $\varepsilon _\sharp$ vanishes on $\fk _\cc$.  Now by Proposition~\ref{dhsprop2}, $\varepsilon([k,X],Y)+\varepsilon(X,[k,Y]) = 0 \,$ for all $k\in \mathfrak{k} _\cc , X,Y \in \mathfrak{g} _\cc$. $\square$

\end{enumerate}
\end{pf}

\begin{lem}\label{dhslem1} If $\mathcal{D}$ is a G-invariant almost Dirac structure on $G/K$, then 
$\mathcal{D}$ is G-invariant if and only if $\pi ^\star \DD$ is G-invariant
\end{lem}
\begin{pf} 

\noi 
Let $X + d\pi_e^\ast\eta \in \pi ^\star \DD _e$.  This means that $d\pi_eX +\eta \in \DD _{\pi (e)}$.   Then $g \cdot (X + d\pi_e^\ast\eta) = (dL_g)_eX + (dL_g)_e^{-\ast} \circ (d\pi)_e^\ast \eta
=(dL_g)_eX + d\pi_g^\ast \circ (dl_g)_{\pi(e)}^{-\ast} \eta$, which is in $\pi ^\star \DD _g$ if and only if $d\pi_g \circ (dL_g)_eX + (dl_g)_e^{-\ast} \eta \in \mathcal{D} _{\pi(g)}$.  But  $d\pi_g \circ (dL_g)_eX + (dl_g)_e^{-\ast} \eta = (dl_g)_{\pi(e)} \circ d\pi_eX + (dl_g)_e^{-\ast} \eta = g 
\cdot (d\pi_eX + \eta)$.  Thus, $g \cdot (d\pi_eX + \eta) \in \mathcal{D} _{\pi(g)}$ if and only if $g \cdot (X + d\pi_e^\ast \eta ) \in \pi ^\star \DD _g$.  Therefore $\pi ^\star \DD$ is G-invariant if and only if $\mathcal{D}$ is G-invariant. $\square$
\end{pf}

\begin{lem}\label{dhslem2} Let $L \subset \cbla$ be any Dirac subalgebra containing $\fk _\cc$, and let $D \subset (\mathfrak{g}  /\mathfrak{k})_\cc \oplus (\mathfrak{g} /\mathfrak{k})_\cc ^\ast$  be any K-invariant linear Dirac structure.  Then 
\begin{enumerate}
\item $\pi^\star(\pi_\star L) = L$,
\item $\pi_\star(\pi^\star D) = D$, and
\item $\pi_\star L$ is K-invariant.  More generally, $L$ is $K$-stable if and only if $\pi _\star L$ is $K$-stable.  
\end{enumerate}
\end{lem}

\begin{pf}Here $\pi$ is used to denote $d\pi_e$.
\begin{enumerate}
\item By definition, $\pi^\star(\pi_\star L) = \{X + \pi^\ast \eta \st \pi X + \eta \in \pi_\star L\}$.  
But $\pi X + \eta \in \pi_\star L$ if and only if 
$X + \pi^\ast \eta \in L$.  So $\pi^\star(\pi_\star L) = \{X + \pi^\ast \eta \, |\, X + \pi^\ast \eta \in L \} $ \\
$= \{X + \omega \in L\, |\, \omega = \pi^\ast \eta \, for\, some\; \eta \}$
$ =  \{X + \omega \in L \st \omega \in im(\pi^\ast ) \}$.
  But $im \pi^\ast = Ann(\fk _\cc)$, so since $proj_{\mathfrak{g} _\cc ^\ast}L \subset Ann(\fk _\cc)$, we have $\pi^\star(\pi_\star L) = L$.

\item  By definition, $\pi_\star(\pi^\star D) = \{\pi X + \eta | X + \pi^\ast \eta \in \pi^\star D \}$.  Since $X + \pi^\ast\eta \in \pi^\star D$ if and only if $\pi X + \eta \in D$, we find that
 $\pi_\star(\pi^\star D) = \{\pi X + \eta | \pi X + \eta \in D \} = D$ because $\pi$ is surjective.  

\item If we consider the Ad representation (i.e. $Ad + Ad^{-\ast}$ representation) on $\cgts $ , then\\
$d\pi_e \circ Ad(k) = (dl_k)_{\pi (e)} \circ d\pi_e$, $Ad(k)^\ast \circ (d\pi_e)^\ast = (d\pi_e)^\ast \circ (dl_k)^\ast$, and $(d\pi_e)^\ast \circ (dl_k)_e^{-\ast} = Ad(k)^{-\ast} \circ (d\pi)_e^\ast$.  We get $k \cdot (X + (d\pi)_e^\ast\eta) = Ad(k)X + Ad(k)^{-\ast} \circ (d\pi)_e^\ast\eta = Ad(k)X + (d\pi_e)^\ast \circ (dl_k)_e^{-\ast}\eta$, which is in $L$ if and only if $d\pi_e \circ Ad(k) + (dl_k)_e^{-\ast}\eta \in \pi _\star L$.  Of course $d\pi_e \circ Ad(k) + (dl_k)_e^{-\ast}\eta = (dl_k)_{\pi (e)} \circ d\pi_e X + (dl_k)_e^{-\ast}\eta = k\cdot (d\pi_e X + \eta)$.  So $k\cdot (X + (d\pi_e)^\ast \eta) \in L $ if and only if $  k\cdot (d\pi_e X + \eta) \in \pi _\star L$.  The result is that $L$ is K-stable if and only if $\pi _\star L$ is K-stable.  But $L$ is necessarily $K$-stable because $\fk _\cc\subset L$. $\square$


\end{enumerate}
\end{pf}

\noi Now we prove Theorem~\ref{dhsthm}:

\begin{pf} First let $L = L(E,\varepsilon) \subset \cbla $ be a Dirac subalgebra containing $\fk _\cc$, and let $D = \pi _\star L$.  By part 3 of Lemma~\ref{dhslem2}, D is K-invariant and so defines a G-invariant almost Dirac structure $\mathcal{D}$ on $G/K$.  Lemma~\ref{dhslem1} implies that since $\DD$ is equivariant, $\pi ^\star \DD$ is a left-invariant almost Dirac structure on $G$.  Now by Lemma~\ref{dhslem2}, $L = \pi ^\star \pi _\star L = \pi ^\star D = (\pi ^\star \DD )_e$.  Hence, $\pi ^\star \DD$ is the left-invariant almost Dirac structure determined by $L$.  By Proposition~\ref{dhsprop3}, since $L$ is a subalgebra, $\pi ^\star \DD$ is integrable.  Then by Proposition \ref{pdsprop2}, $\DD$ is integrable.  This shows that from the linear Dirac structure $L$ we obtain an equivariant Dirac stucture $\DD$ on G/K, and $(\pi ^\star \DD ) _e = L$.\\

\noi In the other direction, an equivariant Dirac structure $\DD$ on $G/K$ yields $L = (\pi ^\star \DD )_e$.  Since $\DD$ is an equivariant Dirac structure, so is $\pi ^\star \DD$ by Lemma~\ref{dhslem1}. Hence, by Proposition~\ref{dhsprop3}, $L \subset \cbla $ is a Dirac subalgebra.  It follows from the definition of pullback $\pi ^\star$ that $\fk _\cc \subset L$.  This shows that from an equivariant Dirac structure $\DD$ on $G/K$, we obtain a Linear Dirac structure $L$,  which contains $\fk$.  To see that the correspondence is a bijection, we need only to observe that $\pi _\star L = \pi _\star (\pi ^\star \DD )_e  = \pi _\star \pi ^\star (\DD _{\pi (e)}) = \DD _{\pi (e)}$ by Lemma~\ref{dhslem2}.  This also shows that $L$ is $K$-stable.  Naturally, $\DD _{\pi (e)}$ is the subspace of $\cgts $ which determines the Dirac structure $\DD$.  This establishes the bijection.  \\

\noi The description in terms of pairs $(E,\vep)$ follows directly from Proposition~\ref{dhsprop2} and the fact that any linear Dirac structure, $L \subset \cbla$ is of the form $L= L(E,\vep)$. $\square$

\end{pf}

\subsubsection{Real Dirac Structures}

\noi The real equivariant Dirac structures on $G/K$ are those $L \subset \cbla $ such that $L = \overline L$ or equivalently, those $L$ such that $L = D_\cc$ for some $D \subset \fg \oplus \fg ^*$.  But by considering only $\VV _G$, $\fg$, and $\fk$ instead of their complexifications, the following theorem follows in exactly the same way as Theorem \ref{dhsthm}.

\begin{thm} 
 Let $G$ be a Lie group and $K$ be a closed, connected subgroup.  There is a bijection between the G-invariant (almost) Dirac structures on
  $G/K$ and the set of (almost) Dirac subalgebras $L$ of
  $\bla$ containing $\mathfrak k $. 
Any Dirac subalgebra of $\bla$ is of the form
  $L(E,\varepsilon)$. The G-invariant Dirac stuctures on G/K are thus parameterized by
  pairs $(E, \varepsilon) $, where E is a Lie subalgebra of $\mathfrak
  g $ containing $\mathfrak k $, $d_E \varepsilon =0$, and $\vep$ vanishes on $\mathfrak k$.
\end{thm}

\begin{cor}\label{c} By correspondence given in Theorem~\ref{dhsthm} between Linear Dirac structures
  $L = L(E,\varepsilon) \subset$ \gts and equivariant Dirac
  structures on G/K,
\begin{enumerate}
\item The symplectic structures are all $L(\mathfrak{g} , \varepsilon )$ such that $Ker (\tep ) = \mathfrak{k}$.
\item The presymplectic structures are all $L(\mathfrak{g}, \varepsilon )$.
\item The Poisson structures are all $L(E , \varepsilon )$ such that  $Ker (\tep ) = \mathfrak{k}$.
\item There is a bijection between G-invariant real Dirac structures on G/K and G-invariant
  presymplectic structures on the spaces H/K for connected subgroups $H$ of $G$
  containing K.  With this bijection, a linear Dirac structure $L(E,\varepsilon)$ corresponds to
  a presymplectic structure on H/K, where H is the connected Lie subgroup
  corresponding to Lie subalgebra E.  Only when $Ker (\tep ) = \mathfrak{k}$ does $\varepsilon$ gives a symplectic structure.  Hence there is a bijection between G-invariant real Dirac structures on G/K and G-invariant
  presymplectic structures on the spaces H/K for connected subgroups H
  containing K.
\end{enumerate}
\end{cor}
\begin{pf}  
\begin{enumerate}
\item Let $L = L(E,\vep) \subset \bla$ be a subalgebra containing $\fk$ such that $Ker \, \vep _\sharp = \fk$, and let $\mathcal{L} = L(TG, \tilde \varepsilon)$ be the left-invariant Dirac structure on $G$ determined by L.  Then $D = d\pi _\star L = L(\mathfrak{g} / \mathfrak{k} , \omega _0 )$
 gives equivariant Dirac structure 
$\mathcal{D} = L(T(G/K) , \omega )$.
  By Proposition~\ref{d},
 $L(TG, \pi ^* \omega ) =  \pi ^\star \mathcal{D} = \mathcal{L}= L(TG,\tilde \varepsilon)$,
 so $\tilde \varepsilon = \pi ^* \omega$. 
 By Prop~\ref{dhsprop3},
 $\tilde \varepsilon$ is closed and left-invariant.  This implies that
 $0 = d_E \tilde \varepsilon = d \tilde \varepsilon = d \pi ^* \omega = \pi ^* d \omega$, whence $d\om = 0$ since $d\pi ^*$ is injective.  To show that $\om$ is non-degenerate, it is enough to show that $\om _0$ is non-degenerate, by G-invariance.  But $\vep _\sharp = (d\pi ^* \om _0 )_\sharp = (d\pi )^* (\om _0)_\sharp d\pi $.  Since $d\pi$ is surjective and $d\pi ^*$ is injective, $Ker( (\om _0)_\sharp ) = 0 $ if and only if $Ker (\tep ) = \fk$.  
 Conversely, given symplectic $\omega$, the above argument may be run backwards to show that for $\varepsilon = d\pi ^* \omega _0$, one has $Ker \tep = \mathfrak{k} $.
\item This follows in the same manner as part 1. The only difference is, using the notation from part 1 of this proof, $\om$ may not be non-degenerate.  

\item Dirac structures $L(E,\omega) \subset \VV _{G/K}$ correspond to $L(d\pi \inv E, \pi ^* \omega )_e \subset$ \gts .  Visibly, $\omega$ is non-degenerate if and only if $\omega _{\pi (e)} $ is non-degenerate, which happens exactly when $Ker ( ({d\pi _e ^* \omega})_\sharp ) = Ker (d\pi _e ) = \mathfrak{k} $.\\

\noi Now it suffices to show that Dirac structures $L(E,\omega) \subset \VV _{G/K}$ with $\omega$ non-degenerate are exactly those corresponding to Poisson structures.  This should make sense because a Poisson structures gives an integrable distribution $E$ and a symplectic structure $\omega$ on the leaves of the foliation determined by $E$.  \\

To see how this works, let $L(E,\omega) \subset \VV _{G/K}$ be a Dirac structure and $\omega$ be non-degenerate. 
Since $T^*(G/K) \twoheadrightarrow E^*$ and $\om _\sharp : E \tilde \lra E^*$, it follows that $pr_{T^* (G/K)} L(E,\om) = T^* (G/K)$.  Thus, $L(E,\om) $ is of the form $L(\be , T^* G/K)$ for some 2-form $\be$ on $G/K$.  It is a result of \cite{gua} and \cite{vai} that almost Dirac structures of the form $L(\be , T^*(G/K))$ are Dirac structures precisely when $\be$ is a Poisson bivector.

\item Follows easily. $\square$ 
\end{enumerate}
\end{pf}

\section{Generalized Complex Structures on Homogeneous Spaces}\label{GCclassification}

\noi We delineate the conditions for $L \subset \cbla $ or a Dirac pair $(E, \vep)$ to represent an generalized complex structure in Corollary~\ref{f} and Proposition~\ref{b}.  Again we emphasize that in addition to the closed integrability condition, we now require the genericity condition $L \cap \overline L = \fk _\cc$.  
\subsection{Classification}

\noi As a corollary to Theorem \ref{dhsthm}, we have: 

\begin{cor}\label{f} Let $G$ be a Lie group and $K$ a closed, connected subgroup.  If $L \subset \cbla$ is a linear (almost) Dirac structure containing $\csts$, then L represents a
  generalized (almost) complex structure on $G/K$ if and only if $L \cap \overline{L} =
  \mathfrak k _{\mathbb C} $. 
\end{cor}

\begin{pf} If $\mathcal D $ is a equivariant complex Dirac
  structure on G/K, then $\mathcal D$ is a generalized complex
  structure if and only if $\mathcal D \cap \overline{\mathcal D} =
  0$, which happens exactly when $\mathcal D _{\pi (e)} \cap
  \overline{\mathcal D } _{\pi (e)} = 0$.  Note that $\pi ^\star $ may be applied to
  arbitrary subspaces and not just maximal isotropic subspaces, so we observe immediately that $\pi
  ^\star (0) = \mathfrak k _{\mathbb C}$.  On the other hand, if V is any nonzero subspace of $
  (\mathfrak g / \mathfrak k )_{\mathbb C} \oplus ((\mathfrak g /
  \mathfrak k )_{\mathbb C} )^* $, it follows from the definition of $\pi ^\star$ together with surjectivity of $d\pi _e$ ( and injectivity of its dual) that if $\pi ^\star V = \csts$, $V \subset \fg / \fk$, whence $\csts = \pi ^\star V = \pi \inv V$ and so $V=0$.  We conclude that $\pi ^\star V = \csts$ precisely when $V=0$.  \\
  
  \noi For
  any two subspaces $V,W \subset  
  (\mathfrak g / \mathfrak k )_{\mathbb C} \oplus (\mathfrak g /
  \mathfrak k )_{\mathbb C} ^* $, $\pi ^\star (V \cap W) \subset
  \pi^\star V $ implies that $\pi ^\star (V \cap W ) = \pi ^\star V
  \cap \pi ^\star W $.  Finally, the observation that $\pi ^\star
  \overline {V}  = \overline{\pi ^\star V} $ completes the proof. $\square$

\end{pf}


\begin{prop}\label{19octa} In the correspondence of Theorem~\ref{dhsthm} and
  Corollary~\ref{f}, 
\begin{enumerate}

\item The complex structures are given by all $L(E,0 ) \subset \cts \oplus \cts ^* $ such that $E +\overline{E} = \cts$ and $E \cap \overline{E} = \csts $.  Thus there is a bijection between G-invariant complex structures on $G/K$ and subalgebras $E \subset \fg _\cc$ such that $E + \overline{E} = \cts$ and $E \cap \overline{E} = \csts $.  This correspondence can be extended to a bijection between $G$-invariant almost complex structures on $G/K$ and subspaces $E \subset \fg _\cc$ such that $E + \overline E = \fg _\cc$ and $E \cap \overline E = \fk _\cc$.  
\item The symplectic structures are all $L(\cts , \varepsilon ) $ such that $Ker (\tep ) = \csts $ and $\varepsilon$ is purely imaginary (i.e. $\varepsilon = i\omega$ for some real 2-form $\omega$).

\end{enumerate}

\end{prop}

\begin{pf}
\noindent We know that complex structures are exactly generalized complex structures of the form $L(E,0)$.  Part (1) now follows from Theorem~\ref{dhsthm} and
  Corollary~\ref{f}.\\

\noi (2) This follows in the same way as part (1) of Corollary~\ref{c}. $\square$

\end{pf}

\begin{rem} The condition that a complex structure be G-invariant means that the distribution it defines in $T(G/K) \otimes \mathbb{C}$ is G-invariant.  This is equivalent to requiring that each $l_g : G/K \longrightarrow G/K $ is holomorphic.  This equivalence follows directly from the fact that a map $f: M \longrightarrow N$ of complex manifolds is holomorphic if and only if $df_p (T_p ^{1,0}M) \subset T_{f(p)} ^{1,0}N $ for all $p \in M$, where $T^{1,0}M$ and $T^{1,0}N$ are the holomorphic tangent bundles.  Here D = $T^{1,0} (G/K) \oplus T^{* , (0,1)} (G/K)$.  This is also true for almost complex structures.  
\end{rem}

\begin{prop} In the correspondence of Theorem~\ref{dhsthm} and
  Corollary~\ref{f}: \begin{enumerate}
\item Any $L(E,i\omega ) $, where $\omega$ is the
  restriction to E of a real 2-form on $\mathfrak{g}$, gives a
  symplectic structure on some H/K for some subgroup H of
  G. Specifically, If D is the subalgebra of $\mathfrak g$ such that
  $D_{\mathbb C} = E \cap \overline{E} $, then H is the connected Lie subgroup
  of G with Lie algebra $\mathfrak k$.  
\item The symplectic structures are all $L(E, i \om)$, where $\om \in \wedge ^2 \fg ^*$.
\item Any generalized complex structure $L(E,\vep ) $ gives a complex
  structure $L(E,0) $ on G/H, where H is the Lie subgroup from part 1, as long as $H$ is closed.
\end{enumerate}
\end{prop}
\begin{pf}
\begin{enumerate}
\item In order for $L \cap \overline{L} = \csts$, $Ker (\om _\sharp) _{|E\cap \overline{E} }  = \csts $.  $E\cap \overline{E} = D_\mathbb{C} $ for some subalgebra $D\subset \mathfrak{g} $ contining $\fk$ because $\overline{ E \cap \overline E} = E \cap \overline E$.  Since D is a Lie subalgebra, it corresponds to some Lie subgroup H containing K.  $\omega$ is clearly non-degenerate on $D / \mathfrak{k} $ and closed.  Therefore it determines a symplectic form on H/K.

\item It is clear that these are the symplectic structures.  The condition that $Ker \,  \om _\sharp = \fk _\cc$ is equivalent to non-degeneracy of the symplectic form, but it is also equivalent to $L \cap \overline L = \csts$.  
\item This follows from Proposition~\ref{19octa} . $\square $ 
\end{enumerate}
\end{pf}

\subsection{Analysis of Conditions on E and $\varepsilon$}
Although Corollary~\ref{f} gives a description of equivariant generalized complex
structures on G/K, it is often more useful to describe them in terms
of a subalgebra $E$ and $\varepsilon \in \wedge ^2 E^* $.\\

\noindent By $\mathbb{C}$-linear extension, $\mathfrak{g} ^* \hookrightarrow (\cts )^*$, and this gives an isomorphism $ (\mathfrak{g} ^*)_\mathbb{C} = \mathfrak{g} ^* \oplus i\mathfrak{g} ^* \simeq (\cts) ^*$.
  If $\alpha \in \ctsd$ and $X \in \cts$, it may be easily
verified that $\overline{\alpha} ({X} ) = \overline{ \alpha
  \overline{ (X)} }$, where $\alpha \mapsto \overline{\alpha}$ denotes
conjugation.  It follows immediately that $Ann(\overline{E}) =
\overline{Ann(E)}$.  
  We define $\overline{\varepsilon} (X,Y) =
\overline{\varepsilon (\overline{X}, \overline{Y} )}$.

\begin{prop} $\overline{L} = L(\overline{E}, \overline{\varepsilon} )$.
\end{prop}
\begin{pf} The result follows by observing the following equalities. \\
$\overline L = \{ X + \xi \st X \in \overline E \; , \;\overline \xi _{|E} = \tep \overline X \}$\\
$= \{ X + \xi \st X \in \overline E \; , $ for all $Y \in E  \; \overline \xi  (Y) = \vep (\overline X , Y) \}$ \\
$= \{ X + \xi \st X \in \overline E \; , \; $ for all $Y \in \overline E \; \overline {\xi (Y)} = \overline \xi  (\overline Y) = \vep (\overline X , \overline Y) \}$\\
$=  \{ X + \xi \st X \in \overline E \; , \; $ for all $Y \in  \overline E  \; {\xi (Y)} = \overline{ \vep (\overline X , \overline Y) } \}$\\
$= \{ X + \xi \st X \in \overline E \; , \; $ for all $Y \in  \overline E  \; {\xi (Y)} = \overline \vep (X,Y) \}$\\
$= L(\overline E, \overline \vep)$. $\square$ 

\end{pf}

\noi {\bf Notation :} We let the symbol $\setminus $ denote the difference of two sets. 

\begin{prop}\label{a} $L \cap \overline{L} = \csts$ if and only if
\begin{enumerate}
\item $E +\overline{E} = \cts $, and

\item $X \in (\overline{E} \cap E )\setminus \csts$ and $X + \xi \in L$ implies $\overline{X} + \overline{\xi} \notin L$ (i.e. $pr_{\fg \oplus \fg ^* \lra \fg} L \cap \overline L = \csts$).
\end{enumerate}
\end{prop}

\begin{pf} The first condition is equivalent to  $Ann(E) \cap Ann(\overline{E}) = 0$, since $Ann(V) \cap Ann(W) = Ann(V+W)$ for any subspaces V,W.
First suppose that $L \cap \overline{L} = \csts$.  $Ann(E) \subset L$
and $\overline{Ann(E)} = Ann{\overline E} \subset \overline{L} $.
This implies that $Ann(E) \cap Ann(\overline{E} ) \subset L\cap
\overline{L} = \csts$.  But $Ann(E) \cap Ann(\overline{E} ) \subset
\mathfrak g_{\mathbb C} ^* $ and so $Ann(E) \cap Ann(\overline{E}) =
0$.  Part (2) is obviously true.\\

\noi Now suppose that (1) and (2) are true.  In order to show that $L \cap \overline L \subset \csts$, let $A = X + \xi \in L \cap \overline L$.  We know that $X \in E \cap \overline E$, and (2) implies that $X \in \fk _\cc$.  Hence, $\tep X = 0$, which implies that $\xi  \in Ann(E)$.  Also $\overline X + \overline \xi \in L$, so $\overline \xi \in Ann(E)$ and therefore $\xi \in \overline {Ann(E)} = Ann (\overline E)$.  But now $\xi \in Ann(E) \cap Ann(\overline E) = Ann (E + \overline E) = Ann(\cts) =0$.  Therefore $A= X$.  By condition (2), $X \in \csts$.  Thus, $L \cap \overline L = \csts$.  $\square$

\end{pf}

\begin{prop}~\label{g} Suppose that $E + \overline{E} = \cts$.  Then $L \cap \overline{L} = \csts$ if and only if for all $X \in E \cap \overline{E}$, if $\varepsilon(X,Y) - \overline{\varepsilon (\overline{X} , \overline{Y} ) } =0  $ for all $Y \in E \cap \overline{E}$, then $ X \in \csts$.

\end{prop}
\begin{pf}Suppose that $L \cap \overline{L} = \csts$, and let $X \in E \cap \overline E \setminus \csts$.  If $\psi , \xi \in \fg ^*$ such that $\xi _{|E} = \tep X$ and $\psi _{|E} = \tep \overline X $, then the condition on $\vep$ asks that there exists $Y \in E \cap \overline E$ such that $\vep (X,Y) \neq \overline {\vep(\overline X , \overline Y )}$ can be expressed as $\overline \xi - \psi \notin Ann(E \cap \overline E)$.  This is equivalent to $\tep X (Y) \neq \overline { \tep \overline X (\overline Y)}$ or in other words, $\xi (Y) \neq \overline \psi (Y)$.  This means exactly that $(\overline \xi - \psi ) (Y) \neq 0$, or $ \overline \xi - \psi \notin Ann(E \cap \overline E)$.  Suppose now, for the sake of contradiction, that $\overline \xi - \psi \in Ann(E \cap \overline E) = Ann E + Ann \overline E$.  Then $\overline \xi - \psi = \alpha + \beta \in Ann(E) + Ann(\overline{E} )$, and $\overline{\beta} \in \overline{Ann(\overline{E} )} = Ann(E) \subset L$, which implies that $X + \xi - \overline{\beta} \in L$.  In this case, $\overline X + \overline \xi - \be \in \overline L$.  However,  $\overline X + \overline \xi - \be = \overline X + \psi + \al + \be - \be = \overline X + \psi + \al \in L$, whence $\overline X + \overline \xi - \be \in L \cap \overline L = \csts$.  This is a contradiction, because it was assumed that $X \in E \cap \overline E \setminus \csts$.  Therefore, the fact that $L \cap \overline L = \csts$ implies the desired condition on $\vep$.  \\


\noi Now suppose that for all $ X \in E \cap \overline{E} $, if $\varepsilon(X,Y) - \overline{\varepsilon (\overline{X} , \overline{Y} ) } =0  \; $ for all $ Y \in E \cap \overline{E}$, then $ X \in \csts$. Suppose that $X + \xi \in L \cap \overline L$ with
 $X \in E \cap \overline E$.  We aim to show that $X \in \csts$.  This is sufficient by Proposition~\ref{a}.  We know that
 $\xi _{|E} = \tep X$ and $\overline \xi _{|E} = \tep (\overline X)$ because $X + \xi $ and $\overline X + \overline \xi $ both lie in $L$.  If $X \notin \csts$, then there exists $Y \in E \cap \overline E$ such that $\vep (X,Y) \neq \overline{ \vep (\overline X , \overline Y)}$.  However, $\vep (X,Y) = \xi (Y) = \overline { \overline \xi ((\overline Y))} = \overline{ \vep (\overline X , \overline Y)}$, which is a contradiction.  Therefore $X \in \csts$ and $L \cap \overline L = \csts$. $\square$

\end{pf}

\begin{prop}\label{b} Let $G$ be a Lie group and $K$ a closed, connected subgroup of $G$.  There is a bijection between $G$-invariant generalized complex structures on G/K and pairs \\
$(E, \varepsilon) $, E a subalgebra of $\cts$ and $\varepsilon \in \wedge ^2 E^*$, such that 
\begin{enumerate}
\item $\csts \subset E$,
\item $E+ \overline{E} = \cts$,
\item $d_E \varepsilon =0$,
\item $\tep (\mathfrak k ) =0$, and
\item For $X \in E \cap \overline E$, if $\varepsilon(X,Y) - \overline{\varepsilon (\overline{X} , \overline{Y} ) } =0  $ for all $Y \in E \cap \overline{E}$, then $ X \in \csts$.
\end{enumerate}
\end{prop}

\begin{pf} The proof is immediate from Proposition~\ref{g} and~\ref{a}  and
  Corollary~\ref{f}. 
\end{pf}

\begin{Def} For a homogeneous space $G/K$, pairs $(E, \vep )$ satisfying the conditions of Proposition~\ref{b} are called \emph{generalized complex pairs} or \emph{GC-pairs}. 
\end{Def}

\begin{rem} Conditions (1) and (4) of Proposition~\ref{b} are conditions for $L(E, \vep)$ to represent an almost Dirac stucture on $G/K$.  The requirement that this almost Dirac structure is integrable is condition (3) together with the requirement that $E$ is a Lie subalgebra.  Finally, conditions (2) and (5) ensure that $L(E, \vep)$ is a generalized complex structure. 
\end{rem} 
  
\begin{rem}  Condition (5) of Proposition~\ref{b} simply asks that $Ker((\vep - \overline \vep )_\sharp : E \cap \overline E \lra E \cap \overline E ) = \fk _\cc$.  Condition (5) may be stated in yet another way.  We may extend $\vep \in \wedge ^2 E^*$ to some $B \in \wedge ^2 \fg _\cc ^*$.  If $B = B_r + i B _i$ is the decomposition of $B$ into real and imaginary parts, then $\vep = \vep _r + \vep _i$, where $\vep _r$ and $\vep _i$ are the real and imaginary parts of $\vep$ respectively (i.e. the restrictions of $B_r$ and $B_i$).  Then $\varepsilon(X,Y) - \overline{\varepsilon (\overline{X} , \overline{Y} ) } =   \vep _i (X,Y) $.  Since $\tep (\csts) =0 $, $\vep _i $ defines a linear map
 $\hat{\varepsilon _i} : E/\csts \times E/\csts \longrightarrow
E/\csts $.   Condition 5 of Prop~\ref{b} may be restated in the
  following way: $\hat{\varepsilon _i} $ is non-degenerate when restricted to $E \cap \overline E /\csts \times E \cap \overline E /\csts$.

\end{rem}

\subsection{B-transformations}
\noindent Recall that for a 2-form $B$ on a manifold $X$, the map $TX \stackrel{ B_\sharp}{\lra} T^*X \hookrightarrow \VV _X$ can be extended by $0$ on $T^*X$ to give a map $\VV _X \lra \VV _X$.  The exponential of this map is called a $B$-transform and is denoted by $e^B$.  For a Dirac structure $L$ on $X$, applying the Courant algebroid automorphism $e^B$ to $L$ gives another Dirac structure $e^B L$, which we will call the \emph{B-transform of $L$} under the 2-form $B$.  Since we are only considering G-invariant Dirac structures on $G/K$, we only consider those B-transforms which transform
 equivariant Dirac structures to equivariant Dirac structures.

\begin{lem}\label{i} Let $\mathcal{D} $ be an equivariant Dirac structure on G/K and 
 $L(E,\varepsilon ) = \pi ^\star \mathcal{D} _e  \subset
 \cgts$. 
\begin{enumerate}
\item For a 2-form $B \in \Omega ^2 (G/K) $, $\pi ^\star (e^B \mathcal D )
 = e^{\pi ^* B} \pi ^\star \mathcal D $. 
\item If $\eta \in \wedge ^2 Ann(\mathfrak{k} ) \subset  \wedge ^2 \mathfrak g ^* $ and $d\eta = 0$, then $e^\eta L(E,\vep) = \pi ^\star (e^{\om}\mathcal{D})$ for some 2-form $\om$ on $G/K$.
\end{enumerate}

\end{lem}

\begin{pf} If V is a subspace of $\mathfrak g$ or a K-invariant
 subspace of $\mathfrak g / \mathfrak k $, let $\tilde V $ denote the
 corresponding distribution.  Similarly, if $\omega \in \wedge ^2
 \mathfrak g ^*$ or $\wedge ^2 (\mathfrak g / \mathfrak k ) ^*$ and
 K-invariant, let $\tilde \omega$ denote the corresonding 2-form.\\

\noi The first claim is that for a 2-form $B \in \Omega ^2 (G/K)$, $\pi ^\star (e^B \mathcal D )
 = e^{\pi ^* B} \pi ^\star \mathcal D $.  
Suppose that $\mathcal D = L(\tilde F , \tilde \omega ) $.
 Recall that $L(\tilde E , \tilde \varepsilon ) = \pi ^\star \mathcal
 D = L(\pi \inv \tilde F , \pi ^* \tilde \omega ) $.  If $i: \tilde E
 \longrightarrow TG $ and $j: \tilde F \longrightarrow T(G/K)  $ are
 inclusions, then $\pi \circ i = j \circ \pi _{|\tilde E} $.  
$\pi ^\star (e^B L(\tilde F , \tilde \omega )) = \pi ^\star L(\tilde F
 , \tilde \omega + j^*B ) = L(\tilde E , \pi ^* \tilde \omega + \pi ^*
 \circ j^* B ) = L(\tilde E , \tilde \varepsilon + i^* \circ \pi ^* B
 ) = e^{\pi ^* B} L(\tilde E, \tilde \varepsilon ) $.\\

\noi Now let $\eta \in \wedge ^2 Ann(\mathfrak k ) $ be closed, so  $\eta =
d\pi _e ^* B$ for a (unique) $B \in \wedge ^2
 (\mathfrak g / \mathfrak k ) ^*$.  It has been shown already that B
is left K-invariant if and only if $\eta$ is Ad(K)-invariant.
However, since $\eta$ is closed and vanishes on $\mathfrak k$, $\eta$
is automatically Ad(K)-invariant, as was stated in Remark~\ref{kinv}.
Therefore, B yields a G-invariant 2-form $\tilde B$, and $\pi ^*
\tilde B = \tilde \eta $ so that $0 = d\tilde \eta = d \pi ^* \tilde B = \pi
^* d \tilde B  $, which implies $dB = 0 $ since $\pi ^*$ is
injective.  From the first part of this proof, $e^\eta L(E,\varepsilon
) = e^{\tilde \eta} L(\tilde E, \tilde \varepsilon ) _e = e^{\pi ^*
  \tilde B} L(\tilde E , \tilde \varepsilon  ) _e = \pi ^\star
(e^{\tilde B} \mathcal D )_e $. $\square$

\end{pf}


\noi We now provide a sufficient condition so that every B-transformation
  taking an invariant Dirac structure $\mathcal D$ on G/K to another invariant Dirac structure is given by a $G$-invariant 2-form.

\noi {\bf Notation :} Recall that for a Lie algebra $\fg$ over $\cc$, $H^2 (\fg , \cc) $ denotes the Lie algebra cohomology in degree 2.  We will let $Z^2 (\fg , \cc) $ denote degree 2 cocycles in the the complex $Hom _\cc (\wedge ^i \fg , \cc)$ which defines the Lie algebra cohomology \cite{wei}.

\begin{prop} Let $\mathcal{D} $ be an equivariant Dirac structure on G/K and 
 $L(E,\varepsilon ) = \pi ^\star \mathcal{D} _e  \subset
 \cgts$.  If $H^2(\fg , \cc )$ surjects onto $H^2
  (E , \cc) $, then every equivariant B-transform of $\DD$ is of the form
 $e^\eta L(E, \vep ) $ for some $\eta \in  \wedge ^2
 \mathfrak g ^* $ such that $d\eta = 0$.  For instance, if E is semisimple,
  then every B-transform is of this type.
\end{prop}
\begin{pf} First we observe that $H^2 (\fg , \cc) $ surjects onto $H^2 (E , \cc)$ if and only if $Z^2 (\fg , \cc) $ surjects onto $Z^2 (E , \cc)$.  Let $B \in \Omega ^2 (G/K) $ be a closed 2-form such that
 $e^B \mathcal D $ is G-invariant.  By Lemma~\ref{i}, $\pi ^\star (e^B \mathcal D )
 = e^{\pi ^* B} \pi ^\star \mathcal D $. In fact, $e^{\pi ^* B} \pi
 ^\star \mathcal D $ is determined by its value at e.   If $\omega =
 \pi ^* B_e $, then  $(e^{\pi ^* B} \pi
 ^\star \mathcal D )_e = e^\omega L(E,\vep )$.  Since $e^\om L$ is a Dirac structure, $d_E \omega _{|E\times E} = 0
 $, so by assumption, there exists some $\eta \in Z^2(\fg  ,
 \cc )$ which agrees with $\omega$ on $E \times E$.\\

\noi If $E$ is semisimple, then by Remark~\ref{h}, any $\vep \in Z^2 (E , \cc)$ is of the form $\phi \circ [ \; , \; ]$ for some $\phi \in E^*$.  But $\fg ^*$ surjects onto $E^*$, so there is some $\tilde \phi \in \fg ^*$ which restricts to $\phi$.  Therefore, letting $\eta = \tilde \phi \circ [ \; , \; ]$ gives the desired result. $\square$ 
\end{pf}

\subsection{Quotients by Disconnected Subgroups} 
\noi So far we have only considered homogeneous spaces $G/K$, where $K$ is a closed, connected subgroup.  In this case, equivariant Dirac structures are given by Dirac pairs, described in Theorem~\ref{dhsthm}, and generalized complex structures are given by GC-pairs.  When $K$ is disconnected, however, in order for a Dirac subalgebra $L \subset \cbla$ to give a Dirac structure on $G/K$, it must be $K$-invariant.  If $K$ is connected, then $K$-invariance follows from $\fk$-invariance, but this is not necessarily the case when $K$ is disconnected.  

\begin{lem} $G$-invariant generalized complex structures on $G/K$ are given by $K$-stable subalgebras of $\cbla$ containing $\fk _\cc$.  A Dirac Lie subalgebra $L = L(E,\vep) \subset \fg _\cc ^* \rtimes \fg _\cc$ is $K$-invariant if and only if $E$ and $\vep$ are $K$-invariant.
\end{lem}

\begin{pf} The proof is exactly the same as the proof of Theorem~\ref{dhsthm} except that we must assume that $L$ is $K$-stable, which will no longer automatically follow from $L$ containing $\fk$.  A simple calculation shows that $K$ preserves $L$ if and only if $K$ preserves $E$ and $\vep$.  $\square$
\end{pf}

\noi Thus, the general version of Theorem~\ref{dhsthm} and Proposition~\ref{b} is:
\begin{thm} Let $G$ be a Lie group and $K$ a closed subgroup.  The G-invariant complex Dirac stuctures on G/K are parameterized by $K$-invariant Dirac pairs $(E,\vep)$, and the G-invariant generalized complex structres are parameterized by $K$-invariant GC-pairs.  
\end{thm}

\begin{rem} If $G$ is connected and $\vep = d_E \phi = \phi \circ [ \, , \, ]$ for some $\phi \in E^*$, then $\vep$ is $K$-invariant if and only if $\phi _{|[E,E]}$ is $K$-invariant.  
\end{rem}

\begin{rem}Let $K^0$ be the identity component of $K$.  If $L \supset \fk _\cc$, then we already noted that $K^0$-invariance is automatic.  To check $K$-invariance, one need only check invariance under a discrete subset of $K$, namely invariance under a representative of each coset in $K/K^0$.  
\end{rem}

\section{Quotients of Compact Groups by Connected Subgroups of Maximal Rank}\label{26nov9} 
\noi In this section let G be a compact group and K be a connected subgroup containing a 
 Cartan subgroup $C$.  We wish to classify the generalized complex
 structures on G/K by listing all GC-pairs $(E,\varepsilon)$.  Subsection \ref{20junea} will focus on the case when $K=C$, and the subsequent subsections will consider any $K \supset C$.  We first notice that since $G$ is compact, its Lie algebra $\fg$ is reductive.  Since $\fk $ contains a Cartan subalgebra $\fc = Lie (C)$, $\fk $ contains the center of $\fg$.  Thus, GC-pairs for $G$ are the same as those for $G/Z_G$, which is semisimple.  Henceforth, we will assume that $G$ is semisimple so that the notation is less cumbersome. We will show that E must be a
 parabolic subalgebra and that $\varepsilon $ is exact in the sense of Remark~\ref{h}.  Proposition~\ref{5augb} gives a full explanation.  We also determine in Section~\ref{5augc} that equivariant generalized complex structures on $G/K$ up to B-transform can be thought of as a symplectic structure on a subgroup together with a complex structure on the quotient of $G$ by that subgroup.  Finally, we explain the geometric structure of the moduli of equivariant Dirac structures on $G/K$ in the final subsection (Proposition~\ref{m}). 
 \\

\noindent  Let $\mathfrak{h} = \fc _\cc$ and
 $\mathfrak{l}_E : = E \cap \overline{E}$.  Conjugation on $\cts$ with
 respect to $\mathfrak{g}$ will be denoted by $\sigma$ or $x \mapsto
 \overline{x} $.  Finally, let $\Delta$ denote the set of roots with respect to
  the Cartan subalgebra $\mathfrak{h} $.  \\

\noi {\bf Notation:} Let $\fg$ be any semisimple Lie algebra over $\mathbb R$ with Cartan subalgebra $\fh $ of $\fg _\cc$.  If $\al $ is a root, we denote a root space $(\fg _\cc )_\al $ by $ \fg _{\cc , \al } $ for convenience of notation.

\begin{rem} Any subalgebra E containing $\mathfrak h$ is of the form $E = \mathfrak{h} \oplus
  \bigoplus _{\alpha \in A} \fg _{\cc , \al }$. Furthermore, since
  $[\fg _{\cc , \al } , \fg _{\cc , \be } ] = \fg _{\cc , \al + \be } $ when
  $\alpha \neq -\beta$, subalgebras containing $\mathfrak h$ are in
  bijection with closed subsets of $\Delta$. A subset $A \subset
  \Delta$ is called \emph{closed} if A  has the property that if $\alpha , \beta \in A$ and
  $\alpha + \beta \in \Delta$, then $\alpha + \beta \in A$.
\end{rem}

\begin{Def} A closed subset $A \subset \Delta$ is called
  \emph{symmetric} if $A = -A$, and a subalgebra $\fl \subset \cts$ is called
  \emph{symmetric} if $\fl = \overline \fl$, i.e. if it is defined over the $\mathbb R$.
\end{Def}

\begin{rem} $\fl$ $\supset \fh$ is symmetric if and only if its corresponding subset $A$
  of $\Delta$ is symmetric. Note that $\fl$ may be symmetric but not be a Levi subalgebra.  
\end{rem}

\begin{lem}\label{k} Let $\fl$ by a symmetric subalgebra of $\cts$ containing
  $\mathfrak h$, which corresponds to a symmetric subset $A
  \subset \Delta$.
\begin{enumerate}
\item $[\fl,\fl] = ([\fl,\fl] \cap \mathfrak h ) \oplus \bigoplus _{\alpha \in A} \fg _{\cc , \al }$, and it is a semisimple Lie algebra.  
\item  $ \fl = Z(\fl) \oplus [\fl,\fl]$.

\end{enumerate}
\end{lem}

\begin{pf}
\begin{enumerate} 
\item  Upon the observation that $A$ is itself a root system, this is a consequence of a theorem of Serre (~\cite{hum}, p.99).  

\item The Killing form $\kappa$ is non-degenerate when restricted to
  $\mathfrak h$.  Let $\mathfrak{h} _0 $ be the orthogonal complement
  to $\mathfrak {h} ^\prime := [\fl,\fl] \cap \mathfrak h $ with inner
  product given by the Killing form.  If $h \in \mathfrak{h} _0 $,
  then for any $\alpha \in \Delta$, $0 = \kappa (t_\alpha , h) =
  \alpha (h)$.  Thus, $0= \alpha (h) X_\alpha = [h , X_\alpha ] $, and
  $h \in Z(\fl)$ (since $\mathfrak h $ is abelian).  On the other hand,
  $\fl= \mathfrak{h} _0 \oplus [\fl,\fl]$. So if $x = x_0 + x^\prime \in
  Z(\fl)$, then since $x_0 \in Z(\fl)$, $x^\prime \in Z(\fl)$.  However, by
  part 1, $[\fl,\fl]$ is semisimple, which implies that $x^\prime = 0$.
  Therefore $Z(\fl) = \mathfrak{h} _0$. $\square$
  \end{enumerate}
  \end{pf}

\subsection{Quotients by Cartan Subgroups}\label{20junea}
\noi In this subsection only we consider the case when $K =C$.

\subsubsection{Generalized Complex Structures on G/C}

\begin{prop}
 The subalgebras $E$ of $\cts$ satisfying $E +
\overline{E} = \cts$ and $\mathfrak{h} \subset E$ are exactly the
parabolic subalgebras.
\end{prop}

\begin{pf} Since $E$ is a subalgebra containing $\fh$, $E $ corresponds to some closed subset $A \subset \De$.  
G is compact, so $\sigma$ maps $\fg _{\cc , \al } $
  isomorphically into $\fg _{\cc , -\al }$. Therefore  $E
  + \overline{E} = \cts $ implies that  $ A \cup -A = \Delta$.  Subalgebras of this form
  are precisely the parabolic subalgebras (~\cite{bou} p.174).  $\square$
\end{pf}

\noi For a parabolic subalgebra $E$ corresponding to a subset $A \subset \Delta$, we will let $\fl = \fl _E$ denote the Levi factor of $E$, i.e. the subalgebra corresponding to roots $A \cap -A$.  Let $E$ be a parabolic subalgebra.  To the end of showing that $\varepsilon$ is exact for any pair
$(E,\varepsilon)$ of Prop~\ref{b}, we put forth the following proposition.

\begin{prop}\label{phi} 
 Fix a parabolic subalgebra $E = \mathfrak{h} \oplus
  \bigoplus _{\alpha \in A} \fg _{\cc , \al }$ of $\cts$, and let $V =
  V_E$ be the vector space $V :
  = \{ \varepsilon \in \wedge ^2 E^* \; | \; d_E \varepsilon = 0 \; \;
  and \; \; {\tep }_{|\mathfrak{h}} = 0 \} $. Choose a base $\Delta _0$ for the roots $ \Delta $ such that the corresponding system of positive roots $\Delta_0 ^+ $ lies in $A$ and a basis $\{X_\alpha \} _{\alpha \in \Delta ^+} \cup \{X_{-\alpha } = \sigma (X_\alpha ) \} _{\alpha \in \Delta ^+} $ for $\oplus _{\alpha \in \Delta } \fg _{\cc , \al } $.  Let $\{ X_\alpha ^*  \in \fg _{\cc , \al }\} _{\alpha \in \Delta}$ be a dual basis. 
\begin{enumerate} 
\item Every $\varepsilon \in V$ is of the form  \[ \varepsilon = \sum
  _{\alpha \in A \cap (-A) } c_{\alpha}X_\alpha ^* \wedge X_{-\alpha}
  ^*  \] for some constants $c_\alpha $, where $c_{- \alpha } =
  -c_{\alpha}$.
\item Any $\vep \in V_E$ is exact.  In fact, any $\vep \in V_E$ is of the form $\tilde \phi \circ [ \;  , \; ]$, where $\phi \in ([\mathfrak{l} , \mathfrak{l} ] \cap \mathfrak{h}
  )^* $, and $\phi$ is extended by zero to a linear functional $\tilde \phi$ on $\cts =
  ([\mathfrak{l} , \mathfrak{l} ] \cap \mathfrak{h} ) \oplus
  Z(\mathfrak{l} ) \oplus \bigoplus _{\alpha \in \Delta} \fg _{\cc , \al }
  $.  Thus $V$ is parameterized by $([\fl , \fl ] \cap \fh )^*$.


\end{enumerate}

\end{prop}

\begin{pf}  \begin{enumerate}
\item For any $h \in \mathfrak{h}$, $X \in \fg _{\cc , \al }$, $Y \in
  \fg _{\cc , \be }$, and $\vep \in V$,
 \begin{eqnarray*} 0 &=&  d_E (h,X,Y)  = \varepsilon (X, [Y,h] ) +
  \varepsilon (Y,[h,X] ) \\
&=& \varepsilon (X,- \beta (h ) Y) + \varepsilon
  (Y , \alpha (h) X) = (\beta (h) - \alpha (h) ) \varepsilon (X,Y). 
\end{eqnarray*}
This implies that for $X \in \fg _{\cc , \al } $ and $Y \in \fg _{\cc , \be }$, one has
  $\varepsilon (X,Y) = 0$ unless $\alpha = -\beta $.  Therefore, if
  one chooses a base $\Delta _0$ for the roots $ \Delta $ such that
  $\Delta_0 ^+ \subset A $ and a basis $\{X_\alpha \} _{\alpha \in
  \Delta ^+} \cup \{X_{-\alpha } = \sigma (X_\alpha ) \} _{\alpha \in
  \Delta ^+} $ for $\oplus _{\alpha \in \Delta } \fg _{\cc , \al } $ along
  with a dual basis $\{ X_\alpha ^*  \in \fg _{\cc , \al }\} _{\alpha \in
  \Delta}$, \[ \varepsilon = \sum _{\alpha \in A \cap (-A) }
  c_{\alpha}X_\alpha ^* \wedge X_{-\alpha} ^*  \] for some constants
  $c_\alpha $, where $c_\alpha = -c_{-\alpha}$. \\

\item By Lemma~\ref{k}, $[\fl , \fl]$ is semisimple.  By part (1) of this proposition, $\vep$ only depends on its restriction to $[ \fl , \fl ]$, so $\vep _{|[\fl , \fl ] \times [\fl , \fl]} = \phi \circ [ \; , \;]$ by Remark~\ref{h}.  Since $\vep$ is of the form described in part (1), we can assume that $\phi \in ([\fl , \fl ] \cap \fh )^*$ and has been extended to be $0$ on all of the root spaces $\fg _{\cc , \al } $, $\al \in A \cap -A$.  We can also extend $\phi$ by zero to $\tilde{\phi}$ on $\cts =
  ([\mathfrak{l} , \mathfrak{l} ] \cap \mathfrak{h} ) \oplus
  Z(\mathfrak{l} ) \oplus \bigoplus _{\alpha \in \Delta} \fg _{\cc , \al }
  $.  Then on $E$, $\vep = \tilde \phi \circ [ \; , \;]$.  Clearly $\phi \in ([\fl , \fl ] \cap \fh )^*$ can be chosen freely, so that $V_E = ([\fl , \fl ] \cap \fh )^*$.  $\square$ 

\end{enumerate}

\end{pf}

\noindent For each $\alpha \in \Delta $, let $\ac$ denote the
associated coroot.  That is, if $\alpha = \kappa (t_\alpha , - )$,
where $\kappa $ is the Killing form, then $\ac = 2t_\alpha /
(\alpha , \alpha ) $. For a symmetric subalgebra $\fl \subset \cts$ that contains $\fh$, we denote its root system by $\De (\fl)$.\\

\begin{prop}\label{19junea} Let $\varepsilon$ and $ \phi $ be as in Proposition~\ref{phi}.
  Then Condition 5 of Corollary~\ref{b} is satisfied if and only if
  either of the following equivalent conditions holds:
\begin{enumerate}
\item $Re(c_\alpha ) \neq 0 $ for all $ \alpha \in \Delta(\mathfrak{l}
) $.
\item $Re (\phi (\ac)) \neq 0$ for every coroot $\ac$ with $\alpha
  \in \Delta(\mathfrak{l} )$.
\end{enumerate}
\end{prop} 

\begin{pf}

1.  Condition 5 of Cor~\ref{b} is satisfied exactly when
$\varepsilon(X_\alpha , \sigma(X_\alpha) ) -
\overline{\varepsilon(\sigma (X_\alpha ) , X_\alpha ) } \neq 0 $ for all $\alpha \in \Delta (\fl)$, which happens if and only if $0 \neq
Re(\varepsilon (X_\alpha , \sigma (X\alpha) ) = Re(c_\alpha  ) $ for all $\alpha \in \Delta  (\mathfrak{l} ) $. \\

2. Let $u_\alpha = [X_\alpha , \sigma (X_\alpha ) ]$ for $\alpha \in \Delta ^+ $.  There exists some constant c such that $u_\alpha = cY_\alpha$, where $Y_\alpha$ is the unique vector in $ \fg _{\cc , -\al } $ such that $[X_\alpha , Y_\alpha ] = \ac $.\\

\noi A quick computation shows that c is a real number.  Since condition 5
of Corollary~\ref{b} is satisfied precicely when $Re (\phi (u_\alpha )) \neq
0 $ for all $\alpha \in \Delta (\mathfrak{l} ) $ ($\varepsilon =
\phi \circ [ , ] $), condition 5
of Corollary~\ref{b} is equivalently satisfied if and only if $Re (\phi (\ac )) \neq
0 $ for all $\alpha \in \Delta (\mathfrak{l} ) $ ($\varepsilon =
\phi \circ [ , ] )$. $\square$

\end{pf}

\begin{rem} Since $t_{\alpha + \beta } = t_\alpha + t_\beta $, we have $\ac +
  \check{\beta} = \frac{1}{ (\alpha +\beta , \alpha + \beta )}
  ((\alpha , \alpha )\ac + (\beta , \beta )\bc )$.  The set of $\phi \in
  ([\mathfrak{l} , \mathfrak{l} ] \cap \mathfrak{h} )^*$ satisfying
  the condition of Proposition~\ref{19junea} is guaranteed to be nonempty; letting 
  $c_\alpha \in \mathbb{R} ^+ $ for all $\alpha \in \Delta _0 (\mathfrak l) := \De (\fl) \cap \De _0$ provides a $\phi \in
  ([\mathfrak{l} , \mathfrak{l} ] \cap \mathfrak{h} )^*$ satisfying
  the condition of Proposition~\ref{19junea}.  

\end{rem}

\begin{cor}\label{26nov5} 

The equivariant generalized complex structures on G/C are
 parameterized by pairs $(E, \phi )$ where E is a parabolic
 subalgebra of $\cts$ and $\phi \in ([\fl ,\fl ] \cap \mathfrak{h} )^*$
 is such that $Re(\phi(\ac)) \neq 0$ for all coroots $\ac$
 with $\alpha \in \Delta (\mathfrak l ) $.
 Moreover, if we fix a Borel
 subalgebra $\mathfrak b $ containing $\fh$, then the equivariant generalized complex structures, up to conjugacy by an automorphism of $\cts$, are parameterized by  pairs $(S,\phi )$ where S varies
 over subsets of simple roots for  $\mathfrak{b} $ and $\phi$ is as above.  
\end{cor}

\begin{pf}
The result is now direct consequences of the previous results.  The only
observation that must be made is that any  parabolic subalgebra is conjugate by
 an automorphism of $\cts$ to some parabolic subalgebra containing $\mathfrak
 b$. This proves the final assertion. $\square$
\end{pf}

\subsubsection{Real Dirac Structures}

\begin{lem}\label{psi} There is a bijection between equivariant real Dirac
  structures on G/C and pairs $(E , \psi ) $ of a subalgebra E
  containing $\mathfrak{c} $ and $\psi \in ([E,E] \cap \mathfrak{c}
  )^*$ or, equivalently, pairs   $(\mathfrak{l} , \phi ) $ of a symmetric
  subalgebra $\mathfrak{l} = E_\cc $ containing $\mathfrak{h} $ and $\phi \in
  ([\mathfrak{l},\mathfrak{l}] \cap \mathfrak{h} )^*$ such that $\phi
  ([E,E] )\subset \mathbb{R} $.

\end{lem}

\begin{pf}
 Subalgebras E containing $
\mathfrak{c} $ are in bijection with subalgebras $\mathfrak{l} $
containing $\mathfrak{h} = \fc _\cc $ such that $\overline{\mathfrak{l}}
= \mathfrak{l} $ (by sending E to $\mathfrak{l} = E_\mathbb{C} $).
These are the symmetric subalgebras containing $\mathfrak h$ and correspond to subsets $S
\subset \Delta$ which are themselves root systems.  \\

\noi Any $\varepsilon \in \wedge ^2 E^*$ may be extended $\mathbb{C}$-linearly to $\varepsilon
_\mathbb{C} \in \wedge ^2 \mathfrak{l} ^* $, and $L(\mathfrak{l} ,
\varepsilon _\mathbb{C} ) $ is a complex Dirac structure.
The proof of Proposition~\ref{phi} is still valid, even though $\fl$ is symmetric but not necessarily a Levi subalgebra.  The result is still that $\varepsilon _\mathbb{C} = \tilde{\phi} \circ
[\; , \; ] $, and $\phi ([E,E] ) \subset
\mathbb{R} $. We may write $([E,E] \cap \mathfrak{c} )_\mathbb{C} = [E,E]_\mathbb{C} \cap \mathfrak{k} _\mathbb{C} = [E_\mathbb{C} , E_\mathbb{C} ] \cap \mathfrak{k} _\mathbb{C} = [\mathfrak{l} , \mathfrak{l} ] \cap \mathfrak{h} $.  Thus, if $\psi = \phi _{|[E,E]\cap \mathfrak{c} }$, $\phi = \psi _\mathbb{C} $, by which we mean the $\mathbb{C}$-linear extension of $\psi$.  Additionally, $\overline{Z(\mathfrak{l} ) \oplus \oplus _{\alpha \in \Delta } \fg _{\cc , \al } } =  Z(\mathfrak{l} ) \oplus \oplus _{\alpha \in \Delta } \fg _{\cc , \al } $, whence $ Z(\mathfrak{l} ) \oplus \oplus _{\alpha \in \Delta } \fg _{\cc , \al } = F_\mathbb{C}$ for some $F\subset \mathfrak{g} $.  Namely, $F = Z(\mathfrak{c} ) \oplus \oplus_{\alpha \in\Delta ^+}U_\alpha$, where $U_\alpha = \mathbb{R}$-$span (X_\alpha + \sigma (X_\alpha ) , iX_\alpha - i\sigma (X_\alpha ) )$.  Clearly $[\mathfrak{c} , F] \subset F$.  We may extend $\psi$ by 0 on F to $\tilde{\psi} \in \mathfrak{g} ^* $.  Then $\tilde{\phi} = \tilde{\psi} _\mathbb{C} $ and $\varepsilon = \tilde{\psi} \circ [\; , \; ] $.  A Dirac structure, therefore, gives a pair $(\mathfrak{l} , \psi ) $
of a symmetric subalgebra $\mathfrak{l} $ containing $\mathfrak{h} $ and
$\psi \in ([E,E] \cap \mathfrak{c} )^*$.Conversely, given such a pair
$(\mathfrak{l} , \psi )$, we may extend $\psi$ to $\tilde{\psi}$ as
before and let E be such that $E_\mathbb{C} = \mathfrak{l} $.  This
gives $L(E, \tilde{\psi} \circ [\; ,\;] ) $ which obviously
corresponds to a Dirac structure. $\square $ 
\end{pf}

\begin{rem} $\phi
  ([E,E] )\subset \mathbb{R} $ is equivalent to saying that $\phi (\ac
  ) = c_\alpha $ is purely imaginary for all $\alpha \in \Delta
  (\mathfrak l) $. 

\end{rem}

\subsection{Quotients by Connected Maximal Rank Subgroups}
\noindent Let G be a compact, semisimple Lie Group and K be a subgroup containing a Cartan
subgroup $C$ of G as before.  Let $\fc = Lie (C)$ and $\fh = \fc _\cc$.  The subalgebra
$\mathfrak{k } _\mathbb C $ corresponds to some root system $\Delta (\fk _\cc) \subset
\Delta$, which will be fixed throughout this section.  For a
subalgebra E containing $\mathfrak k _\mathbb C $, denote $\mathfrak l := E
\cap \overline E $, and let $A = A_\mathfrak l \subset \Delta$ be the root system of $\mathfrak l$.  The results of Subsection~\ref{20junea} easily generalize to any $K \supset C$. 

\begin{prop}\label{5augb} Let $G$ be compact and $K \supset C$ a closed, connected subgroup. 
\begin{enumerate}
\item Equivariant real Dirac structures on G/K correspond bijectively
  to pairs $(E,\psi ) $, where E is a subalgebra containing
  $\mathfrak{k} $ and $\psi \in  Ann_{([E,E] \cap \mathfrak{c})^* }
  ([\mathfrak{k} , \mathfrak{k} ] \cap \mathfrak{c} ) $.
  Equivalently, these Dirac structures can be described by pairs
  $(\mathfrak{l} , \phi ) $, where $\mathfrak{l} $ is a symmetric
  subalgebra containing $\mathfrak{k} _\mathbb{C} $ and $ \phi \in
  Ann_{([\mathfrak{l} ,\mathfrak{l}] \cap \mathfrak{h})^*}
  ([\mathfrak{k}_\mathbb{C} ,\mathfrak{k}_\mathbb{C}] \cap
  \mathfrak{h} ) $ such that $\phi ([\mathfrak k, \mathfrak k ] \cap
  \mathfrak c ) \subset \mathbb R $.
\item The set of invariant generalized complex structures on $G/K$
  correspond bijectively to pairs $(E,\phi ) $ where E is a parabolic subalgebra
  containing $\mathfrak k _ \mathbb C $.  $\phi \in
  Ann_{([\mathfrak{l} ,\mathfrak{l}] \cap \mathfrak{h})^*}
  ([\mathfrak{k}_\mathbb{C} ,\mathfrak{k}_\mathbb{C}] \cap
  \mathfrak{h} ) $, and $Re(\phi (\ac)) \neq 0 $ for all coroots
  $\ac$ with $\alpha \in A \setminus \Delta (\fk _\cc) $.

\end{enumerate}
\end{prop} 

\begin{pf} 

\begin{enumerate}
\item Dirac pairs $(E,\vep)$ are the same as for $G/C$ except that we require $\fk \subset E$ and $\tep (\fk)=0$. 
Any such $\vep$ is, therefore, of the form $\vep = \tilde \phi \circ [ \; , \; ]$ as in Corollary~\ref{phi}.


Letting $ Ann_{([\mathfrak{l} ,\mathfrak{l}] \cap
  \mathfrak{h})^*} ([\mathfrak{k}_\mathbb{C} ,\mathfrak{k}_\mathbb{C}]
\cap \mathfrak{h} ): = \{ \alpha \in  ([\mathfrak{l} ,\mathfrak{l}]
\cap \mathfrak{h})^* \; | \; \alpha ([\mathfrak{k}_\mathbb{C}
  ,\mathfrak{k}_\mathbb{C}] \cap \mathfrak{h} ) = 0 \} $, \\
  $\phi \in Ann_{([\mathfrak{l} ,\mathfrak{l}] \cap
  \mathfrak{h})^*} ([\mathfrak{k}_\mathbb{C} ,\mathfrak{k}_\mathbb{C}]
\cap \mathfrak{h} )$ 
$\iff$ $\phi
\in Ann_{([E,E] \cap \mathfrak{c})^* } ([\mathfrak{k} , \mathfrak{k} ]
\cap \mathfrak{c} ) $ $\iff$ $\phi
\in Ann_{([E,E] \cap \mathfrak{c})^* } ([{E} , {\mathfrak k} ]
\cap \mathfrak{c} ) $ $\iff$ $\tep (\fk)= 0$.  

\item Any $L(E,\varepsilon ) \subset \cgts$ which provides an
  invariant Dirac structure on $G/K$ also provides one on $G/C$, and $\tep $
  vanishes on $\mathfrak k $ if and only if  $\phi \in
  Ann_{([\mathfrak{l} ,\mathfrak{l}] \cap \mathfrak{h})^*}
  ([\mathfrak{k}_\mathbb{C} ,\mathfrak{k}_\mathbb{C}] \cap
  \mathfrak{h} ) $. Finally, condition 5 of Corollary~\ref{b} is met if and
  only if $Re(\phi (\ac)) \neq 0 $ for all $\alpha \in \Delta(\mathfrak
  l ) \setminus \Delta (\fk _\cc) $. $\square$

\end{enumerate}

\end{pf}

\subsection{B-Transforms}\label{5augc}

\begin{lem}\label{26nov4} Every real equivariant Dirac structure on G/K is the B-transform of some
  $L(E,0) \subset$\gts.  Therefore the equivalence class, under B-transformations, of invariant real Dirac structures on G/K is parameterized by subalgebras E of $\mathfrak{g} $ containing $\mathfrak{k} $

\end{lem}

\begin{pf}
\noindent For a Dirac structure given by $L(E,\varepsilon ) \subset $
\gts we've already seen that $\varepsilon = i^* B $ for $B =
\tilde{\psi} \circ [\; ,\; ]$, where $i: E\hookrightarrow
\mathfrak{g} $ is inclusion.  Thus, $e^B L(E,0) = L(E,i^* B ) =
L(E,\varepsilon) $, and any invariant Dirac structure is equivalent
via a B-transformation to some $L(E,0) $. $\square$
\end{pf}

\begin{prop} Let $G$ be compact and $K$ contain a Cartan subgroup.  The following data are equivalent:
\begin{enumerate}
\item $G$-invariant \gcstrs on G/K up to B-transform.
\item Triples consisting of a connected subgroup $H$ of $G$ containing $K$, a G-invariant complex structure on $G/H$, and an $H$-invariant symplectic structure on $H/K$.

\end{enumerate}

\end{prop}

\begin{pf}
We first show that a choice of a parabolic subalgebra $E$ containing $\fk _\cc$ is the same as a choice of a closed, connected subgroup $H$ containing $K$ and a $G$-invariant complex structure on $G/H$.  \\

\noi Let us first show that for any Lie subalgebra $\mathfrak H \subset \fg$ containing $\mathfrak k$, $\fH$ is its own normalizer in $\cts$. Because the normalizer $n_\cts (\fH _\cc) \supset \fH _\cc$ is a direct sum of $\fh = \fk _\cc$ and some root spaces, it is not difficult to see that $n_\fg (\fH) _\cc = n_\cts (\fH _\cc) = \fH _\cc$ and that therefore $n _\fg (\fH) = \fH$.  Now given such a Lie subalgebra $\fH$, the identity component of $N_G (\fH)$ is a closed, connected subgroup with Lie algebra $N_\fg (\fH) = \fH$.  Therefore, any Lie subgroup of $\fg$ containing $\fk$ is the Lie algebra of a closed connected subgroup $H$ of $G$.  \\

\noi Given a parabolic subalgebra $E \supset \fk _\cc$, we find that $E \cap \overline E = \fH _\cc $ for some Lie subalgebra $\fk \subset \fH \subset \fg$.  We just showed that $\fH = Lie \, H$ for some closed, connected Lie subgroup $H$ of $G$.  Then $E$ defines a $G$-invariant complex structure on $G/H$.  In the other direction, given $H$ and a $G$-invariant complex structure on $G/H$, this is simply a subalgebra $E \subset \cts$ such that $E \cap \overline E = \fH _\cc$.  It is clear that these constructions are inverses of each other.  This shows that a choice of a parabolic subalgebra containing $\fk _\cc$ is the same as the choice of a complex structure on the quotient of $G$ by a closed, connected subgroup $H \supset K$.\\

\noi It remains to show that once we have chosen a parabolic subalgebra $E \supset \fk _\cc$, a GC-pair $(E, \vep)$ is the same as a symplectic structure on $H/K$, with $H$ constructed from $E$ as in the previous paragraph.  By Lemma~\ref{26nov4}, up to $B$-transform, every GC-pair is of the form $(E, i \tilde \vphi \circ [ \; , \; ] ) $, where $ \vphi \in ([\fH  , \fH]\cap \fc ) ^*$ is extended $\cc$-linearly to the complexification, then extended by $0$ on root spaces to $\vep = \tilde \phi \in ([E \cap \overline E , E \cap \overline E] \cap \fh )^*$, where $ \vep$ is nondegenerate on $\fH / \fk$. Thus the data for GC-pairs is $\vphi \in  Ann _{ ([\fH  , \fH]\cap \fc ) ^* } ([\fk , \fk ]_\cc \cap \fc )$ such that $\vep$ (as defined above) is nondegenerate on $ \fH / \fk$.  \\

\noi Since $H$ is a closed subgroup of the compact group $G$, $H \supset K$ is itself a compact (hence reductive) group with $K$ containing a Cartan subgroup.  So Remark~\ref{26nov5} and Proposition~\ref{5augb} apply to generalized complex structures on $H/K$.  If we restrict our attention to the symplectic structures, it is clear that these are given by $\vphi \in   Ann _{ ([\fH  , \fH]\cap \fc ) ^* } ([\fk , \fk ]_\cc \cap \fc )$  such that $\vep = \tilde \vphi \circ [ \; , \; ]$ is nondegenerate on $ \fH / \fk$. 
Therefore the choice of $\vep$ is the same as a choice of an $H$-invariant symplectic structure on $H/K$. $\square$



\end{pf}

\subsection{Moduli of Complex Dirac structures on G/K}\label{2Auga}

\noindent We fix some notation.  Let $\gcs$ denote the set of
generalized complex structures on G/K, $K\supset C$.  Let $\cds$ denote the set of
complex Dirac structures, and let $\ds$ denote the set of real Dirac
structures.  We denote by $\ggcs$, $\gcds$, and $\gds$ the ones that
are G-invariant.\\

\noindent The goal of this section is to explain the geometric
structure of the space $\gcds$.  We do this by expressing $\gcds$ as a disjoint union of Euclidean spaces with closure relations.  We know that $\gds $ embeds into $\gcds$ by complexification and that $\ggcs \subset \gcds$.
Additionally, the set $\mathcal S :=$ $\{$subalgebras of $\cts$
containing $\mathfrak k _\cc$ $\}$ embeds into $\gcds$ by sending $E
\mapsto L(E,0)$.\\

\noindent For abstract reasons, we know that $\gcds$ is a variety.
The orthogonal group $ O = O(\mathfrak g \oplus \mathfrak g ^* ,
\langle\, , \, \rangle ) $ acts transitively on maximal isotropic subspaces. The
set of maximal isotropic subspaces $\mathcal L$ is the quotient of $ O$ by
the stabilizer of $\cts$ and so is itself a variety.  The maximal
isotropic subspaces L which contain $\mathfrak k _\cc$ and for which $[L,L]
\subset L$ form a closed algebraic set in this variety.  Therefore we
may think of $\gcds$ as a closed subvariety of $\mathcal L$, which is itself
a closed subvariety of the Grassmanian $Gr(dim(\mathfrak g ) ,\mathfrak g
\oplus \mathfrak g ^* ) $.\\

\noindent In the cases $ su_2 /\mathfrak k $ and $su_3
/\mathfrak k$ (for Cartan subalgebras $\mathfrak k$), $\gcds$ will be
described explicitly.  To describe $\gcds$ more generally is more
difficult.  However, we observe that 
\[ \gcds = \bigsqcup _{E \in \mathcal{S}  } \mathcal{O} _E \]
where $\mathcal{O}_E := \{ L(E,\vep ) \in \gcds \} $. Clearly $L(E,\vep ) \in
\gcds $ if and only if $ \vep \in V_E$ (in the notation of Proposition~\ref{phi}). Therefore $\mathcal{O} _E \simeq V_E$,
which is a complex vector space.  This is not very enlightening,
however, since it is possible that $\overline{\mathcal{O} _E } \cap
\mathcal{O} _F \neq \emptyset$ for some $E,F \in \mathfrak A $. This description will be
  enhanced by stating $dim(\orb _E )$ for each $E \in \mathfrak A$
  and all of the closure relations for the $\orb _E$'s.  \\

 \noindent A subalgebra $ E \in \mathfrak A$ corresponds to a subset
  $A \subset \Delta$.  $ E = E_0 \oplus E^\prime $ (direct sum as
  vector spaces) where $E_0$ is the
  subalgebra corresponding to the symmetric subset $A_0
  = A \cap -A$ and $E^\prime $ is the direct sum of the root spaces
  for $A \setminus A_0$.  Recall that $\Delta (\fk _\cc) \subset \Delta$ is the subset of roots corresponding to $\fk _\cc \supset \fh$.

\begin{prop} With $E = E_0 \oplus E^\prime \in \mathfrak A $ as above, 
\begin{enumerate}
\item $E^\prime$ is an ideal of $E$, and $A_0$ is a root subsystem.
\item $dim( \orb _E ) = rank A_0 $, and therefore $\orb _E \simeq \cc
  ^{rank(A_0) - rank(\Delta (\fk _\cc))}$.
\end{enumerate}
\end{prop}
\begin{pf}
\begin{enumerate} 
\item Since $A_0$ is the intersection of two closed subsets, it is closed.  The fact that reflections leave $A_0$ invariant follows from the fact that for non-proportional roots $\alpha$ and $\beta$, the $\alpha$-string through $\beta$ is unbroken. It is now easily verified that $A_0$ satisfies all of the root system axioms.  \\

\noi To check that $E^\prime$ is an ideal in E, let $\alpha \in A_0$ and
  $\beta \in A \setminus A_0$ such that $\gamma = \alpha + \beta$ is a root.  If
  $\gamma \in A_0$, then $\beta = \gamma  -\alpha \in A_0$ (since $A_0$
  is symmetric and closed) which is a contradiction. Therefore $[E_0 ,
  E^\prime ] \subset E^\prime$.  Now suppose that $\alpha , \beta \in A \setminus A_0$ and
  $\gamma = \alpha + \beta \in A$. If $\gamma \in A_0$, then $-\gamma
  \in A$ and $-\beta = \alpha - \gamma \in A$ which contradicts the
  fact that $\beta \in A\setminus A_0$.  This proves that $[E_0 , E^\prime]
  \subset E^\prime$ and $[E^\prime , E^\prime ] \subset E^\prime$.

\item This follows from Lemmas~\ref{k} and \ref{phi}. $\square$
\end{enumerate}
\end{pf}

\begin{prop}\label{l} $\,$ 
\begin{enumerate}
\item There is a continuous action of $W_\fk := Ann_{\fh ^*} ([\fk _\cc , \fk _\cc] \cap \fh)$ on $\gcds$, and the orbits of $W_\fk$ are all
  $\orb _E$ for $E \in \mathfrak A$. The action of $W_\fk$ is $\phi : L \mapsto e^{ \phi \circ [ \; , \;]} L$.  These are complex
 B-transforms, i.e. the 2-form is allowed to be complex.  
\item $\overline{\orb _E} \cap \orb _F \neq \emptyset$ if and only if
  $\overline{\orb _F } \subset \overline{\orb _E} $.
\end{enumerate}

\end{prop}

\begin{pf}
\begin{enumerate}
\item Any $\vep \in V_E = \{ \om \in \wedge ^2 E^* \st d_E \om =0 \; and \; \om _\sharp (\fk ) = 0 \} $ is of the form $\phi \circ
  [ \, , \, ] $ for some $\phi \in Ann_ {(\mathfrak h \cap [E_0 , E_0] )^*} ([\fk _\cc , \fk _\cc ] \cap \fh 
  )$, and $\phi$ may be extended to $\tilde \phi \in Ann _{ \mathfrak h ^* } ([\fk _\cc , \fk _\cc ] \cap \fh)  $
  as in Lemmas~\ref{k}, \ref{phi}. Thus, any $L(E,\vep )$ is the B-transform of $E
  = L(E,0)$.  B-transforms do not change the subalgebra E, so orbits
  are exactly all $\orb _E$'s.
\item If $x \in \overline{\orb _E} \cap \orb _F $.  There is a sequence $s_n \in \orb
  _E$ converging to $x$.  For any $g\in G$, $gs_n$ is a sequence in $\orb _E$ converging
  to $gx$.  Thus, $gx \in \overline{\orb _E} $ for all $g \in W_\fk$, and $\orb _F \subset \overline{\orb _E} $ because $\orb _F$
  is an orbit for $W_\fk$.

\end{enumerate}
\end{pf}

\begin{prop}\label{m} Let $G$ be compact and $K$ contain a Cartan subgroup.  If $E = E_0 \oplus E^\prime$ as above, then $L(F,\om) \in \cc \mathcal D ^G _{G/K}$ lies in $
  \overline{\orb _E }$ if and only if $F_0$ is a Levi subalgebra of
  $E_0$ and $F^\prime = E^\prime$. Moreover, $\orb _F \cap
  \overline{\orb _E} \neq \emptyset $ if and only if $F_0$ is a Levi subalgebra of
  $E_0$ and $F^\prime = E^\prime$.
\end{prop}

\begin{pf} First we show that if $x \in \gcds$ is not an isolated point, then $x \in
  \overline{\orb _E} $ for some $E \in \mathcal S$, and there is a
  sequence $L(E,\vep _n) $ converging to $x$.  If $x$ is not an isolated point, there is a sequence in $\gcds$
  converging to $x$.  But $\gcds$ is a union of finitely many orbits
  $\{\orb _E\}_{E \in \mathcal S} $.  Since there are finitely many
  sublagebras E in $\mathcal S$, there is a subsequence which lies in
  exactly one of the orbits.\\
  
\noi Suppose that $\orb _F \cap \overline{\orb_E} \neq \emptyset$.
  By part 1 of this lemma and Proposition~\ref{l}, it is enough to
  assume that there is a sequence $L(E,\vep _n ) $ converging to
  $L(F,0)$ and then to show that $F_0$ is a Levi subalgebra of $E_0$
  and that $F^\prime = E^\prime$. Fix a basis $\{ X_\alpha\}_{\alpha
  \in \Delta}$ of $\bigoplus _{\alpha \in \Delta} \fg _{\cc , \al } $ such
  that $[X_\alpha , X_{-\alpha} ] = t_\alpha $ for all $ \alpha \in
  \Delta $. Let $\{X_\alpha ^* \}$ be a dual basis.\\ 

\noindent For a point $L(E,\vep ) \in \gcds$, write $\vep$ in coordinates as $\vep = \sum _{\alpha
  \in A_0} c^\alpha X_\alpha ^* \wedge X_{-\alpha} ^* $.  Let $\pi = \sum _{\alpha \in T}
  \frac{-1}{c^\alpha} X_\alpha \wedge X_{-\alpha} $, where $T= \{ \al \in A_0 \st c^\al \neq 0 \}$.  We have, in fact, expressed $\vep$ as  $\vep \in  \wedge ^ 2 \cts ^* $, which we restrict to $E \times E$.  It therefore makes sense to consider $Im(\tep) \subset \cts ^*$, which is just the span of the dual vectors $\{X^* _\al \} _{\al \in T}$.  We immediately observe that $L(E,\vep ) =
  L(\pi , Ann(E) \oplus Im(\tep) )$
We have a sequence $L(E , \vep _n) \rightarrow L(F,0)$ with $\vep _n = \phi _n \circ [ \, ,
  \, ] $ for a sequence $\phi _n \in ([E_0 , E_0] \cap \mathfrak h )
  ^*$. We may also think of $\vep _n $ as a system $ (c_n ^\alpha )_{\alpha ,
  n}$ with $\al \in A_0$ and $n \in \mathbb N$. For any $\alpha \in A_0$, the sequence $|(c^\alpha _n)|$ may
  converge to infinity.  If it does not, we can choose a bounded
  subsequence and therefore a convergent
  subsequence of $c_n ^\alpha$. Since $A_0$ is finite, there is
  a subsequence $\vep _{n_i} $ for which each  $c^\alpha _{n_i}$ is
  either convergent or goes to infinity.  For ease of notation, assume
  that the original sequence has this property.  In order to get $L(F,0)$ as the limit, for each $\alpha
  \in A_0$, it must be true that either $c^\alpha _n \rightarrow 0$ or
  $|c^\alpha _n| \rightarrow \infty$.  If $c^\al _n \rightarrow k \neq 0$, then since $X_\al + c^\al _n X {-\al} ^*$, $X_{-\al} - c^\al _n X_\al ^* \in L(E , \vep _n)$, we would have that $X _\al + kX_{-\al} ^*$, $X_{-\al} - k X_\al ^* \in L(F, 0)$, which is impossible since $L(F,0)$ is isotropic.  \\

\noindent We define $S : = \{ \alpha \in A_0 \; | \; \phi _n (\ac )
\rightarrow 0 \} = \{ \alpha \in A_0 \; | \; \phi _n (t_\alpha )
\rightarrow 0 \}$.  Note that $\Delta (\fk _\cc)$ is contained in $S$ because $\phi _n $ vanishes on all coroots $\check \al$ for $\al \in \Delta (\fk _\cc)$ and for all $n$.  Since $t_{\alpha + \beta} = t_\alpha + t_\beta $, S is closed, symmetric, and in fact a root subsystem.  If $\alpha _1 , ..., \alpha _r $ is a base for $S$,  suppose $\beta = a_1 \alpha _1 +
...+ a_r \alpha _r $ with all $a_i \in \rr$. then $t_\beta =
a_1t_{\alpha _1} +...+ a_r t_{\alpha _r} $ and $\phi _n (t_\beta)
\rightarrow 0$.  If $W = span(S) \subset span(\Delta)$, then
$W \cap \Delta = S $  It is known that for any subspace $W \subset span (\Delta)$, $W\cap \Delta$ defines a Levi subalgebra ~\cite{bou} (p.178).  Therefore S defines a
Levi subalgebra.  Note that $S= \{ \alpha \in A_0 \; | \; c^\alpha _n
\rightarrow 0 \} $ because $\phi _n (t_\alpha) = \vep _n ([X_\alpha ,
  X_{-\alpha} ])$.  We may replace $\vep _n$ with a sequence $\vep _n
  ^\prime $ such that $c_{\alpha, n} ^\prime = 0 $ for all $ n$
  whenever $\alpha \in S $ and $c_{\alpha, n} ^\prime = c_{\alpha , n} $ for all $n $ 
whenever $\alpha \in A_0 \setminus S$.  Then $L(E, \vep _n) $
  and $L(E,\vep ^\prime _n )$ both converge to $L(F,0)$.  Note that
  $L(E,\vep ^\prime _n )$ may not be in $\gcds$, but it is in the
  Grassmanian $Gr(\mathfrak g \oplus \mathfrak g ^* , dim(\mathfrak g)
  )$ of which $\gcds$ is a subvariety.  \\

\noindent We readily see that $Im((\vep ^\prime _n)_\sharp ) = span(X_\alpha ^* )_{\alpha \in A_0
  - S} $ for all n, and $Ann(E) = span\{ X_\alpha ^*\}_{\alpha \in \Delta - A}
$.  If $\pi _n ^\prime  = \sum _{\alpha \in A_0 \setminus S}
\frac{-1}{c_{\alpha , n}} X_\al \wedge X_{-\al}$, then $L(E,\vep ^\prime _n ) = L(\pi _n ^\prime, Ann(E) \oplus Im((\vep
  ^\prime _n ) _\sharp )) = L(\pi _n  ^\prime, span\{ X_\alpha ^* \; | \; \alpha \in
(\Delta \setminus A) \cup (A_0 \setminus S))$, and $\pi _n ^\prime \rightarrow 0$ as $n \rightarrow \infty$.  
The limit is $L(0, span\{ X_\alpha ^* \; | \; \alpha \in
(\Delta \setminus A) \cup (A_0 \setminus S)) = L(\mathfrak h + E^\prime \oplus
\mathfrak l _S , 0)$ where $\mathfrak l _S$ is the Levi subalgebra
with root system S.  Now $A^\prime \cup S$ is the root system for a
subalgebra F such that S corresponds to $F_0$ and $F^\prime =
E^\prime$.  That $F = \fh + E^\prime + \fl _S$ is a subalgebra follows from the fact that $E^\prime $ is an ideal in $E$.  This shows that every limit point of $\orb _E$ is of this
form.  \\

\noindent It remains to show that any subalgebra F containing $\fk _\cc$ with $F_0$
a Levi subalgebra of $E_0$ and $F^\prime = E^\prime $ is a limit point
of $\orb _E$. Suppose $F_0$ corresponds to the root
subsystem $S\subset \Delta $ (Recall $\Delta (\fk _\cc) \subset S$).  Since S represents a Levi subalgebra of
$A_0$, it is possible to 
 choose a base $\alpha _1 , ..., \alpha _r $ for the
root subsystem $S$ which extends to a base  $\alpha _1 , ..., \alpha
_r , \alpha _{r+1} , ..., \alpha _n $ of $A_0$.   By Lemma~\ref{k}, Any $\vep $ for which $d_E
\vep = 0 $ is determined by the constants $c_1 = c_{\alpha _1}, ...,
c_n = c_{\alpha _n} $.  Choose a sequence $\vep _n $ such that all
$c_i \in \rr $ and are non-negative, $c_i = 0$ for $i \leq r$ and $c_i
\rightarrow \infty $ for $i > r $. The $t_i = t_{\alpha _i } $ form a
basis of $\mathfrak h \cap [E_0 , E_0] $ with dual basis given by $t_i
^*$,so $\vep = (c_1 , ..., c_n ) = (\sum c_i t_i ^* ) \circ [ \, , \, ] $.\\
 
\noi Then $Im(({\vep _n}) _\sharp ) = span\{ X_\alpha ^* \; | \; \alpha \text{ cannot be
 expressed purely as } n_1\alpha _1 + ...+ n_r \alpha _r \; n_i\geq 0
 \} = span \{ X_\al ^* \st \al \in A_0 \setminus S \} $ for all $n$. Then because $L(E,\vep _n ) = L(\pi _n , Im(({\vep _n }) _\sharp )\oplus Ann(E) ) $
 and $\pi _n \rightarrow 0 $, the limit is $L(0, Im(({\vep _n )_\sharp }
 )\oplus Ann(E) ) = L(F, 0 )$.  $\square$

\end{pf}

\subsubsection{The Moduli for $SU_2 $ and $SU_3$} 
\noi We consider the quotients of $SU_2 $ and $SU_3$ by the standard Cartan subgroups and delineate explicitly the moduli of Dirac structures on these spaces. \\

\begin{lem} When $G = SU_2$ and $K$ is the standard Cartan subgroup, \\
$\gcds = \cc \mathbb P ^1 \cup$$\{$two points$\}$, and 
$\GG \mathcal C ^G _{G/K} = \cc \setminus i\mathbb R \subset \cc \subset \cc \mathbb P ^1 \subset \cc \mathbb P ^1 \cup \{two \; points \}$.
\end{lem}

\begin{pf}  The subalgebras of ${\fs \fl}_2 (\cc ) = ({\fs \fu}_2)_\cc
 $ containing the standard Cartan $\mathfrak h $ are $\mathfrak h $,
 the two Borel subalalgebras $\mathfrak b _1 , \mathfrak b_2 $, and
 all of ${\fs \fl}_2 $. Let $\alpha$ denote one root, so in the notation of Subsection~\ref{2Auga}, $\orb _{{\fs \fl}_2} = \{
 L({\fs \fl}_2 , cX^* _\alpha \wedge X^* _{-\alpha} ) \} \simeq \cc $.  As $c
 \rightarrow \infty $, we get $L(\mathfrak h , 0) $. which gives a
 copy of $\cc \mathbb P ^1 $. Proposition~\ref{m} guarantees that this
 $L(\mathfrak h , 0)$ is the only point in the closure of this orbit
 and that the Borel subalgebras are isolated points.  Thus, in this
 case, $\gcds = \cc \mathbb P ^1 \cup$$\{$two points$\}$.  Proposition~\ref{19junea} states that $\GG \mathcal C ^G _{G/K} = \cc \setminus i\mathbb R \subset \cc \subset \cc \mathbb P ^1 \subset \cc \mathbb P ^1 \cup \{two \; points \}$. $\square$
\end{pf}

\begin{lem} Let $G = SU_3$ and $K$ be the standard Cartan subgroup.  Then $\cc \DD ^G _{G/K}$ is the disjoint union of the following connected components:
\begin{enumerate}
\item There is an isomorphism $\overline{\orb
  _{{\fs \fl}_3 }} \simeq \{ [x,y] \times [u,v] \times [s,t] \in (\cc \mathbb P ^1
  )^3 \; | \; vtx + ytu -yvs = 0 \}$. The closure $\overline{\orb
  _{{\fs \fl}_3 }}$ of ${\orb
  _{{\fs \fl}_3 }}$ consists of ${\orb
  _{{\fs \fl}_3 }}$ together with $\OO _E$ for each Levi subalgebra $E$.  The generalized complex structures are $\mathcal G \mathcal C ^G _{G/K} \cap \overline  {\mathcal O _{{\mathfrak s \mathfrak l} _3}} = {\orb
  _{{\fs \fl}_3 }} \cap (\cc \setminus i \mathbb R )^3 \subset \cc ^3 \subset \cc \mathbb P ^1 $.
 \item  Let $E$ be one of the six proper parabolic subalgebras.  Then $\overline \OO _E \simeq \cc \mathbb P ^1$.  The generalized complex structures are $\GG \mathcal C ^G _{G/K} \cap \overline \OO _E - \cc \setminus i \mathbb R \subset \cc \subset \cc \subset \mathbb P ^1$.
 \item The twelve remaining subalgebras represent isolated points.
\end{enumerate}

\end{lem} 
 
 \begin{pf} 
  
\noindent For $SU_3$, $({\fs \fu}_3)_\cc = {\fs \fl}_3 (\cc ) $.  Let $\mathfrak h $
be the standard torus consisting of the diagonal matrices, and let
$X_{i,j}$ denote the matrix with 1 in the (i,j)-th entry and zeros
elsewhere.  Let $\alpha$ denote the root for which the root space
contains $X_{1,2}$, and let $\beta$ denote the root for which
$X_{2,3}$ is the root space.  Write $\vep \in Z^2 ({\fs \fl}_3 , \cc )$ as
$\vep = c_\alpha X_\alpha ^* \wedge X^* _{-\alpha} + c_\beta X_\beta
^* \wedge X^* _{-\beta} +c_{\alpha +\beta} X_{\alpha +\beta} ^* \wedge
X^* _{-\alpha - \beta}$ with $c_{\alpha + \beta } = c_\alpha + c_\beta
$.  We may think of $\orb _{{\fs \fl}_3} $ as $ \{(c_\alpha , c_\beta , c_{\alpha +
  \beta}) \in \cc ^3 \; | \; c_\alpha + c_\beta - c_{\alpha + \beta} =
0 \}$.  Embedding $\cc$ into $\cc \mathbb P ^1 = \cc \cup \infty$ gives $\cc ^3 \subset (\cc \mathbb P ^1 )^3$.\\

\noindent For a closed subset $\Phi \subset \Delta$ of roots, we denote by $E_\Phi$ the corresponding Lie subalgebra which is the sum of $\fh$ and the root spaces for the roots in $\Phi$.  We will also denote $\OO _{E _{\Phi}}$ simply by $\OO _\Phi$. We make the following identifications.\\ 

\noi $\orb_{\pm \alpha} = L(E_ {\{\pm \alpha\} }, c_\alpha X_\alpha ^* \wedge X_{-\alpha}
^* )$ corresponds to the point $[c_\alpha , 1]\times [1,0]\times [1,0] \subset (\cc
{\mathbb P} ^1 )^3$.\\
$\orb _{\pm \beta} = L(E _{\{\pm \beta\} }, c_\beta X_\beta ^* \wedge X_{-\beta}
^* ) $ corresponds to the point  $[1 , 0] \times [c_\beta,1] \times [1,0] \subset (\cc
{\mathbb P} ^1 )^3 $.\\
$\orb_{\pm (\alpha + \beta)} = L(E _{ \{\pm (\alpha +\beta )\} } , c_{\alpha
  +\beta}  X_{\alpha +\beta } ^* \wedge X_{-\alpha -\beta}
^* ) $ corresponds to the point $ [1,0]\times [1,0]\times [c_{\alpha + \beta},1] \subset (\cc
{\mathbb P} ^1 )^3$.\\
$\OO _\emptyset = L(\mathfrak h , 0 ) $ corresponds to the point $ [1,0 ] \times [1,0] \times [1,0] $.\\

\noindent This is a complete list of the Levi subalgebras in ${\fs \fl}_3
$ that contain $\fh$, all of which lie in $\overline{\orb _{{\fs \fl}_3}} $. A quick computation, using the above identifications for the Levi subalgebras, shows that we can identify $\overline{\orb
  _{{\fs \fl}_3 }}$ with $ \{ [x,y] \times [u,v] \times [s,t] \in (\cc \mathbb P ^1
  )^3 \; | \; vtx + ytu -yvs = 0 \}$.  The generalized complex structures are $\mathcal G \mathcal C ^G _{G/K} \cap \overline  {\mathcal O _{{\mathfrak s \mathfrak l} _3}} = {\orb
  _{{\fs \fl}_3 }} \cap (\cc \setminus i \mathbb R )^3 \subset \cc ^3 \subset \cc \mathbb P ^1 $ by Proposition~\ref{19junea}.\\

\noindent There are six proper parabolic subalgebras: $\pm \{\pm \alpha , \beta , \alpha
+\beta \} $, $\pm \{ \alpha , \pm \beta , \alpha + \beta \} $, and\\
$\pm \{ \alpha , -\beta , \pm (\alpha + \beta ) \} $.  For each of the
parabolics E, $\orb _ E \simeq \cc $, and the limit contains one
point.  Therefore $\overline{\orb _E} \simeq \cc \mathbb P ^1 $.
For example, if E is $\{\pm \alpha , \beta , \alpha + \beta \}$, then
$\orb _E = \{L(E, cX_\alpha ^* \wedge X_{-\alpha} ^* )\} $, which is the
same as $\cc$.  Letting $c \rightarrow \infty$ gives $L(\{\beta ,
\alpha + \beta \} ,0 )$.  For each of these parabolic subalgebras, the generlized complex structuers are $\mathcal G \mathcal C ^G _{G/K} \cap \overline {\mathcal O _E} = \cc \setminus i\mathbb R \subset \cc \subset \cc \mathbb P ^1$ by Proposition \ref{19junea}.  Since generalized complex pairs only occur for parabolic subalgebras, this provides a complete list of the generalized complex structures.  There are six subalgebras which contain root spaces for two
roots but are not Levi subalgebras.  In this way, each is in the
closure of some $\orb _E$ for E parabolic.   Each of these copies of
$\cc \mathbb P ^1 $ is a connected component of the moduli space $\gcds$.\\

\noindent There are six subalgebras which
contain only one root space, and there are six Borel
subalgebras.  These are all isolated points.  $\square$

\end{pf}

\section{Semisimple Orbits}\label{sectionorbits} 
\noi We have given a description of Dirac structures on adjoint orbits when $G$ is compact.  We now attempt to describe \gcstrs on more semisimple orbits in more general groups.  
In the case of a semisimple orbit $\mathcal O _h$ in a real semisimple Lie algebra, we
would like to understand what are GC-pairs $(E, \vep)$ for the homogeneous space $\OO _h = Int (\fg) /Z_G (h)$  (where $Int (\fg) $ is the connected Lie subgroup of $Aut (\fg)$ with Lie algebra $ad \fg$).  GC-pairs turn out to be equivalent to a pair $(A, \phi)$ of a closed subset $A$ of roots and a linear functional on $\check A \cap -\check A \subset \check \Delta$ satisfying some conditions (Theorem~\ref{il} and Corollary~\ref{2octa}).  We go on to describe such closed subsets $A$ to parabolic subalgebras.\\


\noi First we consider the case when $h$ is a regular semisimple
element, i.e. when $\mathfrak h := Z_\fg (h)$ is a Cartan
subalgebra.  The notation and formulation of statements is less burdensome for regular elements, but the proofs are essentially the same.  As we will see, the results for general semisimple orbits follow immediately once we have done the regular case. \\

\noi Throughout this section, fix a Cartan involution $\theta$ of
$\fg$.  Since any Cartan subalgebra is $Int(\fg)$-conjugate to a
$\theta$-stable one, we may assume that $\mathfrak h $ is
$\theta$-stable so that $\mathfrak h = \mathfrak t \oplus \mathfrak e
$ is the Cartan decomposition of $\mathfrak h $.  Let $x \mapsto \bar{x} $
or $\sigma$ denote conjugation in $\cts$ with respect to $\fg$.   Also
$\sigma$ will denote conjugation with respect to roots: $(\sigma
\alpha ) (h) = \overline {\alpha (\sigma h)} $ for $\alpha \in \Delta
= \Delta (\cts , \mathfrak h _\cc)$. The involution $\theta$ extends
to a $\cc$-linear map on $\cts$, also denoted by $\theta$.  Since
$\mathfrak h $ is $\theta$-stable, $\theta$ permutes the roots by
$(\theta \alpha ) (h) = \alpha (\theta h ) $.  It is the case that
$\theta = -\sigma$ on $\mathfrak h _\cc ^*$, and $\sigma_{|\mathfrak e
   ^*} = 1$, $\sigma_{|\mathfrak t
   ^*} = -1$.\\

\noi We begin with the (simple version) of the main theorem of this section.  The full version is Corollary \ref{2octa}, which addresses the case when $h$ is an arbitrary semisimple element.  

\begin{thm}\label{il} Let $h $ be a regular semisimple element in a real
  semisimple Lie algebra $\fg$ as above.  Then the equivariant generalized complex
  structures on the adjoint orbit $\mathcal O _h$ are given by pairs
  $(A , \phi ) $ where $A \subset \Delta $ is a closed subset of roots such that $A \cup \sigma A = \Delta $ and $A \cap \sigma A \subset -A$, and $\phi$ is a linear functional on $span_\cc (\check A \cap -\check{A})$ satisfying $\phi(\check \alpha ) \neq \overline{ \phi
  (\check{\sigma(\alpha)})} $ for all $ \alpha \in A \cap \sigma A$.

\end{thm}

\begin{pf} Equivariant generalized complex structures on $\mathcal O  _h =
  G/Z_G(h)$, $G= Int(\fg)$, are given by pairs GC-pairs $(E,\vep)$.  Since $\mathfrak h _\cc \subset E$, $E = \mathfrak h _\cc
  \oplus \bigoplus _{\alpha \in A \subset \Delta} (\cts)_\alpha $ for
  some closed subset A of $\Delta$.  Since $\sigma$ maps $(\cts)_\alpha$ isomorphically to $(\cts)_{\sigma \alpha}$,  Conditions (1) and (2) of
  Proposition~\ref{b} are satisfied if and only if $A \cup \sigma A =
  \Delta$.  Now let $S= A \cap -A$.  Then $\mathfrak h ' := span_\cc
  \check S \subset \mathfrak h _\cc $.  We claim that $d_E \vep = 0 $
  and $\tep (\mathfrak h _\cc ) = 0 $ if and only if $\vep = \tilde
  \phi \circ [ \; , \; ]$, where $\tilde \phi$ is an
  extension--extended to be 0 on $(\mathfrak h ')^\perp $ ($\perp$
  determined by the Killing form $\kappa$)--of some $\phi \in
  (\mathfrak h ' )^*$.  \\

\noi Clearly any such $\phi$ gives  $\vep = \tilde
  \phi \circ [ \; , \; ]$ satisfying $d_E \vep = 0 $
  and $\tep (\mathfrak h _\cc ) = 0 $.  It remains to show that any
  such $\vep$ is of this form.  \\

\noi Just as in Proposition~\ref{phi}, if $X \in (\cts)_\alpha$ and $Y\in
(\cts)_\beta$, then $\vep (X,Y) \neq 0 $ only if $\beta = -\alpha$.
Therefore if we let $E' = \mathfrak h ' \oplus \bigoplus _{\alpha \in
  S} (\cts) _\alpha$, $\vep$ is determined by $\omega := \vep _{|E'
  \otimes E' }$.  It is still true that $\omega _\sharp (\mathfrak h ' ) = 0
$ and $d_{E'} \omega = 0 $.  However, now there is the advantage that
$E'$ is semisimple, which implies that $H^2 (E' , \cc) = 0$.  Then
since $\omega \in Z^2 (E' , \cc)$, $\omega = \phi \circ [ \; , \; ]
$ for some $\phi \in (E')^*$.  However, $\phi$ must vanish on each
$(\cts)_\alpha$ because $(\cts)_\alpha = [ \mathfrak h ',
  (\cts)_\alpha ] = $.  Thus, $\phi $ is determined by $\phi
_{|\mathfrak h ' } $.  \\

\noi Finally, we claim that condition (5) of Proposition~\ref{b} are 
satisfied if and only if the following two conditions are met:\\
a.) $A \cap \sigma A \subset -A$ and \\
b.) $\phi(\check \alpha ) \neq \overline{ \phi
  (\check{\sigma(\alpha)})} \; $ for all $ \alpha \in A \cap \sigma A$.\\

\noi If a.) and b.) are true, let  $X = h + \sum _{\alpha \in A
  \cap \sigma A} X_\alpha \notin \mathfrak h _\cc$.  Since $X \notin
  \mathfrak h _\cc$, $X _\alpha \neq 0 $ for some $\alpha \in A \cap
  \sigma A \subset -A$.  Choose $Y \in (\cts)_{-\alpha }$ such that $[X,Y] =
  \check \alpha $.  This entails that $\vep (X,Y) = \phi (\check
  \alpha )$.  On the other hand, $\overline{\vep (\sigma X , \sigma
  Y)}= \overline{ \phi (\sigma (\check \alpha))} = 
\overline{ \phi (\check {\sigma ( \alpha)})} \neq \phi (\check \al) $ by b.).  So
  condition 5 of Proposition~\ref{b} is satisfied.  \\

\noi Now suppose that condition 5 of Proposition~\ref{b} is satisfied and
let $X \in (\cts)_\alpha$ for some $\alpha \in A \cap \sigma A$.  Condition 5
implies that there exists $Y \in E \cap \overline E$ such that $\phi([X,Y]) \neq \overline {\phi (\sigma [X,Y])}$.
  Sine $E = \fh \oplus \bigoplus _{\al \in A} (\cts)_\al$, we may assume that $Y \in (\cts)_{-\al}$ since $\phi$ (or more precisely $\tilde \phi$) vanishes on the root spaces.  Hence, $0\neq Y \in E \cap \overline E \cap (\cts)_{-\al}$ and so $(\cts)_{-\al} \in E$ and $-\al \in A$.  This shows that $A \cap \sigma A \subset -A$.  Also, $Y$ can be chosen such that $[X,Y] = \check \al$, so $\phi (\check \al) \neq \overline { \phi (\sigma \check \al) } = \overline { \phi ( \check {\sigma \al} )}$. $\square$


\end{pf} 

\noi There is a couple of easy examples:\\
1.) If $\mathfrak h $ is a split Cartan subalgebra, then any
generalized complex pair $(E, \vep )$ must satisfy $E = \cts$ since
$\sigma \alpha = \alpha $ for all roots $ \alpha \in \Delta$.\\
2.) If $\fg$ is compact, then this is the example described in Section~\ref{26nov9}.  \\

\noi If $h \in \fg$ is an arbitrary semisimple element, then $Z :=
Z_\cts (h)$ is a Levi subalgebra containing a Cartan subalgebra $
\mathfrak h$.  Then the subalgebra $Z$ corresponds to some subset
$\Lambda \subset \Delta (\cts , \mathfrak h )$.  The situation in this
case is almost exactly the same except that we require $\Lambda
\subset A$.  This can be stated explicitly in the following theorem.

\begin{cor}\label{2octa}  Let $h $ be a semisimple element in a real
  semisimple Lie algebra $\fg$ as above. As above, let $\Lambda
  \subset \Delta (\cts , \mathfrak h ) $ be the set of roots
  corresponding to $Z(h) \supset \mathfrak h $.   The equivariant generalized complex
  structures on the adjoint orbit $\mathcal O _h$ are given by $Z_{G}(h)$-invariant pairs
  $(A , \phi ) $, where $\Lambda \subset A \subset \Delta $ is a closed subset of roots
  and $\phi \in (span_\cc \{ \check \alpha \; | \; \alpha \in A \cap
  -A \} )^*$, satisfying:
\begin{enumerate}
\item $A \cup \sigma A = \Delta $,
\item $A \cap \sigma A \subset -A$,  
\item $\phi(\check \alpha ) \neq \overline{ \phi
  (\check{\sigma(\alpha)})}$ for all $\alpha \in A \cap \sigma A
  \setminus \Lambda$, and
\item $\phi (\check \alpha ) = 0 $ for all $\alpha \in \Lambda$.  
\end{enumerate}
\end{cor}

\begin{pf} GC-pairs $(E, \vep ) $ correspond to pairs $(A, \phi ) $ as
  in Theorem~\ref{il}.  The only additional requirements are that
  $\Lambda \subset A$ and  $\phi (\check \alpha ) = 0 $ for all
  $\alpha \in \Lambda$. $\square$
\end{pf} 

\noi In the remainder of the section we analyze conditions (1) and (2) of Theorem~\ref{2octa} and relate them to $\theta$-stable parabolic subalgebras.  

\begin{Def} Any closed subset $A \subset \Delta$ will be called a
  \emph{generalized complex subset} of roots if $A \cup \sigma A = \Delta$ and $A \cap \sigma A \subset - A$.
\end{Def}

\noi Generalized complex subsets correspond to subalgebras $E $ of $
\cts$ which occur in generalized complex pairs $(E, \vep )  $ by Corollary~\ref{2octa}.
It would now be helpful to have some description of generalized
complex subsets of roots.

\begin{lem} If A is a generalized complex subset , then $A \cup \theta
  A $ is the closed subset of roots corresponding to a parabolic
  subalgebra containing $\mathfrak h _\cc $.  
\end{lem}

\begin{pf}  A closed subset, $ \Phi$, of roots corresponds to a parabolic
  subalgebra containing $\mathfrak h _\cc $ if and only if $ \Phi \cup
  - \Phi = \Delta$ ~\cite{bou}.  Closed subsets of roots of this type
  are called parabolic.  Let $\Phi = A \cup \theta A $.  First we see
  that $\Phi \cup - \Phi = A \cup \theta A \cup -A \cup -\theta A
  \subset A \cup - \theta A = A \cup \sigma A = \Delta $.  \\

\noi Therefore it only needs to be shown that $\Phi $ is a closed
subset of $\Delta$.  If $\alpha$, $\beta$ $\in \Phi$ and $\alpha +
\beta = \gamma \in \Delta$, we must show that $\gamma \in \Phi$.
Since $A$ and $\theta A$ are closed, this reduces to the case when
$\alpha \in A $ and $\beta \in \theta A$.  In this case, if $\gamma
\in A \subset \Phi$, then there is nothing to show.  Otherwise,
$\gamma \in \sigma A$ because $A \cup \sigma A = \Delta$.  Then
$-\beta \in \sigma A$, and $\gamma \in \sigma A$ so that $\alpha =
\gamma - \beta \in \sigma A$.  Hence, $\alpha \in A \cap \sigma A
\subset -A$. Therefore $\pm \alpha \in A$.  But also $\sigma \alpha
\in A \cap \sigma A \subset -A $ implies $-\sigma \alpha \in A$,
whence $\pm \alpha \; , \pm \sigma \alpha \in A$.  Then $\alpha \in
\theta A$ and so $\alpha + \beta \in \theta A + \theta A \subset
\theta A$.  $\square$  

\end{pf} 

\noi If $A$ is a generalized complex subset of $\Delta$, then $\Phi =
A \cup \theta A $ is  parabolic and $\Phi = \Gamma \sqcup \Psi$,
where $\Gamma = \Phi \cap - \Phi$ and $\Psi = \Phi \setminus \Gamma$.
It is clear that if $\Phi$ corresponds to a parabolic subalgebra $
\mathfrak p \subset \cts$, then $\Gamma$ corresponds to a Levi factor
$\mathfrak l $ of $\mathfrak p $ and $\Psi $ corresponds to the
nilpotent radical $\mathfrak u$ of $\mathfrak p$: $\mathfrak p =
\mathfrak l \oplus \mathfrak u$.  \\

\noi We easily see that $\Gamma = (A \cup \theta A ) \cap (-A \cup
-\theta A) = (A \cap -A) \cup (A \cap -\theta A)\cup (\theta A \cap
-A) \cup (\theta A \cap -\theta A )$.  But $A \cap -\theta A = A \cap
\sigma A \subset A \cap -A$ and similarly, $\theta A \cap -A \subset A
\cap -A$.  Therefore $\Gamma = (A \cap - A ) \cup \theta(A \cap
-A)$.  Since the Levi subalgebra $\fl$ is reductive,

\[
\mathfrak l = \bigoplus _{i=1} ^n S_i  \oplus Z(\mathfrak l ), 
\]
where each $S_i$ is a simple Lie algebra and $Z(\mathfrak l ) \subset
\mathfrak h$.  Each $S_i$ corresponds to a
subset $\Gamma _i \subset \Gamma$.  Obviously $\Gamma _i = (A \cap - A
)\cap \Gamma _i  \cup \theta (A \cap - A ) \cap \Gamma _i $.  The
following lemma will demonstrate that for each $i$, either $\Gamma _i
\subset (A \cap - A )$ or $\Gamma _i \subset \theta (A \cap - A ) $.\\

\noi Obviously $\Phi$ is $\theta$-stable, and $\Gamma = (A \cap -A) \cup \theta (A \cap -A)$, so $\Gamma$ is also $\theta$-stable. Hence, $\Psi$ is also $\theta$-stable.  Of course, this means that $\mathfrak p$, $\mathfrak u$, and
$\mathfrak l $ are also all $\theta$-stable.  We know that $\Delta = A
\cup \sigma A$ so that $-\Psi \cap \Phi = \emptyset $ implies that
$-\Psi \subset \sigma A$.  Applying $\sigma$ gives
$\Psi \subset A$ because $\Psi$ is $\theta$-stable.

\begin{lem}\label{lk} Let $\fg$ be a simple complex Lie algebra containing
  Cartan subalgebra $\mathfrak h$.  If $P$ and $Q$ are two reductive
  subalgebras containing $\mathfrak h$ such that $P+Q = \fg$, then $P=
  \fg$ or $Q= \fg$. Equivalently, if $\Delta$ is an irreducible root
  system and $X, Y$ are closed symmetric subsets of $\Delta = \Delta
  (\fg , \mathfrak h)$ such that
  $X\cup Y = \Delta$, then either $X= \Delta$ or $Y = \Delta$.

\end{lem}

\begin{pf} We prove the latter statement.  The first step is to show
  that $X^c \perp Y^c$, where $X^c = \Delta \setminus X = Y \setminus
  X$.  Let $\alpha \in X$, $\beta \in X^c$.  Suppose that $\alpha
  +\beta = \gamma \in \Delta$.  If $\gamma \in X$, then $\beta =
  \gamma -\alpha $.  But $X $ is symmetric, so $\alpha \in X
  \Longrightarrow -\alpha \in X$.  Since $X$ is closed, $\beta =
  \gamma - \alpha \in (X+X) \cap \Delta \subset X$, which is a
  contradiction.  Therefore, $\Delta \cap ( X + X^c) \subset X^c$ and
  similarly,  $\Delta \cap ( Y + Y^c) \subset Y^c$.  If $\alpha \in
  X^c$ and $\beta \in Y^c$, then $\alpha +\beta \in Y + Y^c \subset
  Y^c$ and $\alpha +\beta \in X+ X^c \subset X^c$, which means that
  $\alpha + \beta \in X^c \cap Y^c = \Delta ^c = \emptyset$.
  Therefore $(X^c + Y^c) \cap \Delta = \emptyset$.  If $\alpha \in
  X^c$, $\beta \in Y^c$, then if $ ( \; , \; )$ denotes the usual inner
  product on roots, $(\alpha , \beta ) \geq 0$ because otherwise,
  $\alpha +\beta $ would be a root.  Thus, $(X^c , Y^c) \geq 0$.  But
  $0 \leq (X^c , Y^c ) = (X^c , -Y^c) = -(X^c , Y^c) \leq 0$.  So
  $(X^c , Y^c )= 0$.  \\

\noi If $Y^c = 0 $, then $Y = \Delta$ and we have proven what we
wanted.  If, on the other hand, $Y^c \neq 0 $, then $Z:= (Y^c)^\perp
\cap \Delta $ are the roots corresponding to  some Levi subalgebra $\mathfrak l $
because if $V$ is the vector space spanned by $\Delta$, $W\cap \Delta$
defines a Levi subalgebra for any subspace $ W \subset V $ ~\cite{bou}.  But since
$X^c \perp Y^c$, $X^c \subset Z$, whence $Z^c \subset X $.  \\

\noi There is
some parabolic subalgebra $\mathfrak p = \mathfrak l \oplus \mathfrak
u$, where $\mathfrak u $ is the nilpotent radical of $\mathfrak p$.
Then $\mathfrak u$ is the direct some of root spaces for some subset
of roots $U \subset \Delta$, $\mathfrak p$ corresponds to the roots $Z
\sqcup U$, and $\Delta = Z \sqcup U \sqcup -U$.  Then since $Z^c
\subset X $, $\pm U \subset X$. \\

\noi In order to show that $X = \Delta$, we must show that $X $
contains all simple roots (with simple roots determined by a Borel
subalgebra $\mathfrak b$ containing $\mathfrak u$ and contained in
$\mathfrak p$).  In this situation, the Levi factor $\mathfrak l$
corresponds to some subset of simple roots. Let $\alpha _ 0$ be a
simple root.   Since $\fg$ is simple,
the Dynkin diagram for $\Delta$ is a connected graph, which is a
tree.  Either $\alpha _0 \in U \subset X$, or it is possible to choose a string $\alpha _0
, \alpha _1 , \ldots , \alpha _k$ of simple roots such that $\alpha _k
\in U$ and $\alpha _0
, \alpha _1 , \ldots , \alpha _{k-1} \in Z$ and such that $(\alpha _i
, \alpha _j ) < 0$ if $j = i+1$ and  $(\alpha _i
, \alpha _j ) = 0$ if  $j > i+1$.  This implies that $\alpha _0 +
\alpha_1 + ...+ \alpha _k \in \Delta$ and also that both $\alpha _1 + ...+ \alpha _k
\in U$ and $\alpha_0 +\alpha _1 + ...+ \alpha _k
\in U$.  Hence $\pm (\alpha _1 + ...+ \alpha _k ) \in Z^c \subset X$,
from which it follows that $\alpha _0 = (\alpha_0 +\alpha _1 + ...+
\alpha _k )- (\alpha _1 + ...+ \alpha _k ) \in X + X \subset X$.
Therefore every simple root lies in $X$.  Because $X$ contains all
simple roots and $X=-X$, $X = \Delta$ as desired. $\square$

\end{pf}
 
\noi Lemma~\ref{lk} ensures that for each $i$, either $\Gamma _i
\subset (A \cap - A )$ or $\Gamma _i \subset \theta (A \cap - A ) $.
 If, in the above notation, $\mathcal S $ denotes the set of summands $\{S_1 ,
...S_n \}$ or $\{\Gamma _1 , ..., \Gamma _n \}$, $\theta$ is a
permutation of $\mathcal S$, and we will set $\mathcal S ^\theta = \{
S \in \mathcal S \; | \; \theta S = S \} $.  More generally, for a parabolic $\fp$, $\mathcal S _\fp$ will denote the set of semisimple summands in $\fp$.

\begin{prop} With the above notation, a generalized complex subset $A$ of $\Delta$ is
  equivalent to the following data:
\begin{enumerate}
\item A $\theta$-stable parabolic subalgebra $\mathfrak p \supset
\mathfrak h $,
\item Subsets $T,R$ of $\mathcal S _\fp $ such that $T \sqcup \theta T \sqcup R = \mathcal S _\fp$,
\item Reductive subalgebras $\fq _i $ of $S_i$ with $\mathfrak h \cap S_i
\subset
\fq _i \subset S_i$, $S_i \in \theta T$.  (The corresponding
set of roots in $\Delta (S_i , \mathfrak h \cap S_i ) $ will be
denoted by $A_i$).  
\end{enumerate}
\end{prop}

\begin{pf}  Given $A$, let $\fp = A \cup \theta A$, let $T = \{ \Gamma _i \; | \; \Gamma _i \subset
A \cap -A \; \; but \; \; \Gamma _i \nsubseteq \theta (A \cap - A) \}
$, and let $R= \{ \Gamma _i \; | \; \Gamma _i \subset (A \cap -
A) \cap  \theta (A \cap - A) \} $. Finally, let $\fq _i $ be the
subalgebra containing $S_i \cap \mathfrak h $ and the root spaces
$\Gamma _i \cap A \cap -A $ for $S_i \in T$.\\

\noi Conversely, given the data $\Phi $ , $T$, $R$, $\fq _i $,
$\Phi = \Gamma \sqcup \Psi$, and $\Gamma = \sqcup _i \Gamma _i $.  Let
\[ A = \Psi \cup \{ \Gamma _i \; | \; \Gamma _i \in T \cup R \} \cup \{
A _i \; | \; S_i \in \theta T \}. \]
Then $\sigma A = -\Psi \cup -\{ \Gamma _i \; | \; \Gamma _i \in \theta
T \cup R \} \cup \{
 \theta A _i \; | \; S_i \in \theta T \} $ and $\theta A = \Psi \cup \{ \Gamma _i \; | \; \Gamma _i \in \theta
T \cup R\} \cup \{
A _i \; | \; S_i \in \theta T \} $ .  Clearly $A \cup \theta A = \Phi$,
$ A \cup \sigma A = \Delta$, and $A \cap \sigma A = (\bigcup A _i ) \cup (\bigcup \theta A _i ) \cup R \subset -A$.  This
means that A is a generalized complex subset.  $\square$

\end{pf}

\section{Nilpotent Orbits}\label{sectionnilorbits} 
\noi We classify generalized complex structures on real nilpotent orbits in $\fs \fl _n (\rr)$ by showing that all are $B$-transforms of symplectic structures.  We conjecture that the same is true in any split semisimple Lie algebra over $\rr$.  We reduce this conjecture to the case of distinguished nilpotent elements in simple Lie algebras.  

\subsection{Nilpotent orbits in $\mathfrak{sl}_n (\rr) $}

\noi In this subsection $G = SL_n (\mathbb R)$ and $\fg = \fs \fl _n (\mathbb R)$.

\begin{prop}\label{ij} Let $\mathfrak g = \sln (\rr)$, and let $e\in \mathfrak g$ be any nilpotent element.  If $E \subset \cts$
  is any Lie subalgebra such that $Z_{\cts} (e) \subset E$ and $E +
  \overline E = \cts$, then $E= \cts$.  
\end{prop} 

\begin{pf}First we prove this for the case when $e $ is
  regular.  Since $e$ is regular, there is a standard triple $\{e,h,g\} \subset
\mathfrak g$ with the following properties.  There is a split Cartan
subalgebra $\mathfrak h \subset \mathfrak g $ so that
\[
\mathfrak g = \mathfrak h \oplus \bigoplus _{\alpha \in \Delta }
\mathfrak g _\alpha
\]
There is a Borel $\mathfrak b \subset \mathfrak h $ such that $h\in
\mathfrak h $ and $\alpha (h) = 2$ for all simple roots $\alpha$,
which are determined by $\mathfrak b$.  Also, for every simple root $\alpha$,
$pr _{\mathfrak g  _\alpha} e \neq 0 $.  \\

\noi Now let E be a Lie subalgebra such that  $Z_{\cts} (e) \subset E
$ and $E +\overline E = \cts $.  If dim$E$ = d, we can view $E
\in Gr_d (\cts) $.  Define $L = \lim _{ s \to - \infty} Ad(exp(sh))E
\; \in Gr_d (\cts) $.  Since $Ad(exp(sh))Z_{\cts} (e) = Z_{\cts} (e) $
for all s, $Z_{\cts} (e) \subset L $.  Let $\gamma$ be the unique
root of maximum height.  Since $E + \overline E = \cts$, there is a
vector $v \in E $ whose projection, $pr_{(\cts)_{-\gamma}}
v$, onto the root space $(\cts)_{-\gamma}$ is nonzero. Since
$(\cts)_{-\gamma} $ is the eigenspace on which h has the least eigenvalue, $\lim _{ s \to -
  \infty} Ad(exp(sh))v = x  \in (\cts)_{-\gamma}$ in
$\mathbb P (\cts) $.  Hence, L
contains $Z_{\cts} (e) $ and $(\cts)_{-\gamma}$ and is stable under $ad_h$. \\

\noi To prove the proposition, it therefore suffices to show that for
any Lie subalgebra $L$ of $\cts$ satisfying $ad_h (L) \subset L$,
$Z_{\cts} (e) \subset L$, and $\cts _{-\gamma} \subset L$ must equal
all of $\cts$.  In fact, we show that we can reduce to the case when
$h \in L$.  If we prove this for all such $L$ which also contain h,
then let $L$ be such a Lie algebra which is $h-stable $ but does not necessarily
contain h.  Then $L^\prime := <h>\oplus L $. satisfies all of the
necessary conditions and contains h.  Therefore, $L^\prime = \cts$.
Since $L$ is $h-stable$, $L$ is a non-zero ideal in $\cts$ which was
assumed to be simple. In fact, only $\fg$ was assumed to be simple,
but $\fg$ simple and split implies that $\fg$ is not complex and
therefore $\cts$ is simple.   Therefore $L = \cts$.  Hence, it suffices to
show that if $L \subset \cts$ is a Lie subalgebra such that $h \in L
$, $Z_{\cts} (e) \subset L$, and $\cts _{-\gamma} \subset L$, then $L=
\cts$.  To this end, suppose that $L$ is any such
subalgebra. \\

\noi There is a decomposition $L = \mathfrak m \oplus \mathfrak n $,
where $\mathfrak n $ is the nilpotent radical of $L$ and $\mathfrak m
$ is a reductive subalgebra.  Since $h \in L$ is semisimple and $\mathfrak n $ is an
ideal, $ad_h$ acts diagonally on $L$ and $\mathfrak n$.  Hence $L$
decomposes as a direct sum of its $ad_h$-eigenspaces, $L_i$ ($i \in
2\mathbb Z $).  Also $\mathfrak n $ decomposes into a direct sum of
its $ad_h$-eigenspaces $\mathfrak n _i $.  \\

\noi First we show that $\mathfrak n _0 = 0 $. By definition,
$\mathfrak n _0  \subset Z_{\cts} (h) =  \mathfrak h _\cc $.
Therefore any all elements in $\mathfrak n _0 $ act diagonally on $L$
and $\mathfrak n $.  Since $\mathfrak n $ is nilpotent, this implies
that $\mathfrak n_0 \subset Z(L) \subset Z_L(e, {\cts}_{-\gamma} )
\subset Z_{\cts} (e ,{\cts}_{-\gamma} ) $.  Thus, to show that
$\mathfrak n _0 = 0$, it suffices to show that $ Z_{\cts} (e
,{\cts}_{-\gamma} ) = 0 $.  \\

\noi Since $e$ is regular, we may assume that $\mathfrak h$ is the
diagonal matrices in $\sln$ and 

\[ e = \sum _{i = 1} ^{n-1} X_{i, i+1}
\] 
where $X_{i,j}$ is the matrix with 1 in the $ij$-th entry and 0
elsewhere.  Then $X_{n, 1} $ is a basis for ${\cts}_{-\gamma} $ if we
take the standard Borel  (in fact the only Borel containing $e$).  It
is easily checked now that $ Z_{\cts} (e
,{\cts}_{-\gamma} ) = 0 $.

\noi This shows that $  Z_{\cts} (e
,{\cts}_{-\gamma} ) = 0 $ and therefore $\mathfrak n_0 = 0$.  
If $0 \neq x \in \mathfrak n _i $ for $i < 0$, then since $Ker(ad_e)
\cap \oplus_{i\leq 0 } {\cts}_i =0$, $0\neq ad_e ^k x \in \mathfrak n
_0 = 0 $, which gives a contradiction.  Therefore $\mathfrak n_i =0$
for $i \leq 0 $ and so 
\[ 
\mathfrak n = \bigoplus _{i > 0 } \mathfrak n _i \;.
\]

\noi However, since $\gamma$ is the highest root, 
\[ 
[X_{-\gamma},
  \cts] \subset \bigoplus _{i \leq 0 } (\cts)_i \;.
\]

\noi Since $X_{-\gamma} \in L$ and $\mathfrak n $ is an ideal in $L$, this
means that $[X_{-\ga} , \mathfrak n ] = 0$.  \\

\noi We have seen that $\mathfrak n  \subset Z_{\cts} (X_{-\gamma} )$.
If $\mathfrak n \neq 0$, then there is some nonzero $x \in \mathfrak
n$.  $\mathfrak n $ is $ad_e$-stable since it is an ideal of L. But
since $e$ is nilpotent, $ad_e ^k x = 0 $ for some minimal $k$. Then
$0 \neq ad_e ^{k-1} x \in \mathfrak n \cap Z_L (e) \subset Z_L
(e,X_{-\gamma} ) \subset  Z_{\cts} (e, {\cts} _{-\gamma} ) = 0 $.
  This is a contradiction.  Therefore $\mathfrak n = 0$.  \\

\noi Because $\mathfrak n = 0$, $L = \mathfrak m$, which is
reductive.  Hence, $L = S \oplus Z$, where $Z$ is the center of $L$
and $S$ is semisimple.  Since $[h,e] = 2e$, $e \in S$.  Because $S$ is
semisimple, there is a standard triple $\{ e, h^\prime , f^\prime \}
\subset S \subset \cts$.  However, when we view $\cts$ as an
$\mathfrak{sl} _2 $-module for the triple $\{ e, h^\prime , f^\prime
\} $, $\cts$ is generated by $Z_\cts (e) $.  Since   $Z_\cts (e) $ and
$\{ e, h^\prime , f^\prime \} $ are contained in L, $L= \cts$.

\noi We have proven the theorem for the case when $ e $ is regular.  If
$e$ is not regular, Subsection~\ref{27nov1} shows that the theorem still
holds.  We need only observe that if $\mathfrak l $ is a Levi
subalgebra of $\sln$, then it is a direct sum of its center and semisimple Lie
algebras, each isomorphic to some $\mathfrak {sl} _m $.  $\square$ 
\end{pf}

\begin{cor} Every equivariant generalized complex structure on $\fs \fl _n (\mathbb R)$ is a B-transform of a symplectic structure. 
\end{cor}
\begin{pf} This now follows from Proposition~\ref{ij} because generalized complex structures of the form $L(TM_\cc , \vep)$ are precisely the B-transforms of symplectic structures.  $\square$
\end{pf} 

\noi With this in mind, the following proposition describes \gcstrs on $\OO _e$ for any nilpotent $e \in \fs \fl _n (\rr)$. 

\begin{prop} The equivariant generalized complex structures on a
  nilpotent orbit $\mathcal O _e $ in $\sln (\rr) $ are given by an 
  open set in $Z(Z_{\fg _ \cc} (e) )^{Z_G (e)}$.  Specifically, the equivariant
  generalized complex structures are parameterized by all $t = t_r +
  it_i \in Z(Z_{\fg _\cc (e) })$ such that  $Z_{\fg _\cc} (t _i) = Z_{\fs \fl _n \cc}(e) $ and such that $t$ is $Z_G(e)$-invariant.
\end{prop}  

\begin{pf} Let $\fg = \fs \fl _n (\rr)$.  Proposition~\ref{ij} shows that any GC-pair $(E, \vep ) $
  must in fact satisfy $E = \cts$.  Since $\cts$ is semisimple, $H^2
  (\cts , \cc ) = 0 $, whence $\vep = \phi \circ [ \; , \; ] $ for
  some $\phi \in \cts ^* $.  Any such $\vep $ satisfies $d_E \vep = 0
  $.  But $\vep $ must also satisfy $\vep ( Z_{\cts} (e) , \cts ) =
  0$.  Using the Killing form, $\cts \simeq \cts ^*$, so we may
  identify $\phi $ with some $t\in \cts$
\begin{eqnarray*} & \; &
\phi ([Z_{\cts} (e) , \cts
  ] ) = 0 \\
 &\iff& \kappa (t, [Z_{\cts} (e) , \cts
  ] ) = 0 \\
&\iff& \kappa ([t , Z_{\cts} (e) ] , \cts) = 0 \\
&\iff & [t ,
  Z_{\cts} (e) ] = 0 \\
&\iff & t \in Z( Z_{\cts} (e) ) 
\end{eqnarray*}  
\noi Therefore, $\vep
  $ satisfies condition 4 of Proposition~\ref{b} if and only if $t \in Z(
  Z_{\cts} (e) ) $.\\

  \noi For $(E , \phi \circ [ \; , \; ] )$ to be a  GC-pair, condition
  5 of Proposition~\ref{b} must also be satisfied.  Breaking up $\phi =
  \phi _r + i \phi _i $ into real and imaginary parts, condition 5 is equivalent to requiring that for all $x \notin Z_{\cts} (e)
  $, there exists $y \in \cts$ such that $\phi _i [x,y] \neq 0 $.  In
  other words, we require that $\phi _i \circ [ \; , \; ]$ be
  nondegenerate on $\wedge ^2 \cts \, /Z_{\cts}(e)$.  Again, when $\phi _i $ is identified with $t_i$ via
  the Killing form, this is equivalent to asking that for each $x \notin Z_{\cts} (e)
  $, there exists $y \in \cts$ such that $\kappa ([t, x] ,y ) \neq
  0$.  But since $\kappa$ is nondegenerate, this happens exactly when
  $[t_i , x ] \neq 0 $ for all $x \notin Z_{\cts} (e) $.  This happens
  if and only if $Z_{\cts} (t _i) \subset Z_{\cts} (e) $.  However,
  $t_i \in    Z( Z_{\cts} (e) )$; hence  $Z_{\cts} (t _i) = Z_{\cts} (e) $. $\square$

\end{pf}

\subsection{Reduction to Distinguished Orbits in Simple Lie Algebras}\label{27nov1}

\noi We wish to extend the results in the previous section to nilpotent
orbits in arbitrary split semisimple Lie algebras.  \\

\noi For brevity, if $\fg$ is a split semisimple real Lie algebra and $e$ is a
nilpotent element, let $P(\fg , e)$ denote the following statement:\\

If $E$ is a subalgebra of $\cts$ such that $Z_\cts (e) \subset E$ and
$E + \overline E = \cts$, then $E= \cts$.  \\

\begin{Conj} $P(\fg, e)$ is true for any split semisimple Lie algebra $\fg$ and any nilpotent $e \in \fg$.  
\end{Conj}

\noi The following results show that it suffices to prove the conjecture for distinguished nilpotent orbits in simple, split Lie algebras.  

\begin{lem} In order to prove $P(\fg , e) $ for any
  split semisimple Lie algebra $\fg$ and any nilpotent $e \in \fg$, it
  is enough to prove the result when $\fg $ is
  simple.  

\end{lem}

\begin{pf} First assume that $P(\fg , e ) $ is true whenever $\fg $ is
  simple and $e\in \fg$ is nilpotent.  Now let $\mathfrak g$ be any split semisimple Lie algebra so that $\fg = \oplus \fg _i$
 is a direct sum of split simple Lie algebras.  One can complete $e$ to a standard $\fs \fl _2$ triple $\{ e,
 h,f\}$ such that $e= \sum e_i$, $f = \sum f_i$, $h = \sum h_i$, and each $\{e_i, h_i, f_i \}$ is a
 standard triple in $\fg _i$.\\

\noindent Let $\pi _i : \cts \longrightarrow (\fg _i )_\cc $ denote the
projection map.  Clearly $Z_{(\fg _i )_\cc } (e_i) \subset Z_\fg (e)
\subset E$,
whence $Z_{(\fg _i )_\cc } (e_i) \subset \pi _i (E) $.  We know that $(\fg _i )_\cc = \overline{(\fg _i
  )_\cc } $ and $E+\overline{E} = \cts$, from which it follows that
$\pi _i E + \overline{ \pi _i E} = {(\fg _i )_\cc } $.  But it was
assumed that $P(\fg _i , e_i )$.  This implies that $\pi _i E = {(\fg _i )_\cc }
$ because $\fg _i $ is split simple.\\

\noindent Consequently, if $x_i \in {(\fg _i )_\cc }$, there exists $x
\in E $ such that $\pi _i x = x_i $.  In particular, there exists
$f^\prime \in E$ such that $pr_{(\cts)_i} f^\prime = f_i$.  Since $ e_i \in Z_\fg (e) $,
$h_i = [e_i , f_i ] = [e_i , f^\prime] \in E $.  Then it is also the
case that $-2f_i = [h_i , f_i ] = [h_i , f^\prime ] \in E$.  We may
conclude that $f_i \in E$, for each $f_i \in E$.  The centralizer
$Z_\cts (e) $ of e and $\{e,h,f\} $ generate $\cts$ as a subalgebra.
Therefore $E = \cts$.   $\square$

\end{pf}

\noi We now aim to show that in order to prove $P(\fg , e)$ for arbitrary
nilpotent elements $e$, it is enough to show this in the case when $e$
is a distinguished nilpotent in $\cts$. \\

\noi Define a \emph{split Levi subalgebra} of $\cts$ to be any Levi
subalgebra of $\cts$ that contains a split Cartan subalgebra of
$\fg$.  Note that any split Levi subalgebra $\mathfrak l $ satisfies
$\mathfrak l = \overline {\mathfrak l}$.

\begin{lem} Let $\fg$ be a split semisimple Lie algebra.  If $e \in
  \fg$ is a nilpotent element that is not distinguished, then  there
  is a proper split Levi subalgebra containing $e$.  
\end{lem}

\begin{pf} Let $\{e,h,f\} \subset \fg$ be a standard triple.
 There is a Cartan involution $\theta$ of $\fg$ such that $\theta e =
 -f $, $\theta f= -e $, and $\theta h = -h $.  This is possible when $\{e,h,f\}$ span 
$\fg$
 and also in general because any Cartan involution of a semisimple subalgebra may be
 extended to a Cartan involution of the entire Lie algebra (see 9.4.1 ~\cite{com}).  Let $\fg =
  \mathfrak k \oplus \mathfrak p $ be the corresponding Cartan
  decomposition.  By the theory of $\mathfrak {sl} _2 $ representations,
 $e$ is distinguished in $\cts$ if and only if $Z_{\cts} (e,h) =
 Z_{\cts} (e,h,f)  $ is $0$.  Since $e$ is not distinguished, there exists
 a nonzero $x \in Z_{\cts} (e,h) $.  Since $\{ e,h,f\} \subset \fg$,
 this means that there exists a non-zero $y \in Z_\fg (e,h) $.  If we
 let $u = y - \theta y $, then $u \in \mathfrak p$ and is therefore
 diagonalizable.  Note that $0 = \theta ([y,h] ) =- [\theta y , h]$ so
 that $\theta y \in Z_\fg (h)$.  Also, $0 = \theta ([y, f] ) =
 -[\theta y , e ]$ so that $\theta y \in Z_\fg (e,h) $ and therefore
 $u \in Z_\fg (e,h) $.  Notice that $u$ lies in some  maximal abelian
 subspace $\mathfrak v  $ of $\mathfrak p $. All maximal abelian
 subspaces of $\mathfrak p $ are conjugate by the group $Int(\fg)$ of
 inner automorphisms.  Therefore $\mathfrak v $ is conjugate to a
 split Cartan subalgebra and must itself be a split Cartan subalgebra.
 Since $u$ is semisimple, $Z_\cts (u) $ is a Levi subalgebra.  It is
 split since $\mathfrak v \subset Z_\cts (u) $, and finally $e \in
 Z_\cts (u) $.  So $e$ lies in a a split Levi subalgebra.  If $u \neq
 0$, then $Z_\cts (u) \neq \cts$, and $e $ lies in a proper split Levi
 subalgebra. \\

\noi It remains to show that there is a $y \in Z_\fg (e,h) $ such that
$u = y - \theta y \neq 0 $.  In other words, we need to show that
$Z_\fg (e,h )$ is not contained in $ \mathfrak k$.  To see this, we
can embed $\mathfrak g $ into some $\sln$ by the adjoint map.
Henceforth in this proof, we will view $\fg$ as a subalgebra of
$\sln$.  If $x \in \sln$ is a nilpotent element which is in Jordan
form, it is possible to choose a standard triple $\{x,y,z \} $ whose
semisimple element $z$ is a diagonal matrix.  It can be shown that any
semisimple element $s \in Z_{\sln} (x,z) $ has real eigenvalues.  Now,
if $Z_\fg (e,h ) \subset \mathfrak k$, let $0 \neq y \in Z_\fg (e,h) $.
Since $y \in \mathfrak k $, $y$ has purely imaginary eigenvalues.
However, there is an automorphism $g \in Aut(\slnr ) $ such that $g.e
= x$ and $g.h = z$, where $x$ is in Jordan form and $z$ is a diagonal
matrix.  But then $g.y \in Z_{\sln} (x,z)$ is a semisimple element with
purely imaginary eigenvalues.  This is a contradiction. \\

\noi Therefore $Z_\fg (e,h) $ is not contained in $\mathfrak k$, and this completes
the proof.  $ \square $

\end{pf}

\begin{lem} Let $\fg$ is a split semisimple Lie algebra.  Any
  nilpotent $ e \in \fg$ is contained in a split Levi subalgebra of
  $\cts$ that is minimal among all Levi subalgebras containging $e$.
\end{lem}

\begin{pf} We proceed by induction on dim$(\fg)$.  The case when dim$\fg
  = 3 $ is trivial.  Let $e \in \fg$ as above.  If $e$ is
  distinguished, then we are done.  If not, then there is a split Levi
  subalgebra $\mathfrak l _\cc$ containing $e$ and a split Cartan
  $\mathfrak h$.  Here $\mathfrak l
  _\cc $ is the complexification of some $\mathfrak l \subset \fg$.
  Since $\mathfrak l _\cc $ is reductive, $\mathfrak l = [ \mathfrak l
  , \mathfrak l ] \oplus Z(\mathfrak l ) $, and $e \in \mathfrak l
  ^\prime := [\mathfrak l ,
  \mathfrak l] $.  If $e$ is distinguished in $\mathfrak l ^\prime $
  (which is a split real form of $\mathfrak l ^\prime _\cc$), then
  $\mathfrak l ^\prime _ \cc$ is minimal.  Otherwise, there exists a
  split Levi subalgebra $\mathfrak m $ of $\mathfrak l ^\prime $ with
  $e \in \mathfrak m \varsubsetneq \mathfrak l ^\prime $ in which $e$
  is distinguished.  This follows by induction hypothesis.  The Levi subalgebra $\mathfrak
  m $ contains a split Cartan $\mathfrak s$, and we wish to show that
  $\mathfrak s \oplus Z(\mathfrak l ) $ is a split Cartan for $\fg$.
  Note that $\mathfrak h \cap \mathfrak l ^\prime $ is a split Cartan
  subalgebra of $\mathfrak l ^\prime $.  But $\mathfrak s$ is also a
  split Cartan subalgebra of $\mathfrak l ^\prime $.  Therefore there
  exists $g \in Int (\mathfrak l ^\prime )$ such that $g.\mathfrak s =
  \mathfrak h \cap \mathfrak l ^\prime $.  We may view $g$ as an
  element of $Int (\mathfrak g )$. so that $g.(\mathfrak s \oplus
  Z(\mathfrak l ) ) = g. \mathfrak s + g.Z(\mathfrak l ) = \mathfrak h
  \cap \mathfrak l ^\prime + Z(\mathfrak l ) = \mathfrak h $. Now
  Since $\mathfrak s \oplus Z(\mathfrak l ) $ and $\mathfrak h $ are
  conjugate by a memeber of $Int(\mathfrak g) $,  $\mathfrak s \oplus
  Z(\mathfrak l ) $ is a split Cartan subalgebra of $\mathfrak g$.
  Then $\mathfrak m _\cc \oplus Z(\mathfrak l )$ is a split Levi
  subalgebra in which $e$ is distinguished.  That is to say,
  $\mathfrak m _\cc \oplus Z(\mathfrak l )$ is a minimal Levi
  containing $e$, and it is split. $\square$
\end{pf} 

\begin{prop} Let $\fg$ be a split semisimple Lie algebra, and let $e
  \in \fg$ be any nilpotent element.  If $\mathfrak l$ is a minimal
  split Levi subalgebra containing $e$, then $P([\mathfrak l ,
  \mathfrak l ] , e)$ implies $ P(\fg , e ) $.  Therefore to prove $P(\fg , e) $
  for any split semisimple $\fg$ and nilpotent $e \in \fg$, it suffices to prove that for any split semisimple Lie algebra $\fs$ and distinguished nilpotent $x \in \fs$, $P(\fs , x)$.
\end{prop} 

\begin{pf} Assume that $P([\mathfrak l ,
  \mathfrak l ] , e)$, and let $E$ be a subalgebra of $\cts$ such that
  $Z_\cts (e) \subset E$ and $E + \overline E  = \cts$.  Let $\mathfrak h \subset \mathfrak l $ be a split Cartan
  subalgebra.  There exits $h \in \mathfrak h $ such that $\mathfrak l
  $ is the 0-eigenspace for $ad_h$.  But since $e \in \mathfrak l $,
  $h \in Z_\cts (e) \subset E$.  Hence E decomoposes as a direct sum
  of eigenspaces for $ad_h$, and $pr_{\mathfrak l _\cc} E = E \cap
  \mathfrak l _\cc$.  Here $pr_{\fl _\cc}$ is projection onto the 0-eigenspace, $\fl _\cc$, for $ad_h$.  This implies that $pr_{\mathfrak l _\cc } E $ is
  a subalgebra of $\mathfrak l _\cc $ containing $Z_{\mathfrak l _\cc
  } (e) $ and such that  $pr_{\mathfrak l _\cc } E + \overline {
  pr_{\mathfrak l _\cc } E } = \mathfrak l _\cc$.  \\
  
  \noi Note that $\fl _\cc = [\fl , \fl]_\cc \oplus Z(\fl)_\cc$, and $Z(\fl) \subset Z_{\fg} (e) \subset E$.  Hence, $( E \cap [\fl , \fl ]_\cc ) + \overline {( E \cap [\fl , \fl ]_\cc ) } = [\fl , \fl ] _\cc$.  In conjunction with the fact that $Z_{[\fl , \fl ]_\cc } \subset E \cap [\fl , \fl ] _\cc$ and $P([\fl , \fl] , e)$ implies that $E \cap [\fl , \fl ]_\cc = [\fl , \fl ] _\cc$ and therefore $\fl _\cc \subset E$. This means that $e$ lies in some semisimple Lie subalgebra of $E$, namely $[\fl , \fl]$.  It follows that there is a standard triple $\{e,h,f \} \subset E$.  But $f$ and $Z_\cts (e)$ together generate all of $\cts$.  $\square$

\end{pf}

\section{Riemannian Symmetric Spaces}\label{sectionriemann}

\noi Let $(M, Q)$ be a Riemannian symmetric space $M$ with metric $Q$, and let ${\bf G } = I(M)_0$ denote the identity component of the isometry group.  We fix a point $p_0 \in M$, and let $K \subset {\bf G }$ be the subgroup fixing $p_0$.  We know that $\fk = Lie \, K$ is the set of fixed points of an involutive isometry $\theta$ of $\fg$.  It is known that $\fg$ decomposes as a direct sum of Lie algebras $\fl _c $, $\fl _n$, $\fl _e$, each of which is fixed by $\theta$, where $(\fl _c , \theta )$, $(\fl _n , \theta)$, $(\fl _e, \theta )$ are orthogonal symmetric pairs of the compact, non-compact, and Euclidean type respectively.  Each of $\fl _c$ and $\fl _n $ decompose into a direct sum of irreducible orthogonal symmetric pairs (in the sense of Definition~\ref{26nov5}).  Thus, $\fg$ is a direct sum of ideals $\fg = \fl _e \oplus \bigoplus _{i = 1} ^n \fg _i$, where each $(\fg _i , \theta)$ is an orthogonal symmetric pair of the compact or non-compact type. Thus, the Lie algebra $\bla$ is a direct sum, which we denote by $\bla = \fu _e \oplus \bigoplus _{i=1} ^n \fu _i$.  The main result of this section is the following theorem.

\begin{thm}\label{26nov8} With the above notation, any $G$-invariant generalized complex structure $L$ , when viewed, as a subalgebra $L \subset \bla$ is a direct sum $L = L_e \oplus \bigoplus _{i=1} ^n L _i \subset (\fu _e \oplus \bigoplus _{i=1} ^n \fu _i )_\cc$.  Each $L_i$ represents either a complex structure or a B-transform of a symplectic structure.  When $M$ is simply connected, the generalized complex structure is a product of generalized complex structures on $M  = M _e \times \Pi _{i=1} ^n M _i$.
\end{thm}

\noi In Subsections~\ref{26nov6} and \ref{26nov7} we consider an arbitrary Riemannian symmetric pair $(G,K)$ with $G/K \simeq M$.  We assume that $G$ is connected and acts on $M$ by isometries, so we have a map $\tau : G \lra I_0 (M)$.  Just in case $G$ doesn't act effectively, let $N = Ker (\tau ) $, and $G/K \simeq (G/N)/(K/N) \simeq M$.  We can show that $K/N$ is compact \cite{hel} and replace $(G,K)$ by $(G/N, K/N)$.  We therefore assume from the beginning that $G$ acts effectively.  Hence, $K$ is compact.  Also $\fk$ contains no non-zero ideal of $\fg$ because if it did contain such an ideal $\mathfrak i$ the connected subgroup $I$ corresponding to $\mathfrak i$ would be normal and therefore act trivially on $G/K$, which would contradict that $G$ acts effectively.  We now have a closed embedding $G \hookrightarrow I_0(M)$, where $ U_0 \subset K \subset U $ (where $U$ is the subgroup fixing $eK \in G/K =M$). 


\subsection{Irreducible Semisimple Symmetric Spaces}\label{26nov6}

\noi If $(G,K) $ is a Riemannian
Symmetric pair with G semisimple, then the Lie algebra $\fg$ has a
decomposition $\fg =
\mathfrak k \oplus \mathfrak p $ into
$\pm1$-eigenspaces for an involution $\theta$ corresponding to the
pair $(G , K) $.  In other words, $(\fg
, \theta ) $ is an orthogonal symmetric pair with decomposition $\fg =
\mathfrak k \oplus \mathfrak p $.

\begin{Def}\label{26nov5} An orthogonal symmetric pair  $(\fg , \theta) $ is
called \emph{irreducible} if $\fg$ is semisimple, $\mathfrak k$
contains no nonzero ideal of $\fg$, and $\mathfrak p $ is an
irreducible $\mathfrak k$-module ($\mathfrak k$, $\mathfrak p$ as
above). A Riemannian symmetric space $(G,K)$ is called \emph{irreducible}, if
G is semisimple and the corresponding orthogonal symmetric pair is
irreducible.  
\end{Def}

\noi Generalized complex structures on $G/K$ are given by GC-pairs $(E,\vep) $.  The conditions on E are that $E$ is
a subalgebra containing $\csts$ and $E + \overline E = \cts$.  
  It is known that if $\fg$ is noncompact, then $\theta$ is a
Cartan involution of $\fg$.  If $\fg$ is compact, then the dual
orthogonal symmetric pair $(\fg ^\bullet = \mathfrak k \oplus i \mathfrak p
, \theta ^\bullet ) $ is non-compact, and $\theta ^\bullet $ is a Cartan
involution.  GC-pairs $(E,\vep)$ for $(\fg , \theta )$ are the same as
those for  $(\fg ^\bullet = \mathfrak k \oplus i \mathfrak p
, \theta ^\bullet ) $ except that condition 5 of Proposition~\ref{b} is different
in the two cases.  Moreover, when it is convenient for finding
candidates for GC pairs $(E,\vep)$, one may assume
that $(\fg , \theta )$ is of noncompact type or of compact type.

\begin{lem}\label{feb19a} Let $\fg = \mathfrak k \oplus \mathfrak p $ be a Cartan
  decomposition of a real semisimple Lie algebra $\fg$.  Suppose that
  $\mathfrak p$ is an irreducible $\mathfrak k$-module.  Then if $E$
  is any subalgebra of $\cts$ containing $\csts$, then either $E = \cts$ or $E = \csts \oplus \mathfrak a $, where $\mathfrak a$ is an
  irreducible $\mathfrak k$-submodule of $\mathfrak p _\cc$ and
  $\mathfrak a \oplus \bar {\mathfrak a } = \mathfrak p _\cc $.

\end{lem}

\begin{pf} Since $\csts \subset E$, $E = \csts \oplus \mathfrak a $
  for some subspace $\mathfrak a $ of $\mathfrak p _\cc$.  The
  subspace $\mathfrak a $ must therefore be a $\mathfrak k$-module.  But because $\overline {\mathfrak
  a }$ is
  also a $\mathfrak k$-module, so is $\mathfrak a \cap \bar {\mathfrak
  a}$.  However, since $\mathfrak a \cap \bar {\mathfrak
  a}$ is stable under complex conjugation,  $\mathfrak a \cap \bar {\mathfrak
  a} = V_\cc$ for some subspace $V \subset  \mathfrak a \cap \bar {\mathfrak
  a }\cap \mathfrak p$.  Note that $V$ is also a $\mathfrak k$-module,
  being the intersection of three $\mathfrak k$-modules.  We assumed
  that $\mathfrak p$ was an irreducible $\mathfrak k$-module.  Hence,
  either $V= \mathfrak p _\cc$, in which case $E = \cts$, or $V = 0$,
  in which case $\mathfrak p _\cc = \mathfrak a \oplus \bar {\mathfrak
  a}$.  \\

\noi In the second case, we must show that $\mathfrak a $ is
irreducible.  Suppose that there is a proper submodule $W \subset \mathfrak
a$.  This would mean that $W\oplus \bar W \varsubsetneq \mathfrak p
_\cc$ is a submodule.  Again, though, $W\oplus \bar W $ is stable under conjugation,
so it is the complexification of some $U \subset \mathfrak p$, which
must also be a $\mathfrak k$-module.  Necessarily $U \varsubsetneq
\mathfrak p$, which contradicts the fact that $\mathfrak p$ is
irreducible.  $\square$ 

\end{pf} 

\begin{lem}\label{feb19b} Let $(\fg = \mathfrak k \oplus \mathfrak p , \theta ) $ be
  an irreducible semisimple orthogonal symmetric pair corresponding to the Riemannian symmetric pair $(G, K)$.  Let $K_0$ denote the connected component of $K$. Either $K_0$ is semisimple,
  or $K_0 \simeq S^1 \times K_{ss}$, where $K_{ss}$ is semisimple.  Equivalently, $dim(Z_\mathfrak k ) \in \{ 0 ,1 \} $.  
\end{lem}
\begin{pf} $\mathfrak k$ is compact , hence reductive, so
  $\mathfrak k = Z(\mathfrak k) \oplus [\mathfrak k , \mathfrak k]$, and on the group level, $K_0 = Z_0 K_{ss} $, where $Z_0$
  is the connected component of the center of $K_0$.  Both $
  Z_0$ and $K_{ss}$ are closed.  \\

\noi Let $\mathfrak a \subset \mathfrak p _\cc$ be and irreducible
submodule.  Let $Z = Z(\mathfrak k)$.  Because $\fk$ is compact and contains $Z$, $Z$ acts on $\mathfrak a$ by
simultaneously diagonalizable matrices so that $\mathfrak a$
decomposes as a direct sum of eigenspaces for $Z$.  Then since 
$Z$ commutes with $\mathfrak k _{ss} = [\mathfrak k , \mathfrak k]$, $
\mathfrak k _{ss}$ preserves each of these eigenspaces.  Since
$\mathfrak a$ is irreducible, there can therefore be at most one
eigenspace.  That is, $Z$ acts by a scalar $\lambda \in Z^*$.  on
$\mathfrak a$: $[x,v] = \lambda (x) v $ for all $x \in Z$, $v \in
\mathfrak a$.  If $dim(Z) >1$, there exists a nonzero $x\in
Ker(\lambda) $.  This implies that $x \in Z$ and $[x, \mathfrak a] =
0$.  But also $[x,\overline{\mathfrak a} ] = [\overline x , \overline
  {\mathfrak a} ] = \overline{[x,\mathfrak a]} = 0$ so that $x \in
Z_{\cts}$, which gives a contradiction because $\cts$ is semisimple.
Therefore $dim(Z)\leq 1$.  If $Z=0$, then $\mathfrak k$ is semisimple.
If $dim(Z) = 1$, $Z_0 $(the connected subgroup with Lie algebra $Z$) is a compact connected abelian subgroup of dimension 1
and is therefore isomorphic to $S^1$.  $\square$ 

\end{pf}

\begin{prop}\label{xy} Let $(\fg = \mathfrak k \oplus \mathfrak p , \theta ) $ be
  an irreducible semisimple orthogonal symmetric pair.  Then any GC-pair $(E, \vep) $ for $G/K$ must be of one of
  the following two forms: \begin{enumerate}
\item $E = \cts$ and $\vep = \phi \circ [ \; , \; ] $ for some $\phi \in
  Ann_{\cts} ( [\csts , \csts] \oplus \mathfrak p _\cc ) \simeq
  Z(\csts)^* $.  We decompose $\phi$ as $\phi = \phi _r + i \phi_i$ with $\phi _r , \phi _i
  \in \fg ^*$.  A sufficient and necessary condition for the existence of a GC-pair
  is that $\phi _i \circ [ \; , \; ]$ is symplectic (non-degenerate).
\item $E = \csts \oplus \mathfrak a$ with $ \mathfrak a \oplus
  \overline{\mathfrak a} = \mathfrak p _\cc$ and $\vep = 0$ so that
  $(E,0)$ is a complex structure.  
\end{enumerate}

\noi Therefore on an irreducible semisimple Riemannian symmetric
space $G/K$, the only G-invariant generalized complex structures are complex
structures and complex symplectic structures (i.e. closed, complex
two-forms) for which the imaginary part is symplectic.
\end{prop}

\begin{pf}
\begin{enumerate} 
\item  Recall that $\cts$ semisimple implies $H^2 (\cts
  , \cc) = 0$ ~\cite{wei}.  Thus, $d_{\cts} \vep = 0 $ implies that $\vep = \phi
  \circ [ \; , \; ] $ for some $\phi \in \cts ^*$.  But $\tep (\csts )
  = 0 $ implies $ \phi ([\csts , \cts])= 0$.  Therefore
  $\phi \in Ann_{\cts} ( [\csts , \cts])$.  It is easily seen that
  $[\cts , \csts] = [\csts , \csts] \oplus \mathfrak p _\cc $.  Note however that $\mathfrak k$ is compact, hence reductive, so
  $\mathfrak k = Z(\mathfrak k) \oplus [\mathfrak k , \mathfrak k]$.
  From this it follows that $
Ann_{\cts} ( [\csts , \csts] \oplus \mathfrak p _\cc ) \simeq
  Z(\csts)^* $.  Any such $\vep = \phi \circ [ \; , \; ] $ with $\phi \in
  Ann_{\cts} ( [\csts , \csts] \oplus \mathfrak p _\cc )$ will satisfy
  conditions 1-4 of Proposition~\ref{b}.  In order to satisfy condition 5,
  let $x\in \mathfrak p _\cc$ and suppose that $\phi([x,y]) -
  \overline{\phi[\overline x , \overline y]} = 0  $ for all $y \in
  \cts$. Since, $2i\phi ([x,y])=\phi([x,y]) -
  \overline{\phi[\overline x , \overline y]} $, this means that $\phi
  _i ([x,y])=0 $ for all $y \in \cts$.  Condition 5 is thus satisfied
  if and only if for each $x \in \mathfrak p _\cc$, there is some $y
  \in \cts$ such that $\phi _i ([x,y]) \neq 0$.  In other words, this
  happens when $\phi _i \circ [ \;  , \; ] $ is nondegenerate.

\item We have already seen that If $E \neq \mathfrak p _\cc$, then $E$
  is of this form.  We have $\cts = \csts \oplus \mathfrak a \oplus
  \overline{\mathfrak a} $ and so $\cts ^* \simeq \csts ^* \oplus
  \mathfrak a ^* \oplus
  \overline{\mathfrak a} ^*$, and $ E ^* \simeq  \csts ^* \oplus
  \mathfrak a ^* $.  If $\vep \in \wedge ^2 E^* = \wedge ^2 (\csts ^* \oplus
  \mathfrak a ^* )$ and $d_E \vep = 0$, we may extend $\vep $ by $0$
  on $\cts \times \overline{\mathfrak a} $ to $\vep ^\prime \in \wedge
  ^2 \cts ^* $.  It can easily be checked that since $(\vep ^\prime
  )_\sharp $ vanishes on $\csts \oplus \overline{\mathfrak a} $ that
  $d_{\cts} \vep ^\prime $ vanishes on each of the following sets:
$\csts \times \mathfrak a \times \overline{\mathfrak a}$, $\mathfrak a
  \times \mathfrak a \times \overline{\mathfrak a} $, $\csts \times
  \csts \times \overline{\mathfrak a}$, $\mathfrak a \times
  \overline{\mathfrak a} \times \overline{\mathfrak a} $, $\csts \times
  \overline{\mathfrak a} \times \overline{\mathfrak a} $, and $\overline{\mathfrak a} \times
  \overline{\mathfrak a} \times \overline{\mathfrak a} $.
  Furthermore, $d_{\cts} \vep ^\prime = d_E \vep = 0 $ when restricted
  to $(\csts \oplus \mathfrak a ) ^3 $, and therefore $d_{\cts} \vep
  ^\prime = 0$.  \\

\noi Since $\cts$ is semisimple, $H^2 (\cts , \cc ) = 0$, whence $\vep
^\prime = \phi \circ [ \; , \; ] $ for some $ \phi \in \cts ^* $.
Recall that $\phi \in Ann_{\cts ^*} ([\csts , \csts ] \oplus \mathfrak
p _\cc ) \simeq Z(\csts )^*$.  If $\vep \neq 0 $, then $\phi \neq 0$
and so $\phi (x) \neq 0$ for some $x $ spanning $Z(\csts) $.  Also,
since $\csts = Z(\mathfrak k)_\cc \oplus [\mathfrak k , \mathfrak
  k]_\cc$ and $\vep \neq 0$, the projection $pr_{Z(\csts)} ([\mathfrak
  a , \mathfrak a] ) = Z(\mathfrak k)_\cc$.  But then $Z(\mathfrak
k)_\cc = \overline {  Z(\mathfrak
k)_\cc } = \overline {pr_{Z(\csts)} ([\mathfrak
  a , \mathfrak a] ) } = pr_{Z(\csts)} (\overline {[\mathfrak
  a , \mathfrak a] ) } = pr_{Z(\csts)} ([\overline{\mathfrak
  a } ,\overline{ \mathfrak a }] )$.  In other words, $\vep ^\prime
(\overline{\mathfrak a} , \overline {\mathfrak a} ) = \phi
([\overline{\mathfrak a} , \overline {\mathfrak a} ]) = \phi
(Z(\mathfrak k ) _\cc ) = \phi (Z(\csts) )\neq 0 $.  However, $\vep
^\prime$ was constructed to vanish on $\cts \times \overline{\mathfrak
  a} $, so this is a contradiction.  Therefore $\vep = 0 $.
 $\square $  

\end{enumerate}
\end{pf}

\begin{rem} The proof of Proposition \ref{xy} in fact shows that any G-invariant Dirac structure on an irreducible Riemmanian symmetric space $G/K$ is a complex structure or a (complex) presymplectic structure.  
\end{rem} 

\begin{lem}If there is a G-invariant complex structure on $G/K$, then
  $\fg$ is not a complex Lie algebra.
\end{lem} 

\begin{pf} Suppose that $\fg$ is complex and that $G/K$ has a
  G-invariant complex structure $\mathfrak a \subset \mathfrak p _\cc$ as
  in part (2) of Lemma~\ref{xy}.  Then $\fg = \mathfrak u _\cc = \mathfrak u
  \oplus J \mathfrak u$ for some compact real form $\mathfrak u$ of
  $\fg$.  By the classification of irreducible semisimple
  orthogonal symmetric pairs ~\cite{hel}, $\fg$ is simple, $\mathfrak k = \mathfrak u$,
  and $\mathfrak p = J\mathfrak u$.  It is apparent that $\mathfrak a \subset
  \mathfrak p _\cc = J\mathfrak u \oplus iJ\mathfrak u \ J (\mathfrak
  u \oplus i \mathfrak u) = J \mathfrak u _\cc $.  This implies that
  $J\mathfrak a \subset \mathfrak u _\cc = \mathfrak u \oplus i
  \mathfrak u = \mathfrak k _\cc$.  It is
  also true that $[\mathfrak u _\cc , J\mathfrak a ] = J[\mathfrak u
  _\cc , \mathfrak a] = J[\mathfrak k _\cc , \mathfrak a ] =
  J\mathfrak a$ so that $J\mathfrak a$ is an ideal in $\mathfrak u
  _\cc \simeq \fg$.  This contradicts the fact that $\fg$ is
  simple.$\square$ 

\end{pf}

\subsubsection{Irreducible Riemannian Symmetric Spaces of the Non-Compact Type}

\noi  For this subsection, $(\fg  = \mathfrak k
\oplus \mathfrak p , \theta ) $ will be an irreducible semisimple
orthogonal symmetric pair of non-compact type associated to a symmetric space $G/K$.  When $(\fg , \fk)$ is of the non-compact type, $K$ is connected \cite{hel}, so $\fk$-invariance always implies $K$-invariance.

\begin{thm}\label{feb1a} Let $G/K$ be an irreducible Riemannian symmetric space of the non-compact type.  The following are equivalent:
\begin{enumerate}
\item $G/K$ admits a G-invariant generalized complex structure.
\item $G/K$ is a Hermitian symmetric space.  In particular, $G/K$ is Kahler with G-invariant Kahler structure.
\item $Z_{\mathfrak k} \neq 0 $.
\end{enumerate} 
\end{thm}


\noi Theorem~\ref{feb1a} will be proven in the course of several
lemmas.

\begin{lem}\label{feb1b} G-invariant generalized complex structures on $G/K$ exist
  only when $\mathfrak p _\cc$ is not an irreducible $\csts$-module.  In that case, there
  exists a G-invariant complex structure on $G/K$.

\end{lem}

\begin{pf}
First suppose that $Z_\mathfrak k \neq 0$.  Let $0 \neq z \in Z_{\mathfrak
  k}$.  $ad_z$ acts diagonally on $\mathfrak p _\cc$, so $\mathfrak p
  _\cc $ is a direct sum of eigenspaces for $z$.  Since $z$ commutes
  with $\mathfrak k$, $\mathfrak k$ preserves each eigenspace.  Thus,
  each eigenspace is a $\mathfrak k$-submodule of $\mathfrak p _\cc$.
  We have already seen that these are either all of $\mathfrak p _\cc$
  or complex structures $\mathfrak a \subset \mathfrak p _\cc$ as in
  Lemma~\ref{xy}.  If there exists such an $\mathfrak a$, then
  $\mathfrak p_\cc$ is not irreducible.  On the other hand, if $ad_z$
  has only one eigenvalue, $\lambda$ on 
  $\mathfrak p _\cc $, then $[\mathfrak p _\cc , \mathfrak p _\cc ]
  \subset \mathfrak k _\cc \cap (\mathfrak g _\cc )_{2\lambda} $,
  where $(\mathfrak g _\cc )_{2\lambda} $ is the $2\lambda$-eigenspace
  for $ad_z$.  However, $\mathfrak k _\cc \subset (\mathfrak g _\cc
  )_0 $.  This is possible if and only if $\lambda = 0 $ or
  $[\mathfrak p _\cc , \mathfrak p _\cc ] = 0 $.  If $\lambda = 0 $,
  then $z \in Z_\fg = 0$, which is a contradiction.  Therefore we must
  have $[\mathfrak p _\cc , \mathfrak p _\cc ] = 0 $, whence
  $[\mathfrak p , \mathfrak p ] = 0 $.\\

\noi Now supposing that $[\mathfrak p _\cc , \mathfrak p _\cc ] = 0 $,
the assumption that $Z_k \neq 0 $ and $\mathfrak p _\cc $ is
irreducible will lead to a contradiction.  
Since $(\fg , \theta ) $ is irreducible, either $\fg$ is simple, or
$\fg ^*$ (dual symmetric space, not the vector space dual) is simple.  
If $\fg$ is simple, then $[\mathfrak p , \mathfrak p ] = 0 $ implies
that $[\mathfrak k , \mathfrak k] \oplus \mathfrak p $ is a proper
ideal of $\fg$, which is a contradiction.  Therefore, whenever $Z_\mathfrak k \neq 0 $,
$\mathfrak p _\cc $ is not irreducible, so there exists a complex
structure $\mathfrak a \subset \mathfrak p _\cc$. If $\fg$ is not
simple, then $\fg ^*$ is simple, and a similar argument applies. \\

\noi Now suppose that $Z_\mathfrak k = 0 $.   There can be no
G-invariant symplectic structures because by Lemma~\ref{xy}, $\phi \in
Ann_{\mathfrak g _\cc ^*} ([\mathfrak k _\cc , \mathfrak k_\cc ]
\oplus \mathfrak p _\cc) = Ann_{\fg _\cc ^*} (\mathfrak g _\cc ^*) = 0
  $.  The only possible generalized complex structures are complex
  structures.  If a complex structure exists, as in Lemma~\ref{xy},
  then $\mathfrak p _\cc $ is not irreducible. $\square$
\end{pf}

\begin{lem}\label{feb1d} If $Z_\mathfrak k \neq 0 $, then $\mathfrak p _\cc$ is not
  irreducible and in fact $G/K$ has a G-invariant
  complex structure.
\end{lem}
\begin{pf} This is demonstrated in the proof of Lemma~\ref{feb1b}. $\square$
\end{pf}

\begin{lem}\label{feb1e} If there exists a G-invariant complex
  strcuture on $G/K$, then $G/K$ is a Hermitian symmetric space.  
\end{lem} 

\begin{pf} Let $\mathfrak a \subset \mathfrak p _\cc$ be a complex structure in the sense of Lemma~\ref{xy}.  We know that  $\mathfrak a \oplus \overline{\mathfrak a} =
  \mathfrak p _\cc$.  As usual, $\kappa$ will denote the Killing
  form on $\fg$, which when restricted to ${\mathfrak p \times \mathfrak p}$ is
  positive definite since $\fg$ is non-compact (It would be negative definite if 
$\fg$ were compact, in which case we could replace $\kappa$ by
  $-\kappa$).  $\kappa$ gives rise to its
  complexification $\kappa _\cc$ on $\mathfrak p _\cc \times \mathfrak
  p _\cc$, which will again be denoted by $\kappa$, due to the fact
  that $\kappa _\cc $ is none other than the Killing form on $\cts$.
  Defining $h (x,y) = \kappa (x , \overline y )$ gives a Hermitian
  form on $\mathfrak p _\cc $.  $h _{|\mathfrak a \times \mathfrak a}
  $ is still a Hermitian form.  But since $\mathfrak a $ is a complex
  structure, the projection $\pi : \cts \lra \fg $ provides an
  isomorphism $\pi : \mathfrak a \lra \fp _\cc $, thereby giving a Hermitian
  metric $H$ on $\fg / \fk _\cc = \fp _\cc$ (with respect to the complex structure $J$
  defined by $\mathfrak a$).  Then $H = g + i\omega$, where $g$ and
  $\omega$ are the real and imaginary parts of H, respectively.
  Because $H$ is hermitian, $g$ is positive definite and
  $J-invariant$, and because $\kappa$ is $Ad(K)$-invariant, so are $H$
  and $g$.  \\

\noi Finally, we observe that if $x \in \mathfrak p $ and $k \in \mathfrak k$,
then since $\mathfrak a $ is stable under $\csts$, $[k , x ] -i[k ,
  Jx] = [k , x -iJx ] = [k,x] -iJ[k,x]$.  It follows that $J$ commutes
with $ad(\mathfrak k)$.  This is all that is needed \cite{hel} for $G/K$ to be a
Hermitian symmetric space. $\square$
\end{pf}

\noi Theorem~\ref{feb1a} can now be proved:
\begin{pf} (1 $\Longrightarrow$ 2 ) Suppose that $G/K $ has a non-trivial generalized complex
  strucutres. Then by Lemma~\ref{feb1b}, $G/K$ admits a G-invariant complex structure.  Lemma~\ref{feb1e} ensures that  $G/K$ is Hermitian. The symplectic form $\omega$ associated with this Hermitian form is obviously invariant and is known to be symplectic ~\cite{hel} (i.e. $G/K$ is Kahler).  This is also easily checked by verifying that  since $\omega$ is G-invariant and vanishes on $\mathfrak k\times \fg$, then $d\omega = 0 $.\\


\noi (2 $\Longrightarrow$ 3 ) If $G/K$ is Hermitian, then there is a G-invariant symplectic structure $\phi \circ [ \; , \; ] $, with $\phi \in Z(\mathfrak k )^*$ as in Lemma~\ref{xy}.  Since this is symplectic, $\phi \neq 0$ and therefore $Z_{\mathfrak k} \neq 0 $.\\

\noi (3 $\Longrightarrow$ 1) If $Z_\mathfrak k \neq 0$, then Lemma~\ref{feb1d} implies that $G/K$ has a complex structure. $\square$
\end{pf} 

\begin{prop}\label{17junea}Let $(G,K)$ be an irreducible Riemannian symmetric pair of the non-compact type.  If $G/K$ admits any G-invariant generalized complex
  structures, then $\mathcal {CD} _{G/K} ^G = \cc \mathbb P ^1 \cup
  \{ two \; points \}$.  The two points are complex structures
  $\mathfrak a $, $\overline {\mathfrak a }$ as in Lemma~\ref{xy}, and $\mathcal{GC} _{G/K} ^G = \{ c \in \cc \; | \; Im(c) \neq 0 \} \cup
  \{ \mathfrak a , \overline {\mathfrak a} \} $.  The correspondence is established in the following way.  Fixing a G-invariant symplectic structure on $G/K$ $\omega$, $\cc ^\times \simeq \{L(T(G/K)_\cc , c.\omega ) \; | \; c \in \cc ^\times \}$, and  $T(G/K)_\cc$ corresponds to $ 0 \in \cc$, while $T^*(G/K)_\cc $ corresponds to $ \infty \in \mathbb{CP} ^1 $.  The complex structures $\mathfrak a$ and $\overline {\mathfrak a} $ are isolated points, whereas the new generalized complex structures are deformations of symplectic ones.
\end{prop}
\begin{pf} If $L = L(E,\vep ) $ is a complex Dirac structure, then $E$
  is a subalgebra of $\cts$ containing $\csts$.  By
  Lemma~\ref{feb19a}, $E = \csts$ or $E = \fk _\cc \oplus \mathfrak a$, $E = \fk _\cc  \oplus
  \overline{\mathfrak a} $, or $E = \cts$ as in Lemma~\ref{xy}. \\

\noi If $E = \cts$, then $\vep = \phi \circ [ \; , \; ]$ for some
$\phi \in Z_{\mathfrak k _\cc}^*$ as in Lemma~\ref{xy}.  By
Lemma~\ref{feb19a}, $Z_{\mathfrak k _\cc}$ is 1-dimensional, so If we
fix some $0 \neq \phi _0 \in   Z_{\mathfrak k }^*$, then all possible
pairs $(\cts , \vep )$ giving a complex Dirac structure are given by
$\vep = c\phi_0$ for some $c \in \cc$.  If $c = 0$, $L= T(G/K)$ and if
$c = \infty$, then $L = T^*(G/K)$ just as in
Proposition~\ref{m}.  By Lemma~\ref{xy}, the generalized complex
structures are given by all $c$ such that $ Im(c) \neq 0$.  \\

\noi If $E = \mathfrak k _\cc \oplus \mathfrak a$, then $\mathfrak a$
is a complex structure.  We will see that in this case, $\mathfrak a$
and $\overline { \mathfrak a }$ are the only complex structures.  Let
$\mathfrak b \subset \mathfrak p _\cc$ be another such complex
structure.  Again, we have $\mathfrak b \oplus \overline{\mathfrak b} = \mathfrak p _\cc$.  But
also $\mathfrak b \cap \mathfrak a = 0$ and $\mathfrak b \cap
\overline { \mathfrak a } = 0 $ since they are intersections
of distinct irreducible $\mathfrak k $-modules.  Consequently, $\mathfrak b$ is
the graph of some $\rr$-linear isomorphism $T: \mathfrak a \lra \overline
{ \mathfrak a }$.  For all $ k \in \mathfrak k$, $[k , x + Tx] = [k,x]
+[k,Tx] = [k,x] + T[k,x]$ because $\mathfrak b$ is a $\mathfrak
k$-module.  Hence, $T$ commutes with $ad_\mathfrak k$ and $T$ is in fact an isomorphism of $\mathfrak
k$-modules.  \\

\noi Since $G/K$ admits a G-invariant complex structure, $Z_\mathfrak
k \neq 0 $ by Theorem~\ref{feb1a}.  Let $0 \neq z \in Z_\mathfrak k$.
The proof of Lemma ~\ref{feb1b} shows that $z$ has exactly two
eigenvalues on $\mathfrak p _\cc$, $\lambda$ and $-\lambda$.  The
$\lambda $ eigenspace is $\mathfrak a$ and the $-\lambda$ eigenspace
is $\overline  { \mathfrak a }$.  For any $x\in \mathfrak a$, $-\lambda Tx
= [z,Tx] = T[z,x] = T(\lambda x) = \lambda Tx$.  Since $T$ is an
isomorphism, $\lambda = 0 $, which is impossible because $\fg$ has
trivial center (being semismiple).  Therefore, there exists no such
isomorphism $T$, and $\mathfrak a $ and $\overline { \mathfrak a }$ are
the only complex structures.  $\square$
\end{pf}

\subsubsection{Irreducible Riemannian Symmetric Spaces of the Compact Type} 

\noi Let $(G,K)$ be a Riemannian symmetric pair of the compact type.  Theorems \ref{feb1a} and \ref{17junea} are still true if $K$ is connected.  This is the case if $G$ is simply connected \cite{hel}.  If $K$ is not connected one must check which $\fa \subset \fp _\cc$ and $\om = \phi \circ [ \, , \, ]$ of the previous section are $K$-invariant.

\subsection{General Riemannian Symmetric Spaces}
\subsubsection{Semisimple Riemannian Spaces}\label{26nov7} 

It is known that any semisimple orthogonal symmetric Lie algebra
$(\fg , \theta ) $ is a direct sum of semisimple ideals $\fg _i $
preserved by $\theta$ such that $(\fg _i , \theta _i = \theta _{|\fg
  _i} )$ is itself an orthogonal symmetric pair.  Obviously, the
subpace $\mathfrak k$ fixed by $\theta$ is a direct sum of $\mathfrak
k_i \subset \fg _i $.

\begin{lem}\label{xz} Let $(\fg , \theta)$ be a semisimple orthogonal symmetric pair
  such that $\fg = \oplus \fg _i $
 is a direct sum of semisimple ideals $\fg _i$ such that $(\fg _i , \theta
 _i ) $ are all irreducible orthogonal 
symmetric pairs. Let $\csts \subset E \subset \cts$ be a subalgebra such
 that $E + \overline {E} = \cts$.  
Then $E = \oplus E_i$, where $E_i \subset (\fg _i)_\cc$.  
\end{lem} 
\begin{pf} Fix some $i$.  Let $E_i = pr_{(\fg _i)_\cc} E$. We wish to
 show that $E_i \subset E$.  Because each $(\fg _i)_\cc
$ is an ideal and closed under conjugation, $E_i$ is a subalgebra of
  $(\fg _i)_\cc$ and $E_i + \overline {E_i} = (\fg _i )_\cc$.  Due to Lemma~\ref{xy}
and the fact that $\fg _i$ is irreducible, $E_i = (\csts) _i \oplus
  \mathfrak a _i $, where either $\mathfrak a _i \oplus
  \overline{\mathfrak a _i} = (\mathfrak p _i)_\cc$ or $\mathfrak a _i
  =(\mathfrak p _i)_\cc$.  \\

\noi We know that  $\csts = \oplus (\csts) _i \subset E$, whence $(\csts)
_i \subset E $.  So to show that $E_i \subset E$, it suffices to show
that $\mathfrak a_i \subset E$.  But, again since each $\fg _j  $ is
an ideal, $[(\csts )_i , \mathfrak a _i ] \subset [(\csts) _i , E_i ] =
  [(\csts)_i, E] \subset E$. The only thing that needs to be checked is
  that $\mathfrak a _i = [(\csts) _i , \mathfrak a _i]$.
If $\mathfrak a _i = (\mathfrak p _i)_\cc$,
  then this is obviously true.  Otherwise, by Lemma~\ref{xy}, $\fg _i $ being
  irreducible means that $\mathfrak a _i $ is an irreducible
  $(\csts)_i$-module.  Consequently,  $\mathfrak a _i = [(\csts) _i ,
    \mathfrak a _i]$. $\square$
 \end{pf}

\begin{lem}\label{xa} Let $(\fg , \theta)$ be a semisimple orthogonal symmetric pair
  such that $\fg = \oplus \fg _i $
 is a direct sum of semisimple ideals $\fg _i$ such that $(\fg _i , \theta
 _i ) $ are all irreducible orthogonal 
symmetric pairs. Any GC-pair $(E , \vep) $ is a direct sum of GC-pairs
 $(E_i , \vep _i )$ for $\fg _i$ in the following sense:  $E = \oplus
 E_i $, each $(E_i , \vep _i )$ is a GC-pair for $(\fg _i , \theta _i
 ) $,and $\vep = \vep _1 \oplus \vep _2 ...\oplus \vep_n $.  
\end{lem} 

\begin{pf} Lemma~\ref{xz} shows that $E $ is a direct sum of the
 $E_i$'s.  It only remains to show that $\vep = \vep _1 \oplus \vep _2
 ...\oplus \vep_n $, which would follow if we could show that $\vep
 (\mathfrak a _i , \mathfrak a_j ) = 0$ if $ i \neq j $.  Here
 $E_i = (\csts) _i \oplus \mathfrak a _i $.  We see that
\[
 \vep (\mathfrak a _i , \mathfrak a_j) = \vep ([(\csts) _i , \mathfrak
 a_i ] , \mathfrak a _j) = \vep (\mathfrak a _i , [(\csts )_i ,
 \mathfrak a _j] ) = \vep (\mathfrak a _i , 0) = 0. 
\]   
\flushright { $\square$}
\end{pf}

\noi If $G/K$ is Riemannian symmetric space with G semisimple, then the
orthogonal symmetric pair $(\fg , \theta ) $ decomposes as
$\mathfrak g = \oplus \fg _i $, a direct sum of irreducible orthogonal
symmetric pairs $(\fg _i , \theta _i ) $.  But if G is simply connected,
this means that $G = \prod G_i $, where $G_i$ is the simply connected Lie
group with Lie algebra $\fg _i$.  Then $G/K \simeq \prod G_i / K_i $ as
long as $K$ is also connected.
Thus $G/K$ is a product of irreducible semisimiple Riemannian
symmetric spaces.  

\begin{thm}\label{xx}  Let $G$ be semisimple and $G/K$ a Riemannian symmetric
  space.  
\begin{enumerate}
\item If $G$ is simply connected and $K$ is connected, $G/K$ is then a product of irreducible Riemannian symmetric
  spaces $G_i /K_i $. Any G-invariant
  generalized complex structure on is a product of generalized complex
  structures on the $G_i/K_i$'s.  
\item Even if $G$ is not simply connected, a
  generalized complex structure on $G/K$ is given by a subalgebra $L =
  L(E, \vep ) \subset \cgts$, which may still be thought of as a
  product on the Lie algebra level:  $L = L( E_1 \oplus ...\oplus E_n  , \vep _1
  \oplus ...\oplus \vep n ) = \oplus L(E_i , \vep _i )$.
\end{enumerate} 
\end{thm}

\begin{pf} The proof is immediate from Lemmas ~\ref{xy}, ~\ref{xz},
  and ~\ref{xa}. $\square$
\end{pf}  

\subsubsection{When $G$ is the Isometry Group}

For the remainder of this section, $G$ will refer to the identity component of the isometry group of a Riemannian symmetric space.  Let $\fg = Lie(G)$. The
pair $(G,K)$ yields an orthogonal symmetric pair $( \fg , \theta )$. Then
$\fg$ decomposes into ideals $\fg = \fg _c \oplus \fg _n \oplus \fg _a
$, where $(\fg _c , \theta _c = \theta _{|\fg _c} ) $, $(\fg _n , \theta _n = \theta _{|\fg _n} ) $, and $(\fg _a , \theta _a = \theta _{|\fg _a} ) $ are
orthogonal symmetric pairs of the compact, noncompact, and abelian
type respectively. We have a decomposition $\fg _c = \mathfrak k _c
\oplus \mathfrak p _c $ and similar decompositions for $\fg _n$ and
$\fg _a$.  

\begin{lem}\label{xb} In notation described above, 
\begin{enumerate}
\item Let $E$ be a subalgebra of $\cts$ containing $\csts$ such that
  $E + \overline {E} = \cts$.  Then $E = E_c \oplus E_n \oplus E_a $
  as Lie algebras, where each summand is contained in $\fg _c $, $\fg
  _a$, or $\fg _n$.  
\item If $(E =   E_c \oplus E_n \oplus E_a , \vep)$ is a GC-pair, then
  $\vep = \vep _c \oplus \vep _n \oplus \vep _a $.  
\end{enumerate}
Therefore any GC-pair is of the form  $( E_c \oplus E_n \oplus E_a ,
\vep _c \oplus \vep _n \oplus \vep _a ) $.
\end{lem}

\begin{pf} $\;$
\begin{enumerate}
\item Since each summand is an ideal, $E_c = pr_{(\fg_c)_\cc} E$ is a
  subalgebra containing $(\mathfrak k _c )_\cc$ and such that $E_c +
  \overline {E_c} = (\fg_c)_\cc $.  The summand $\fg _c $ is of compact type, so
  it is semisimple.  We have seen in the proof of Lemma~\ref{xz} that 
$E_c = (\mathfrak k _c)_\cc  \oplus \mathfrak a _c$ and $\mathfrak a
  _c = [\mathfrak a _c ,  (\mathfrak k _c)_\cc ] \subset [ (\mathfrak
  k _c)_\cc , E] \subset [E,E] \subset E$.  Therefore $\mathfrak a _c
  \subset E$ and $E_c \subset E$.  An identical argument shows that
  $E_n \subset E$, whence also $E_a \subset E$.  Then $E$ is a direct
  sum of these.  

\item This would follow if we could show that $\vep ( \mathfrak a _c ,
  E_n \oplus E_a ) = \vep (\mathfrak a_n , E_c \oplus E_a ) = \vep
  (\mathfrak a _a , E_c \oplus E_n ) = 0 $.  We first argue that $\vep ( \mathfrak a _c ,
  E_n \oplus E_a ) = 0$.  Since $\fg _c $, $\fg
  _n$, and $\fg _a$ are ideals. $\vep ( \mathfrak a _c ,
  E_n \oplus E_a ) =  \vep ( [\mathfrak k _c , \mathfrak a _c] ,
  E_n \oplus E_a ) = \vep (\mathfrak a _c , [\mathfrak k _c , E_n
  \oplus E_a ]) = \vep (\mathfrak a _c ,0) = 0$.  An identical
  argument shows that $\vep (\mathfrak a_n , E_c \oplus E_a ) = 0$.
  Now it automatically follows that $\vep (\mathfrak a _a , E_c \oplus
  E_n ) = 0 $.  $\square$
\end{enumerate}
\end{pf}

\noi It is well known that any simply connected Riemannian symmetric
space $M$ can be expressed as a product $M = M_c \times M_n \times
M_a$ of Riemannian symmetric spaces of the compact, noncompact and
abelian types in accordance with the decomposition  $\fg = \fg _c \oplus \fg _n \oplus \fg _a
$.  

\begin{prop} Let $M = G /K $ be a simply connected Riemannian
  symmetric space, where $M = M_c \times M_n \times
M_a$ as above and where $G$ is the identity component of the isometry group.  Any equivariant generalized complex structure on $M$ is a product
of generalized complex structures on $M_c$, $M_n$, and $M_a$.  If,
however, $M$ is any Riemannian symmetric space, then a generalized complex structure $ L \subset \cgts$ may
still be thought of as a product since $L= L(E, \vep ) = L( E_c \oplus
E_n \oplus E_a , \vep _c \oplus \vep _n \oplus \vep _a ) = L(E_c ,
\vep _c ) \oplus L( E_n , \vep _n ) \oplus L(E_a  , \vep _a) $.  
\end{prop}
 
\begin{pf} This follows from Lemma~\ref{xb}. $\square$
\end{pf}

\noi Generalized complex structures on Riemannian symmetric spaces of
semisimple type (compact and noncompact) have already been described.
We now only need to describe generalized complex structures on
Riemannian symmetric spaces of abelian type.  \\

\noi $G/K$ is of the abelian, or Euclidean, type if the associated
orthogonal symmetric pair $(\fg = \mathfrak k \oplus \mathfrak p ,
\theta ) $ satisfies the condition that $\mathfrak p$ is an abelian
ideal of $\fg$. \\

\begin{lem}\label{20nov1} Let $G/K$ be of the Euclidean type, and let  $(\fg = \mathfrak k \oplus \mathfrak p ,
\theta ) $ be the associated orthogonal symmetric pair.  If $(E, \vep
)$ is a GC-pair, then $E = \mathfrak k _\cc \oplus \mathfrak a$, where
$\mathfrak a + \overline{\mathfrak a} = \mathfrak p _\cc $.  In fact, any such $\fa$ is of the form $\fa = \fa _1 \oplus \fa_2$, where $\fa _1 = \overline{ \fa _1}$ and $\overline {\fa _2 } \cap \fa = 0$.   
Furthermore, $Z^2 (E,\cc) \simeq(\wedge ^2 \mathfrak a ^*)^K$, where $Z^2 (E,
\cc) = \{ \vep \in \wedge ^2 E ^* \; | \; d_E \vep = 0 \} $ and
$(\wedge ^2 \mathfrak a ^*)^K$ is the space of K-fixed points in the
K-representation on $\wedge ^2 \mathfrak a ^*$.  
\end{lem} 
\begin{pf}
Obviously, any GC-pair $(E,\vep)$ must satisfy $E = \mathfrak k _\cc
\oplus \mathfrak a$, where $\mathfrak a +\overline{\mathfrak a} =
\mathfrak p _\cc$.  That $\fa = \fa _1 \oplus \fa _2 $ follows from the fact that since $\fk$ is a compactly embedded subalgebra, $\fa$ decomposes as a direct sum $\fa = \oplus _{i = 1}^n V_i$ of $\fk$-submodules.  Letting $I = \{ i \st \overline V_i \cap \fa \subset \fa \}$, $\fa _1 = \oplus _{i \in I} V_i$ and $\fa _2 = \oplus _{i \notin I} V_i$.  \\

\noi Now we address the question of which $\vep$ are
admissable.  To say that $\vep$ vanishes on $\fk$ is simply to say that $\vep \in \wedge ^2 \fa ^*$.  Using the formula for $d_E \vep$ given in Proposition~\ref{dhsprop2}, and checking $d_E \vep$ on each of $\fa \times \fa \times \fa$, $\fa \times \fa \times \fk$, $\fa \times \fk \times \fk$, and $\fk \times \fk \times \fk$, it is easy to see that $d_E \vep =0$ if and only if $\vep$ is $K$-invariant.  $\square$


\end{pf}

\subsection{Real Dirac Structures} 
All of the techniques used to describe generalized complex structures
carry over to the real case, and the results can be summarized in the
following two propositions.

\begin{prop} Let $G/K$ be a semisimple irreducible Riemannian
  symmetric space.  All G-invariant real Dirac structres on $G/K$ are
  presymplectic structures, i.e. are of the form $L(\fg ,
  \vep) $.  Any such $\vep $ is of the form $\vep = \phi \circ [ \; ,
  \; ]$, for some $\phi \in Ann_{\mathfrak k ^*} ( [\mathfrak k,
  \mathfrak k ]\oplus \mathfrak p ) \simeq Z(\mathfrak k ) ^* $.  
\end{prop}

\begin{prop} Let $G/K$ be a Riemannian symmetric space.  Any real
  Dirac structure is a product of Dirac structures in the sense
  described above.     
\end{prop} 



\begin{thebibliography}{10}




\bibitem{bou} Bourbaki, N.,  {\it Lie groups and Lie algebras. Chapters 4-6}, Springer-Verlag, 2002



\bibitem{com} Collingwood, D. and McGovern, W.M., {\it{Nilpotent orbits in semisimple Lie algebras}}, Van Nostrand Reinhold Mathematics Series, 1993

\bibitem{cou} Courant, T., Dirac Manifolds, {\it Tans. Amer. Math. Soc. }, 319 pp. 631-661, 1990

\bibitem{cow} Courant, T.; Weinstein, A.; Beyond Poisson Structures.  {\it Action hamiltoniennes de groupes. Troisieme theoreme de Lie }; vol. 27 of {\it Travaux en Cours}, pp. 39-49; Hermann, Paris, 1988
  


\bibitem{gua} Gualtieri, M., Generalized Complex Geometry (math.DG/0401221 v1)

\bibitem{gua2} Gualtieri, Generalized Complex Geometry (arxiv:math/0703298v2)

\bibitem{hit} Hitchin, N., Generalized Calabi-Yau manifolds (DG/0209099v1)

\bibitem{hel} Helgason, S., {{\it Differential Geometry, Lie Groups and Symmetric Spaces}}, Academic Press, 1978

\bibitem{hum} Humphreys, J.E. {\it Introduction to Lie Algebras and Representation Theory}, Springer-Verlag, 1972

\bibitem{kap}Kapustin, A., A-branes and Noncommutative Geometry (hep-th/0502212)

\bibitem{kna} Knapp, A. W., {\it{Lie Groups Beyond an Introduction, Second Edition} }, Birkhauser, 2002 





\bibitem{vai} Vaisman, I., {\it Lectures on the geometry of Poisson manifolds}, Birkhauser, 1994


\bibitem{wei} Weibel, C.A., {\it An introduction to homological algebra}, Cambridge University Press, 1994

\end{thebibliography}
\end{document}